\documentclass[12pt]{article}

\usepackage{amssymb, graphicx, parskip}

\usepackage{lipsum}
\usepackage[utf8]{inputenc}
\usepackage{titlesec}
\usepackage{enumitem}
\usepackage{amssymb,amsmath,amsthm}
\usepackage{graphicx}
\usepackage{todonotes}
\usepackage[utf8]{inputenc}
\usepackage[english]{babel}
\usepackage{float}
\usepackage{multibib}
\usepackage{appendix}
\usepackage[margin=2cm]{geometry}
\usepackage{hyperref}
\usepackage{ulem,cancel}

\makeatletter
\newcommand{\chapterauthor}[1]{%
  {\parindent0pt\vspace*{-25pt}%
  \linespread{1.1}\large\scshape#1%
  \par\nobreak\vspace*{35pt}}
  \@afterheading%
}
\makeatother

\numberwithin{equation}{section}
\newtheorem{theorem}{Theorem}[section]
\newtheorem*{theorem*}{Theorem}
\newtheorem{corollary}[theorem]{Corollary}
\newtheorem{lemma}[theorem]{Lemma}
\newtheorem{proposition}[theorem]{Proposition}
\newtheorem{definition}[theorem]{Definition}

\theoremstyle{definition}
\newtheorem{remark}[theorem]{Remark}
\newtheorem{example}[theorem]{Example}

\newcites{intro,urn,bp,ce,con}%
{References,
References,
References,
References,
References}

\newcommand{\PP}{\mathbb{P}}

\usepackage{mathtools}
\newcommand{\R}{\right}

\newcommand{\F}{\left}
\newcommand{\M}{\middle}

\newcommand{\Var}{\mathrm{Var}}
\newcommand{\E}{\mathbb{E}}
\newcommand{\Cov}{\mathrm{Cov}}

\newcommand{\bs}[1]{\boldsymbol{#1}}

\newcommand{\e}{\mathrm{e}}

\newcommand\eqindis{\stackrel{\mathclap{\mathrm{d}}}{=}}

\newcommand\conindis{\stackrel{\mathclap{\mathrm{d}}}{\rightarrow}}

\newcommand\coninprob{\stackrel{\mathclap{\mathrm{p}}}{\rightarrow}}

\title{Functional central limit theorems for non-local branching Markov processes}
\author{Christopher B. C. Dean\footnote{Department of Statistics, University of Warwick, Coventry, CV4 7AL, UK. Email: \texttt{\{Christopher.B.C.Dean\}, \{Emma.Horton\}@warwick.ac.uk}} \and Emma Horton$^*$}
\date{\today}

\begin{document}
\maketitle

\begin{abstract}
The aim of this paper is to study the fluctuations of a general class of supercritical branching Markov processes with non-local branching mechanisms. We establish functional central limit theorems and show that the limiting behaviour falls into three regimes, determined by the size of the spectral gap associated with the first-moment semigroup of the branching process. The main novelty is to develop a unified functional fluctuation theory for spatial branching Markov processes with non-local reproduction, allowing a general finite-dimensional spectral structure for the first-moment semigroup, including non-simple leading eigenvalues and nilpotent Jordan-type components. In doing so, we extend the classical small, critical and large fluctuation trichotomy beyond the finite-type and local spatial settings, and obtain limiting processes that capture the covariance structure induced by non-local offspring displacement.

\medskip


\medskip

\end{abstract}

\section{Introduction}\label{sec:intro}
Understanding the asymptotic behaviour of supercritical branching processes has long been a topic of great interest, see \cite{bAH1976,bathreya71,bEHK,bmoments,bSNTE2,bJanson,bCLT} and references therein. Recently, there has been an increased focus on results for branching processes with non-local branching mechanisms due to their relevance in a wide variety of applications, such as neutron transport \cite{bSNTE2, bSNTE1}, cell division \cite{telomeres, bstrongGF, bLLNGF} and the evolution of traits of individuals \cite{bmarguetLLN}. 

\medskip

The aim of this article is to provide a unified characterisation of the fluctuations of a general class of supercritical branching Markov processes with non-local branching mechanisms via functional central limit theorems. We show that there exist three regimes corresponding to the size of the spectral gap of the process, such that each regime results in different limiting behaviour. Results of this kind have been studied in a range of settings, including multi-type branching processes \cite{bathreya71,bJanson}, spatial branching Markov processes \cite{bCLTOU, bCLT,CLTNS}, CMJ processes \cite{CMJ} and superprocesses \cite{bsuperOU, bsuperOUren}. We also mention the recent work \cite{nicolas}, developed independently and in parallel to our own, which proves central limit theorems for non-local branching processes under weaker stability assumptions on the first-moment semigroup, but in a framework restricted to a simple leading eigenvalue and to convergence at fixed times in the small and critical regimes. The results presented in this article extend this existing literature in several directions.
\begin{itemize}
    \item We treat branching Markov processes in which particles move in a general spatial state space and, at branching times, may give birth to offspring at locations different from that of the parent. This simultaneously incorporates spatial motion and genuinely non-local reproduction. In this respect, our results go beyond the multi-type setting of \cite{bJanson}, where the state space is finite, and beyond the spatial branching frameworks of \cite{bCLT,CLTNS}, where branching is local in the spatial variable. The non-local mechanism is not a purely notational extension: it changes the structure of the second moment formulae and hence the covariance terms appearing in the fluctuation limits. In particular, even when the first moment semigroup admits a spectral description analogous to the local case, the second moment and covariance formulae contain new correlation terms arising from the joint law of the offspring locations; these terms disappear or collapse to simpler moments in the local setting.

    \item We allow the leading spectral structure of the first-moment semigroup to be finite-dimensional, rather than assuming a simple Perron eigenvalue. In particular, our framework permits several leading eigenfunctions and, more generally, nilpotent Jordan-type components. This leads to polynomial corrections in the law of large numbers and in the fluctuation scalings, and is not covered by the simple-eigenvalue assumptions in \cite{nicolas,bJanson,bCLT,CLTNS}.

    \item We obtain functional convergence, rather than only at a fixed time. More precisely, in the small and critical regimes (defined below) we prove convergence of finite-dimensional distributions for bounded test functions and functional convergence in the Skorokhod topology for a suitable regularity class of test functions. This strengthens the one-time central limit theorems available in \cite{nicolas,bCLT,CLTNS}. The functional formulation is important because it captures temporal correlations in the limiting fluctuations; in particular, in the small regime it reveals the dependence between the martingale-limit component and the remaining fluctuation component, while in the critical regime it makes visible oscillatory effects arising from complex eigenvalues that are hidden at a single time point.
\end{itemize}
Thus, the novelty of the paper lies not only in each of these extensions individually, but in treating them simultaneously within a single fluctuation theory covering the small, critical and large regimes.

Let us now introduce the setting in which we will be working. Let $E$ be a Polish space and $\dagger \notin E$ a cemetery (absorbing) state. We will write $B(E)$ for the set of complex-valued bounded measurable functions on $E$, $B^+(E)$ for those functions in $B(E)$ that are non-negative, and $B_1(E)$ (resp. $B_1^+(E)$) for those functions in $B(E)$ (resp. $B^+(E)$) that are uniformly bounded by unity. For a complex-valued function, $f$, we write $\bar f$ for its complex-conjugate. We will also sometimes restrict ourselves to the real-valued functions in these sets, in which case we will write $B_{\mathbb R}(E)$, and so on. For any $f \in B(E)$, we extend it to $E\cup \{\dagger\}$ by setting $f(\dagger) = 0$.

We consider a measure-valued stochastic process $X =(X_t, t\geq0)$ given by
\[
X_t := \sum_{i = 1}^{N_t}\delta_{x_i(t)}, \qquad t \ge 0,
\]
whose atoms $\{x_i(t) : i = 1, \dots, N_t\}$ evolve in $E\cup \{\dagger\}$ according to the following dynamics. From an initial position $x \in E$, particles evolve independently in $E$ according to a c\`adl\`ag Markov process $(\xi, \mathbf P_x)$. When at $y \in E$, particles branch at rate $\gamma(y)$, at which point the particle is replaced by a random number, $N$, of offspring at positions $x_1, \dots, x_N$ in $E$. {The law of the offspring locations and their number is denoted by ($\mathcal{Z},\mathcal P_y$), where $y$ denotes the position of the parent particle at the branch time and \[\mathcal{Z}=\sum_{i=1}^N\delta_{x_i}.\]}We will also often use the notation $\mathcal E_y$ for the associated expectation operator. We refer to $X$ as a $(\mathbf P, \gamma, \mathcal P)$-branching Markov process (BMP). We will let $\mathbb P_{\mu}$ denote its law when initiated from $\mu \in \mathcal M(E)$, the space of finite measures on $E$, with $\mathbb E_\mu$ denoting the corresponding expectation operator. 

\medskip

One way to characterise the behaviour of $X$ is via its linear expectation semigroup,
\begin{equation}\label{eq:linear}
  \psi_t[f](x) := \E_{\delta_x}\left[X_t[f]\right], \quad t \ge 0,\, x \in E, \, f \in B(E),
\end{equation}
where, for a measure $\mu$ and a function $f$, we have used the notation
\[
  \mu[f] := \int_E f(x)\mu({\rm d}x).
\]

Indeed, in many cases, the following assumption is satisfied.




\begin{enumerate}[label=(Ha),ref=(Ha)]
\item \label{Ha} There exists real number \(\lambda_1>0\), bounded function $\varphi_1 \in B^+(E)$, probability measure $\tilde\varphi_1$, such that $\tilde\varphi_{1}[\varphi_{1}] = 1$, and, for any $x \in E$, $t \ge 0$ and $f \in B(E)$,
\[
  \psi_t[\varphi_{1}](x) = {\rm e}^{\lambda_1 t}\varphi_{1}(x), \quad \text{and}\quad \tilde\varphi_{1}[\psi_t[f]] = {\rm e}^{\lambda_1 t}\tilde\varphi_{1}[f].
\]
Moreover, 
\begin{equation}
    \label{eq: ha assum}
  \sup_{x \in E, f \in B_1(E)}\big|\e^{-\lambda_1 t}\varphi_1(x)^{-1}\psi_t[f](x) - \tilde\varphi_1[f]\big|
  \to 0, \quad \text{as } t \to \infty.
\end{equation}

\end{enumerate}
 Assumption \ref{Ha} imposes a leading eigentriple for \(\psi\). In particular, for every \(t\geq 0\), \(\e^{\lambda_1 t}\) is an eigenvalue of \(\psi_t\) with eigenmeasure-eigenfunction pair \((\tilde\varphi_1,\varphi_1)\). Furthermore, this eigenvalue is leading, since when \(\psi_t\) is renormalised by the growth rate \(\e^{\lambda_1 t}\), we obtain convergence to a stationary profile given by \(\tilde\varphi_1\). This behaviour can be shown to hold for a large class of processes; we refer the reader to \cite{bNTEbook} and references therein for further details.

\medskip

In the case that $\lambda_1 > 0$, the process $X$ is called {\it supercritical}, since \ref{Ha} stipulates that the expected number of particles in the system grows exponentially in time.
In this regime, one can consider the almost sure behaviour of the system via a strong law of large numbers. Indeed, in \cite[Chapter 12]{bNTEbook}, it was shown that, for \(x \in E\) and \(f\in B_{\mathbb R}^+(E)\) such that \(f/\varphi \in B_{\mathbb R}^+(E)\),
\begin{equation}\label{eq: LLN cor}
\frac{{\rm e}^{-\lambda_1 t}X_t[f]}{\varphi_1(x)} \to \tilde\varphi_1[f]W_\infty, \quad \text{as } t \to \infty,
\end{equation}
almost surely and in $L^2$, where $W_\infty$ is the almost sure limit of the unit mean $\mathbb P_{\delta_x}$-martingale 
\begin{equation}\label{eq:mg cor}
\frac{{\rm e}^{-\lambda_1 t}X_t[\varphi_1]}{\varphi_1(x)}, \quad t \ge 0, \, x \in E.
\end{equation}

The aim of this article is to further characterise the fluctuations of $X$ in the supercritical regime. In the setting of multitype and spatial branching processes, it is known (see \cite{bJanson}, \cite{CLTNS}) that the fluctuations obey one of three possible dynamics depending on the size of the spectral gap of \(\psi\). As we will see in this article, this trichotomy remains true in the non-local setting. We thus require additional assumptions on the spectrum of $\psi$. In the rest of this section, we will present simplified versions of the much more general results stated in the next section, since the notation becomes heavier, although the moral of the results remains the same. We thus introduce a (simplified) assumption that gives the existence of a second eigenvalue of $\psi$. 

\begin{enumerate}[label=(Hb),ref=(Hb)]
\item \label{Hb} There exists a real number \(\lambda_2<\lambda_1\), a bounded function $\varphi_2 \in B_{\mathbb{R}}(E)$, a probability measure $\tilde\varphi_2$ such that $\tilde\varphi_{2}[\varphi_{2}] = 1$, $\tilde\varphi_{i}[\varphi_{j}] = 0$, \(i\neq j\), and, for any $x \in E$, $t \ge 0$ and $f \in B(E)$,
\[
  \psi_t[\varphi_{2}](x) = {\rm e}^{\lambda_2 t}\varphi_{2}(x), \quad \text{and}\quad \tilde\varphi_{2}[\psi_t[f]] = {\rm e}^{\lambda_2 t}\tilde\varphi_{2}[f].
\]
Moreover, 
\begin{equation}
    \label{eq: hb assum}
  {\sup_{x \in E, f \in B_1(E)}t\bigg|\e^{-\lambda_2 t}\psi_t[f](x) - \e^{(\lambda_1-\lambda_2)t}\tilde\varphi_1[f]\varphi_1(x)-\tilde\varphi_2[f]\varphi_2(x)\bigg|
  \to 0, \quad \text{as } t \to \infty.}
\end{equation}
\end{enumerate}

As stated previously, there are three possible dynamics for the leading-order fluctuations depending on the size of the spectral gap of \(\psi\). We say that the process is in the {\it large} regime when $2\lambda_2 - \lambda_1 > 0$, the {\it critical} regime when $2\lambda_2 - \lambda_1 = 0$ and the {\it small} regime when $2\lambda_2 - \lambda_1 < 0$. Under \ref{Ha}, \ref{Hb} and some other assumptions, which we detail in the next section, we obtain the following results. Here and throughout the article, we write $\mathbb D[0, \infty)$ for the space of c\`adl\`ag paths on $[0, \infty)$ under the $J_1$ Skorokhod topology. {In what follows, \(C_{\psi}(E)\) is a sufficiently regular function space to be defined later, moreover \(C_{\psi}(E)\) can be relaxed to \(B(E)\) if one is interested in a weaker form of convergence (in probability and in \(L^2\) for the large regime, and in the sense of finite dimensional distributions for the small and critical regimes).}

\begin{theorem*}[Large Regime (simplified)]
\label{thm: main result large cor}
 Suppose that \(2\lambda_2>\lambda_1\) and let \(f \in C_{\psi}(E)\). Then, for any initial position \(x \in E\), as \(t\rightarrow \infty\)
\begin{equation}
\label{eq: large result new 1 cor}
   \e^{-\lambda_2t}\F(X_{t}[f]-\e^{\lambda_1 t}\tilde\varphi_{1}[f]\varphi_1(x)W_{\infty} \R)\to \tilde\varphi_{2}[f]\widetilde W_{\infty}(x),
\end{equation}
\(\mathbb{P}_{\delta_x}\)-almost surely and in \(L^2\), where \(\widetilde W_{\infty}(x)\) is the \(\mathbb{P}_{\delta_x}\)-almost sure limit of the martingale
    \begin{equation}
        \widetilde M_{t} :=  \e^{-\lambda_2 t}X_t[\varphi_{2}], \quad t\geq 0.
    \end{equation}
\end{theorem*}

\begin{theorem*}[Small Regime (simplified)]
\label{thm: main result small cor}
 Suppose that \(2\lambda_2<\lambda_1\) and let \(f \in C_{\psi}(E)\). Then, for any initial position \(x \in E\), as $n \to \infty$
    \begin{equation}
    \label{eq: main result small cor}
        \e^{-\frac{\lambda_1}{2} (n+t)}\F(X_{n+t}[f]-\e^{\lambda_1 (n+t)}\tilde\varphi_{1}[f]\varphi_1(x)W_{\infty}  \R)\conindis (\varphi_1(x)W_{\infty})^{1/2}Z_S(t), \quad in \quad \mathbb D[0,\infty), 
    \end{equation}
  
    where $Z_S$ is a mean zero Gaussian process that is independent of $W_\infty$.
\end{theorem*}

\begin{theorem*}[Critical Regime (simplified)]
       \label{thm: main result critical cor}
       Suppose that \(2\lambda_2=\lambda_1\) and let $f \in C_{\psi}(E)$. Then, for any initial position \(x \in E\), as $n \to \infty$,
    \begin{equation}
    \label{eq: main result crit simplified}
        n^{-\frac{1}{2}}\e^{-\frac{\lambda_1nt}{2} }(X_{nt}[f]-\e^{\lambda_1 nt}\tilde\varphi_1[f]\varphi_1(x)W_{\infty})\conindis (\varphi_1(x)W_{\infty})^{1/2}Z_C(t), \quad in \quad \mathbb D[0,\infty),
    \end{equation}
    where \(Z_C\) is a scaled Brownian motion that is independent of \(W_{\infty}\).
\end{theorem*}

The classification into these three regimes can be stated in terms of the spectral gap, $\lambda_1 - \lambda_2$. Indeed, the large (resp.\ critical, small) regime corresponds to the spectral gap being less than (resp.\ equal to, greater than) $\lambda_1/2$. Now, regardless of the regime, the leading order behaviour of the system will always be described by $(\lambda_1, \varphi_1, \tilde\varphi_1)$, as in \ref{Ha}. However, the size of the spectral gap dictates how quickly this convergence happens, and hence how the system fluctuates around this first order behaviour. 

In the small regime, the spectral gap is greater than $\lambda_1/2$, and so the convergence seen in \ref{Ha} happens quickly. This gives too little time for second-order fluctuations to manifest. Thus, in the small regime, {\color{black}we see a limit over a small time scale that is stationary}. In the critical regime, the speed of convergence, $\lambda_1/2$, is such that the ``pull'' towards the stationary distribution in \ref{Ha} matches the ``push'' away arising from the second-order fluctuations, thus resulting in Brownian fluctuations over a longer timescale. Finally, if the spectral gap is smaller than $\lambda_1/2$, then the speed of convergence in \ref{Ha} is slow enough to allow the growth of second-order fluctuations. This dominates the stationary behaviour seen in the small and critical regimes, yielding fluctuations that mimic the behaviour of \eqref{eq: LLN cor}. As previously noted, we will present these results in a more general framework that removes the need for simple eigenvalues $\lambda_1$ and $\lambda_2$. However, this heuristic remains true in the more general setting. 

{The rest of the article is set out as follows. In the next section, we give our main results in full generality, along with the necessary assumptions. 
In Section \ref{sec:discussion}, we will spend some time discussing our assumptions and relevant examples. 
Section \ref{sec:proofs2} contains the proofs of the main results. Finally, we provide three appendices containing some useful results that will be employed throughout the article, along with proofs of some technical results. }

\section{Main results}\label{sec: main results general}
In this section, we state the more general versions of the assumptions \ref{Ha} and \ref{Hb} that will subsequently allow us to present our main results, which generalise the simplified versions stated in the introduction.

\subsubsection*{Spectral assumption}

Before we state our main assumption, we require some preliminary definitions. 
\begin{definition}
\label{definition: eigenspace}
For \(\lambda\in\mathbb C\), we say that \(\lambda\) is an eigenvalue of \(\psi\) with finite-dimensional generalised eigenspace \(A_\lambda\subset B(E)\) if there exists a bounded projection
\[
\Pi_\lambda:B(E)\to A_\lambda
\]
such that, writing \(K_\lambda:=\ker \Pi_\lambda\), we have
\begin{equation}\label{eq: direct sum eigenprop}
B(E)=A_\lambda\oplus K_\lambda,
\end{equation}
both \(A_\lambda\) and \(K_\lambda\) are invariant under \(\psi_t\) for every \(t\geq 0\), and
\begin{equation}\label{eq: eigenvalue prop equivalent}
\psi_t|_{A_\lambda}=e^{(\lambda+\mathcal N_\lambda)t},
\end{equation}
where \(\mathcal N_\lambda:A_\lambda\to A_\lambda\) is nilpotent, ie. there exists $p > 0$ such that $\mathcal N_\lambda^p = 0$.
\end{definition}

The above definition entails the existence of a basis \(\varphi_{j}^{(k)}\), \(1\leq j \leq p\), \(1\leq k\leq k_j\), of \(A_{\lambda}\), such that 
\begin{equation*}
    \mathcal{N}_{\lambda}\varphi_{1}^{(k)}=0,\quad 1\leq k \leq k_1, \quad \mathcal{N}_{\lambda}\varphi_{j}^{(k)}=\varphi_{j-1}^{(k)},\quad 2\leq j \leq p,\quad 1\leq k \leq k_j.
\end{equation*}

One can see this as the Jordan basis of the generator of \(\psi_t|_{A_{\lambda}}\); indeed, we shall refer to \(\varphi_{j}^{(k)}\) as an eigenfunction of rank \(j\). Moreover, since \(\Pi_{\lambda}\) is bounded, there exist bounded (with respect to \(\|\cdot \|_{\infty}\)) linear coordinate maps \(\tilde\varphi_{j}^{(k)}:B(E)\rightarrow \mathbb{C}\), such that 
\begin{equation}\label{eq:eigenprops}
\tilde\varphi_{j}^{(k)}[\varphi_{j}^{(k)}]=1, \qquad \tilde\varphi_{j}^{(k)}[\varphi_{\ell}^{(q)}]=0, \, (j,k)\neq (\ell ,q), \qquad \tilde\varphi_{j}^{(k)}[f]=0, \, f \in K_{\lambda}.
\end{equation}

Note that for \(f \in B(E)\), \(x \in E\), \(\overline{\psi_t[f](x)}=\psi_t[\bar{f}](x)\). Thus, if \(\lambda \in \mathbb{C}\setminus {\mathbb{R}}\), by taking the complex conjugate in  Definition \ref{definition: eigenspace} we see that there exists a conjugate eigenspace \(A_{\bar\lambda}=\mathrm{span}\{\bar \varphi_{j}^{(k)},1\leq j \leq p,1\leq k\leq k_j\}\) corresponding to the eigenvalue \(\bar \lambda\) and nilpotent operator \(\mathcal{N}_{\bar \lambda}\) satisfying
\begin{equation}
\label{eq: complex eigenspace}
    \mathcal{N}_{\bar \lambda}\bar\varphi_{1}^{(k)}=0, \; 1\leq k \leq k_1, \qquad \mathcal{N}_{\bar \lambda}\bar\varphi_{j}^{(k)}=\bar \varphi_{j-1}^{(k)},\; 2\leq j \leq p, \, 1\leq k \leq k_j.
\end{equation}

The operator \(\mathcal N_\lambda\) should be understood as the nilpotent part of the
semigroup restricted to the generalised eigenspace \(A_\lambda\). Since \(\mathcal N_\lambda\) is nilpotent, \(e^{t\mathcal N_\lambda}\) is a finite polynomial in
\(t\). The operator \(\mathcal N_\lambda\) therefore records the Jordan block structure
associated with the eigenvalue \(\lambda\): if \(\mathcal N_\lambda=0\), then
\(A_\lambda\) is spanned by genuine eigenfunctions, while if
\(\mathcal N_\lambda\neq 0\), then \(A_\lambda\) also contains generalised
eigenfunctions, and polynomial prefactors appear in the asymptotics.

For notational clarity, when several eigenvalues, \(\lambda_1,\dots,\lambda_m\in\mathbb{C}\), are considered at once, we
write the Jordan basis of \(A_{\lambda_i}\) as
\[
   \varphi^{(k)}_{i,j},
   \qquad 1\leq j\leq p_i,\quad 1\leq k\leq k_{i,j},
\]
where \(i\) indexes the eigenvalue, \(j\) indexes the rank in the Jordan chain, and \(k\) indexes multiplicity within that rank. Thus
\[
   \mathcal N_{\lambda_i}\varphi^{(k)}_{i,1}=0,
   \qquad
   \mathcal N_{\lambda_i}\varphi^{(k)}_{i,j}
   =
   \varphi^{(k)}_{i,j-1},
   \quad j\geq 2.
\]
We also have \(\tilde\varphi_{i,j}^{(k)}[\varphi_{\ell,n}^{(q)}]=0\) for each \((i,j,k)\neq (\ell,n,q)\).

Moreover, when several eigenspaces are considered simultaneously, we write \(\mathcal N\) for the
operator obtained by applying the corresponding nilpotent part on each
generalised eigenspace and killing the remaining component. More precisely, if
\[
        f = \sum_{i=1}^m f_{\lambda_i} + f_0,
        \qquad
        f_{\lambda_i}\in A_{\lambda_i}, \quad
        f_0\in \bigcap_{i=1}^m K_{\lambda_i},
\]
then we define
\[
        \mathcal Nf := \sum_{i=1}^m \mathcal N_{\lambda_i} f_{\lambda_i}.
\]
Thus \(\mathcal N\) is block diagonal across the generalised eigenspaces
\(A_{\lambda_1},\ldots,A_{\lambda_m}\), with block \(\mathcal N_{\lambda_i}\) on
\(A_{\lambda_i}\), and \(\mathcal Nf_0=0\). In particular,
\begin{equation}\label{eq:gen-eigenfunctions}
  \psi_t[\varphi_{i,j}^{(k)}](x) = {\rm e}^{(\lambda_i+\mathcal{N}) t}\varphi_{i,j}^{(k)}(x), \quad \qquad \tilde\varphi_{i,j}^{(k)}[\psi_t[f]] = {\rm e}^{\lambda_i t}\tilde\varphi_{i,j}^{(k)}[{\rm e}^{\mathcal{N}t}f], \quad f\in B(E).
\end{equation}
We can now state our main assumption on the spectrum of \(\psi\).

\paragraph*{Spectral assumption.}
{\begin{enumerate}[label=(H1),ref=(H1)]
\item \label{H1b} There exists an integer \(m\geq 1\), a real number \(\lambda_1>0\), \(\lambda_{2},\dots ,\lambda_{m}\in \mathbb{C}\), that satisfy\footnote{If $m = 1$, the list $\lambda_2 \dots, \lambda_m$ is an empty list.} \(\lambda_1>\mathrm{Re}\lambda_2\geq \dots \geq \mathrm{Re}\lambda_m\), such that \(\lambda_1,\dots,\lambda_m\) are eigenvalues of \(\psi\) in the sense of Definition \ref{definition: eigenspace} such that
\begin{equation}
    \label{eq: h1b assum}
  \sup_{x \in E, f \in B_1(E)}\e^{-\mathrm{Re}\lambda_m t}t\bigg|\psi_t[f](x) - \sum_{i=1}^{m}\e^{(\lambda_i+\mathcal{N})t}\Phi_{i}[f](x)\bigg|
  \to 0, \quad \text{as } t \to \infty,
\end{equation}
where 
\[\Phi_{i}[f](x)=\sum_{j=1}^{p_i}\Phi_{i,j}[f](x) \quad \text{ with } \quad
\Phi_{i,j}[f](x) = \sum_{k=1}^{k_{i,j}}\tilde\varphi_{i,j}^{(k)}[f]\varphi_{i,j}^{(k)}(x). \quad \]
\end{enumerate}}


We reiterate that here \(i\) indexes the eigenvalue, \(j\) indexes the rank in the Jordan chain, and \(k\) indexes multiplicity within that rank. The operator \(\mathcal N\) is the global nilpotent operator obtained by collecting the
nilpotent shifts \(\mathcal N_{\lambda_i}\) along each Jordan chain; on
\(A_{\lambda_i}\) it coincides with \(\mathcal N_{\lambda_i}\), and on the residual
component it is zero.

{Heuristically, \ref{H1b} implies that the spectrum of the semigroup \(\psi\) is dominated by \(m\geq 1\) eigenvalues whose eigenspaces are finite dimensional, with larger values of $m$ yield a better approximation of $\psi$ as a sum over the first $m$ eigenspaces. }

Let us illustrate \ref{H1b} further by considering some specific cases in more detail. First note that when \(m=p_1=k_{1,1}=1\), \ref{H1b} is the natural spectral analogue of \ref{Ha}. In the case that $m = p_1 = 1$ but $k_{1, 1} > 1$, we extend \ref{Ha} to the case where $\lambda_1$ is not simple. This case implies that $k_{1, 1}$ is both the algebraic and geometric multiplicity of $\lambda_1$, with corresponding eigenfunctions given by $\varphi^{(1)}, \dots, \varphi^{(k_{1, 1})}$ (note, we have dropped the subscript for simplicity). In the case that $p_1= 2$, say, the functions $\varphi_{1, 1}^{(1)}, \dots, \varphi_{1, 1}^{(k_{1, 1})}$ are eigenfunctions, but the functions $\varphi_{1, 2}^{(1)}, \dots, \varphi_{1, 2}^{(k_{1, 2})}$ are {generalised eigenfunctions of rank $2$}, meaning $\mathcal N\varphi_{1, 2}^{(k)} \neq 0$ but $\mathcal N^2\varphi_{1, 2}^{(k)} = 0$. In general, if $p_1 \ge 3$ then $\varphi_{1, j}^{(k)}$ is of rank $j$. Note that, in the case that $m=1$ and $p_1 \ge 2$, the geometric multiplicity of $\lambda_1$ is $k_{1, 1}$ and the algebraic multiplicity is $\sum_{j = 1}^{p_1} k_{1, j}$. 
Next, setting $m = 2$ with \(\lambda_1,\lambda_2 \in \mathbb{R}\) and \(p_1 =p_2 = k_{1,1} = k_{2,1}=1\) recovers both \ref{Ha} and \ref{Hb}. {Allowing $m \ge 3$ generalises this to $m$ eigenvalues \(\e^{\lambda_1 t},\dots,\e^{\lambda_m t}\), each with eigenfunctions $\varphi_{i, 1}^{(1)}, \dots, \varphi_{i, 1}^{(k_{i, 1})}$ and generalised eigenfunctions $\varphi_{i, j}^{(1)}, \dots, \varphi_{i, j}^{(k_{i, j})}$ of rank $j$.} The asymptotic given in \eqref{eq: h1b assum} stipulates that, when normalised by the smallest growth rate ${\rm e}^{\mathrm{Re}\lambda_m t}$, the expectation semigroup converges. 

 Recall that each eigenvalue in \(\mathbb{C}\setminus \mathbb{R}\) has an associated conjugate eigenvalue. Thus, since \ref{H1b} captures the discrete part of the spectrum of \(\psi\) whose eigenvalues have real part at least \(\rm{Re}\lambda_m\), and since
\[
\rm{Re}\overline{\lambda_i}=\rm{Re}\lambda_i,
\]
for each \(1\leq i \leq m\), there exists \(1\leq \ell\leq m\) such that \(\overline{\lambda_i}=\lambda_\ell\).

Throughout the paper we choose the Jordan bases compatibly with this conjugation. Thus, for every \(1\leq j\leq p_i\) and \(1\leq k\leq k_{i,j}\), there exist unique \(1\leq n\leq p_\ell\) and \(1\leq q\leq k_{\ell,n}\) such that
\[
    \overline{\varphi_{i,j}^{(k)}}=\varphi_{\ell,n}^{(q)}.
\]
We introduce some further notation. Define
\begin{align*}
    &m_L:=\max\{1\le i \le m: 2\mathrm{Re}\lambda_i>\lambda_1\},\\
    &m_C:=\max\{1\le i \le m: 2\mathrm{Re}\lambda_i\geq \lambda_1\}.
\end{align*}

 {Since the eigenvalues are ordered by their real part, the index $m_L$ (resp.\ $m_C$) represents the last large\footnote{i.e. satisfying $2{\rm Re}\lambda_i > \lambda_1$ (resp. $2{\rm Re}\lambda_i \ge \lambda_1$)} (resp.\ critical) eigenvalue. The case \(m_L=1\) corresponds to the situation in which no non-leading eigenvalue is large. The setting of $m_L=m_C$ corresponds to the existence of no critical eigenvalues.}


\subsubsection*{Branching and moment assumptions.}
We introduce the following boundedness assumptions on the branching rate and offspring distribution.
\begin{enumerate}[label=(H\arabic*), ref=(H\arabic*)]
\setcounter{enumi}{1} 
\item\label{H4} The branching rate \(\gamma \in B_{\mathbb R}^+(E)\). 
\item[(M\(k\))]\label{H2} $\sup_{x \in E}\mathcal{E}_x[N^k]< \infty$.
\end{enumerate}

Assumption \ref{H4} is a convenient sufficient condition ensuring that the branching process does not explode in finite time. It will also allow us to obtain the uniform moment and covariance estimates, in \(x\in E\) and \(f\in B_1(E)\). We expect that \ref{H4} could be replaced by suitable non-explosion and integrability assumptions adapted to the underlying motion and offspring law, but we do not pursue this extension here.

Regarding the moment assumptions, the functional convergence results below are stated under \((M4)\), while some \(L^k\) and finite-dimensional statements require only the corresponding moment assumptions specified in their theorems. 

\subsubsection*{Regularity assumptions.}
We also require a class of functions which we refer to as $C_{\psi}(E)$. We say that \(f\in C_{\psi}(E)\subseteq B(E)\) if there exists an even integer \(k_f \geq 2\), and {\(C_f\geq 0\)}, such that \begin{align}
                 & \sup_{x \in E, 0 \leq t \leq {\color{black}1}}t^{-1/k_f}\ |\mathbf E_{x}[f(\xi_{t})]-f(x)| \leq C_f, \label{eq: tight assum 1112}\\
                &\sup_{x \in E,0\leq t \leq {\color{black}1}} t^{-1}\mathbf E_x[(f(\xi_t)-\mathbf E_x[f(\xi_t)])^{k_f}]\leq C_f,  \label{eq: tight assum 2112}\\
                &\sup_{x \in E,0\leq t \leq {\color{black}1}} t^{-1}\mathbf E_x[(f(\xi_t)-\mathbf E_x[f(\xi_t)])^{2k_f}]\leq C_f, \label{eq: tight assum 3112} 
                \end{align}
            and,
            \begin{equation*}
                \sup_{x \in E}\mathcal{E}_x[N^{2k_f}]< \infty.
            \end{equation*}

The role of \(C_\psi(E)\) is solely to obtain functional convergence. The finite-dimensional versions of the small and critical regime results hold for all bounded measurable test functions \(f\in B(E)\), whereas the functional convergence results require \(f\in C_\psi(E)\) in order to control short-time increments of \(X_t[f]\) uniformly in the initial position.            

The size of the class $C_\psi(E)$ depends on the regularity of the underlying particle motion and the available moment bounds on the offspring distribution. Note that, when only \((M4)\) is assumed, \(C_\psi(E)\) should be understood as the class of functions satisfying \eqref{eq: tight assum 1112}-\eqref{eq: tight assum 3112} with \(k_f=2\). More generally, if a larger value of \(k_f\) is needed to satisfy \eqref{eq: tight assum 1112}-\eqref{eq: tight assum 3112}, then stronger offspring moment assumptions are required. Again, we discuss this assumption further in Section \ref{sec:discussion}.

\subsubsection*{Results}

We now present our main results. We assume that \ref{H1b} and \ref{H4} hold throughout.

\medskip

\begin{theorem}[Large Regime]
\label{thm: main result large}
 Assume that (M\(4\)) holds and let \(f \in C_{\psi}(E)\). 
 Then, for any initial position \(x \in E\), as \(t\rightarrow \infty\)
\begin{equation}
\label{eq: large result new 1}
   \e^{-\mathrm{Re}\lambda_{m_L}t}\F(X_{t}[f]-\sum_{i=1}^{m_L}\sum_{j=1}^{p_i}\sum_{k=1}^{k_{i,j}}\tilde\varphi_{i,j}^{(k)}[\e^{(\lambda_i+\mathcal{N}) t}f]W_{i,j}^{(k)}(x) \R)\to 0, 
\end{equation}
\(\mathbb{P}_{\delta_x}\)-almost surely 
and in \(L^4\), where, for \(1\leq i \leq m_L\), \(1\leq j \leq p_i\), and \(1\leq k \leq k_{i,j}\), \(
W_{i,j}^{(k)}(x)\) is the \(\mathbb{P}_{\delta_x}\)-almost sure limit of the martingale
    \begin{equation}
    \label{eq: general martingales}
        M_{i,j,t}^{(k)} :=  X_t[\e^{-(\lambda_i+\mathcal{N}) t}\varphi_{i,j}^{(k)}], \quad t\geq0,
    \end{equation}
    which is indeed a martingale by \eqref{eq:gen-eigenfunctions}.
\end{theorem}

\begin{theorem}[Large regime $L^k$ convergence]\label{cor:large}
For \(k\geq2\) even, suppose that (M$k$) holds and take $f \in B(E)$. Then \eqref{eq: large result new 1} holds with convergence in $L^k$.
\end{theorem}

\medskip

As a corollary of Theorem \ref{thm: main result large}, we obtain the following strong law of large numbers. Note that this extends \cite[Theorem 12.3]{bNTEbook} to non-simple leading eigenvalues.
\begin{corollary}[LLN]
\label{thm: SLLNs}
    Let \(f \in C_{\psi}(E)\). Under the conditions of Theorem \ref{thm: main result large}, for any initial position \(x \in E\), we have that as $t \to \infty$,
    \begin{equation}
        \e^{-\lambda_1 t}t^{-(p_1-1)}X_t[f]\rightarrow \frac{1}{(p_1-1)!}\sum_{i=1}^{k_{1,1}}\tilde\varphi_{1,1}^{(i)}[\mathcal{N}^{p_1-1}f]W_{1,1}^{(i)}(x), \label{eq: sllns}
    \end{equation}
    \(\mathbb{P}_{\delta_x}\)-almost surely and in \(L^4\). Moreover, for \(f \in B(E)\), under the assumptions of Theorem \ref{cor:large}, \eqref{eq: sllns} holds in \(L^k\).
\end{corollary}

\medskip

\begin{remark}\label{remark : large regime}~
\begin{itemize}
    \item Note that we do not state a separate functional convergence theorem in the large regime since Theorem \ref{thm: main result large} already gives a stronger leading order description: after subtracting the contributions of the large spectral modes, the normalised remainder converges to zero almost surely and in \(L^4\). Thus the dominant large regime fluctuations are captured by the martingale limits \(W^{(k)}_{i,j}(x)\), rather than by a new Gaussian fluctuation process in \(\mathbb D[0,\infty)\). A non-trivial functional result in this regime would concern only smaller order residual fluctuations after this martingale limit expansion has been removed.
    \item The main novelty of these results is that it extends the standard law of large numbers type results, see e.g. \cite[Theorem 12.3]{bNTEbook}, to the case where the spectral assumption \ref{H1b} holds, where the leading large modes may form a finite-dimensional generalised eigenspace with several eigenfunctions and nilpotent Jordan components. Consequently, the large regime approximation is expressed in terms of the family of martingale limits \(W^{(k)}_{i,j}\) and the polynomial corrections generated by \(\mathcal N\). 
    \item Comparing Theorems \ref{thm: main result large} and \ref{thm: main result large cor}, we see that \(f\in C_\psi(E)\) and (M4) are only needed for the pathwise/tightness component of Theorem \ref{thm: main result large}; the \(L^k\) convergence holds for all bounded test functions under the corresponding moment assumption. We also note that in the latter case the joint $k$-th moments of the martingale limits, \(W_{i,j}^{(k)}(x)\), are given by Theorem 2.1 of \cite{Moments_Chris}. In particular, (M\(2\)) implies that the limits are non-degenerate.
\end{itemize}
\end{remark}

\medskip

We now move to the small and critical regimes. In these regimes we will state two results, one detailing functional convergence and another for convergence in finite dimensional distributions. For the former, we will obtain Gaussian limits, which will be fully characterised by their covariance structure. For this, we need the following additional notation.

We first introduce a covariance-type operator for the offspring distribution. For \(f,g \in B(E)\), \(t\geq 0\), $x \in E$, define
\begin{equation}\label{eq:V}
    V[f,g](x) = \gamma(x)\mathcal{E}_{x}\Bigg[\sum_{\substack{i,j =1\\ i \neq j}}^N f(x_i)g(x_j)\Bigg].
\end{equation}
Note that under (M2) and \ref{H4}, for any \(f,g \in B(E)\), \(V[f,g]\) is uniformly bounded in \(x\). The operator \(V\) is the point at which the non-local branching mechanism enters the covariance structure of the aforementioned limiting Gaussian processes. In the local case, the offspring locations
coincide with the parent location, so this term reduces to a local expression
involving products of \(f(x)\) and \(g(x)\). In the present setting, \(V[f,g](x)\)
depends on the joint distribution of the locations of distinct offspring born
from a parent at \(x\), and therefore records the spatial correlations between
siblings created at a branching event.

Next, let
$$\mathcal{W}(x):=\{W_{i,j}^{(k)}(x),\text{ \(1\leq i \leq m_L\), \(1\leq j \leq p_i\), \(1\leq k \leq k_{i,j}\)\}},$$
where we recall that the $W_{i, j}^{(k)}$ were given in Theorem \ref{thm: main result large} and $m_L$ was defined after the statement of \ref{H1b}. 

With this notation in hand, for \(f,g \in B(E)\), $x \in E$, define
\begin{align}
\mathcal{C}^{(1)}_{t}(f,g,\mathcal W(x)) &= \frac{1}{(p_1-1)!}\sum_{i=1}^{k_{1,1}}W_{1,1}^{(i)}(x)\sum_{j,k=1}^{m_{L}}\e^{(\frac{\lambda_1}{2}-\lambda_j)t} \notag \\
&\quad \times \left(\int_{0}^{\infty}\e^{(\lambda_1-\lambda_j-\lambda_k)u}\tilde\varphi_{1,1}^{(i)}\F[\mathcal{N}^{p_1-1} V\F[\e^{-\mathcal{N}(t+u)}\Phi_{j}[f],\e^{-\mathcal N u}\Phi_{k}[g]\R]\R]\mathrm{d}u \right.\notag \\
&\hspace{8cm}\left.-\tilde\varphi_{1,1}^{(i)}[\e^{-\mathcal{N}t}\Phi_j[f]\Phi_k[g]]\right), \label{eq:C1}\\
\mathcal{C}^{(2)}_{t}(f,g,\mathcal W(x))
&=\frac{\e^{-{\lambda_1}t}}{(p_1-1)!}\sum_{i=1}^{k_{1,1}}W_{1,1}^{(i)}(x)\tilde\varphi_{1,1}^{(i)}[\mathcal{N}^{p_1-1}\Cov_{\cdot}(X_{t}[f],X_{t}[g])],\label{eq:C2}
\end{align}
where $\Cov_x$ and $\Var_x$, $x \in E$, denotes, respectively, the covariance and variance with respect to $\mathbb P_{\delta_x}$. In addition, whenever \(m_C<m\) define, for \(f,g\) satisfying $\Phi_i[f] = \Phi_i[g] = 0$, $1 \le i \le m_C$, 
\begin{align}
\mathcal{C}^{(3)}_{t}(f,g,\mathcal W(x))&= \frac{\e^{-\frac{\lambda_1}{2}t}}{(p_1-1)!}\sum_{i=1}^{k_{1,1}}W_{1,1}^{(i)}(x)\bigg(\tilde\varphi_{1,1}^{(i)}[\mathcal{N}^{p_1-1}f\psi_{t}[g]]\notag \\
&\hspace{3cm}+\int_0^{\infty}\e^{-\lambda_1 u}\tilde\varphi_{1,1}^{(i)}\F[\mathcal{N}^{p_1-1}V[\psi_u[f],\psi_{t+u}[g]]\R]\mathrm{d}u\bigg). \label{eq:C3}
\end{align}

\medskip

{\color{black}Note that the integral in the definition of \(\mathcal{C}^{(1)}\) converges since, for \(1\leq i,j \leq m_L\), we have \(\mathrm{Re}\lambda_i+\mathrm{Re}\lambda_j>\lambda_1\). The integral in the definition of \(\mathcal{C}^{(3)}\) converges since \(\Phi_i[f] = \Phi_i[g] =0, 1 \le i \le m_C\) and \(m_C<m\), so the growth of \(\psi_t[f]\) and \(\psi_t[g]\) is tied to the growth of some small\footnote{i.e. satisfying $2{\rm Re}\lambda_i < \lambda_1$} eigenvalues by \ref{H1b}.} 

\medskip

We also set 
\begin{equation}
f_1^{(t)} = \sum_{i=1}^{m_{L}}\sum_{j=1}^{p_i}\sum_{k=1}^{k_{i,j}}\tilde\varphi_{i,j}^{(k)}[f]\e^{(\frac{\lambda_1}{2}-\lambda_i-\mathcal{N})t}\varphi_{i,j}^{(k)}, \quad f_2^{(t)} = f-f_1^{(t)},
\label{eq:decomp}
\end{equation}
with the simplification $f_1 = f_1^{(0)}$, $f_2 = f_2^{(0)}$. Thus \(f_1=f_1^{(0)}\) is the projection of \(f\) onto the large spectral modes, while \(f_2\) is the remaining component. The time-dependent version \(f_1^{(t)}\) incorporates the rescaling of these large modes that appears in the small regime covariance.

\medskip

\begin{theorem}[Small Regime]
\label{thm: main result small}
 Assume that (M\(4\)) holds and let \(f \in C_{\psi}(E)\). Furthermore, assume that \(m_{L}=m_{C}<m\). Then, for any initial position \(x \in E\), as $n \to \infty$ 
    \begin{equation}
    \label{eq: main result small}
        \e^{-\frac{\lambda_1}{2} (n+t)}n^{-\frac{p_1-1}{2}}\F(X_{n+t}[f]-\sum_{i=1}^{m_L}\sum_{j=1}^{p_i}\sum_{k=1}^{k_{i,j}}\tilde\varphi_{i,j}^{(k)}[\e^{(\lambda_i+\mathcal{N})(n +t)}f]W_{i,j}^{(k)}(x) \R)\conindis Z_S^f(t),
    \end{equation}
    in \(\mathbb{D}[0,\infty)\), where, for \(f,g \in C_{\psi}(E)\), conditionally on the collection \(\mathcal{W}(x)\), \(Z_S^f\) and \(Z_S^g\) are jointly complex mean-zero Gaussian processes, with covariance, for \(0\leq r \leq t < \infty\), 
\begin{align}
\E[Z_S^f(r)Z_S^g(t)|\mathcal{W}(x)] =  \mathcal{C}^{(1)}_{t-r}(f_1,g_1,\mathcal W(x)) - \mathcal{C}^{(2)}_{t-r}(f_1^{(t-r)}, g_2,\mathcal W(x)) + \mathcal{C}^{(3)}_{t-r}(f_2, g_2,\mathcal W(x)), \label{eq: small fluct 1}\\
         \E[Z_S^f(r)\bar{Z}_S^g(t)|\mathcal{W}(x)] =  \mathcal{C}^{(1)}_{t-r}(f_1,\bar g_1,\mathcal W(x)) - \mathcal{C}^{(2)}_{t-r}(f_1^{(t-r)},\bar g_2,\mathcal W(x)) + \mathcal{C}^{(3)}_{t-r}(f_2,\bar g_2,\mathcal W(x)).\label{eq: small fluct 2}
    \end{align}
  \end{theorem}  
  
  \begin{theorem}[Small regime f.d.d.]\label{cor:small}
  Assume that (M\(4\)) holds and let $f \in B(E)$. Then \eqref{eq: main result small} holds in the sense of finite dimensional distributions. 
  \end{theorem}

\begin{remark}~\label{rem:small}
\begin{itemize}
\item {Note that \(m_L=m_C<m\) corresponds to the setting of no critical eigenvalues and at least one small eigenvalue. As we will see in Theorem \ref{thm: main result crit}, the presence of critical eigenvalues leads to different asymptotic behaviour. If \(m_L\neq m_C<m\), then Theorems \ref{thm: main result small} and \ref{cor:small} can still be recovered, but we must restrict to test functions \(f\) that also satisfy
\begin{equation*}
    \Phi_i[f]=0,\quad m_L+1\leq i \leq m_C.
\end{equation*}
In words, the projection of \(f\) into the critical eigenspaces is 0.} 
\item {If the eigenvalues in  \ref{H1b} only consist of large and critical eigenvalues, then Theorems \ref{thm: main result small} and \ref{cor:small} can still be recovered as long as the gap between $\lambda_m$ and the next element of the spectrum of \(\psi\) is sufficiently large. Formally, this is the case if one can take a ``dummy" eigenvalue $\lambda_{m+1}$ with $\lambda_{m+1}< \lambda_1/2$ and trivial eigenspace \(A_{\lambda_{m+1}}=\{0\}\). Any such trivial eigenspace \(A_{\lambda_{m+1}}\) will satisfy Definition \ref{definition: eigenspace}. The non-trivial property required is that \eqref{eq: h1b assum} holds under the renormalisation \(\e^{-\lambda_{m+1}t}t\). As expected, the consequences of the theorems are independent of the ``dummy" eigenvalue chosen, since it does not appear in the covariance functions \(\mathcal{C}^{(1)},\mathcal{C}^{(2)},\mathcal{C}^{(3)}.\)}
\item {The small regime limit can be decomposed into two components. Fluctuations of \(X[\varphi_{i,j}^{(k)}]\) around their martingale limits, \(W_{i,j}^{(k)}(x)\), whose covariance is described by \(\mathcal{C}^{(1)}\), and fluctuations from ``small" functions $f$ satisfying $\Phi_i[f] = 0, 1 \le i \le m_L$ whose covariance is described by \(\mathcal{C}^{(3)}\). Then, \(\mathcal{C}^{(2)}\) describes the covariance between these two components. Note that \(\mathcal{C}^{(2)}_0=0\) regardless of \(f\) and \(g\). Thus, if one only considers convergence in distribution, then these components appear independent and an understanding of the correlation between them is lost. This highlights the importance of the functional convergence result.}
\end{itemize}
\end{remark}

\medskip

To state our final theorem for the critical regime, we require some more notation. {For \(m_L+1\leq i \leq m_C\), define
\begin{align*}
    {\rm Ei}(\lambda_i) = \left\{f\in B(E): \Phi_j[f]=0 \text{ for all } m_L+1\leq j\leq m_C,\ j\neq i \right\}, 
\end{align*}
where \(\mathrm{Ei}_1(\lambda_i)\) is the restriction to such functions bounded by unity. Furthermore, let \[\mathrm{Ei}(\Lambda_C) = \bigcup_{\substack{m_{L}+1\leq i \leq m_C}}\mathrm{Ei}(\lambda_i).\]
In words, \(\mathrm{Ei}(\Lambda_C)\) consists of all functions that have non-zero projection in at most one critical eigenspace. For \(f\in \mathrm{Ei}(\Lambda_C)\), if there is a unique index \(\nu_f\) such that \(f \in \mathrm{Ei}(\lambda_{\nu_f})\), set \(\lambda_f:=\lambda_{\nu_f}\) and define
\begin{align*}
p_f:= \max \{i\geq 1: \mathcal{N}^{i-1}\Phi_{\nu_f}[f] \neq 0\}.
\end{align*}
Otherwise, take \(\lambda_f = \lambda_1/2\), \(p_f=1\), and define \(\Phi_{\nu_f}[f]=0\). Thus, if \(f \in \mathrm{Ei}(\Lambda_C)\) has non-zero projection in exactly one critical eigenspace, then \(\lambda_f\) is the eigenvalue of this eigenspace. In particular, if \(\lambda_1/2\) is the unique critical eigenvalue, then \(\mathrm{Ei}(\Lambda_C)=B(E)\).} 

Finally, let \(C_{\psi}(\Lambda_C) =\mathrm{Ei}(\Lambda_C) \cap C_{\psi}(E) \), and for \(f,g\in \mathrm{Ei}(\Lambda_C)\), \(r,t\geq 0\), let
\begin{equation*}
    \mathcal{C}^{(4)}_{r,t}(f,g) := \bs1_{\lambda_f = \bar \lambda_g}\int_0^{r}\frac{(r-v)^{p_f-1}v^{p_1-1}(t-v)^{p_g-1}}{(p_f-1)!(p_1-1)!(p_g-1)!}\mathrm{d}v.
\end{equation*}

\medskip

\begin{theorem}[Critical Regime]
\label{thm: main result crit}
Assume that (M\(4\)) holds and let \(f \in C_{\psi}(\Lambda_C)\). {Also assume that \(p_1=1\) in \ref{H1b}}. Then, for any initial position \(x \in E\), as $n \to \infty$
    \begin{equation}
    \label{eq: main result crit}
        \e^{-\lambda_f nt}n^{-\frac{2p_f+p_1-2}{2}}\F(X_{nt}[f]-\sum_{i=1}^{m_L}\sum_{j=1}^{p_i}\sum_{k=1}^{k_{i,j}}\tilde\varphi_{i,j}^{(k)}[\e^{(\lambda_i+\mathcal{N})n t}f]W_{i,j}^{(k)}(x) \R)\conindis Z_C^f(t)
    \end{equation}
    in \(\mathbb{D}[0,\infty)\), where, for \(f,g\in C_{\psi}(\Lambda_C)\), conditionally on the collection \(\mathcal{W}(x)\), \(Z_C^f\) and \(Z_C^g\) are jointly complex mean-zero Gaussian processes, with covariance, for \(0\leq r \leq t < \infty\),
\begin{align}
&\E[Z_C^f(r)Z_C^g(t)|\mathcal{W}(x)] = \mathcal{C}^{(4)}_{r,t}(f,g)\sum_{i=1}^{k_{1,1}}\tilde\varphi_{1,1}^{(i)}[\mathcal{N}^{p_1-1}V[\mathcal{N}^{p_f-1}\Phi_{\nu_f}[f],\mathcal{N}^{p_{g}-1}\Phi_{\nu_g}[g]]]W_{1,1}^{(i)}(x),\label{eq: critical covariance main 1}\\
         &\E[Z_C^f(r)\bar{Z}_C^g(t)|\mathcal{W}(x)] =  \mathcal{C}^{(4)}_{r,t}(f,\bar g)\sum_{i=1}^{k_{1,1}}\tilde\varphi_{1,1}^{(i)}[\mathcal{N}^{p_1-1}V[\mathcal{N}^{p_f-1}\Phi_{\nu_f}[f],\mathcal{N}^{p_{g}-1}\Phi_{\nu_{\bar g}}[\bar {g}]]]W_{1,1}^{(i)}(x). \label{eq: critical covariance main 2}
    \end{align}
  \end{theorem}
  
  \begin{theorem}[Critical regime f.d.d.]\label{cor:crit}
  Suppose (M\(4\)) holds and take \(f\in \mathrm{Ei}(\Lambda_C)\). Then \eqref{eq: main result crit} holds in the sense of finite dimensional distributions.
  \end{theorem}
  {Note that for functional convergence in the critical regime we require that the leading eigenvalue is of rank 1, \(p_1=1\). This is a technical assumption required due to our method of proof of tightness which can be relaxed if one is only interested in finite dimensional distributions.  If, for \(f\in \mathrm{Ei}(\Lambda_C)\), \(\lambda_f\neq \lambda_1/2\), then the exponential renormalization in \eqref{eq: main result crit} oscillates in the complex plane. This has interesting consequences and is the reason we are required to work with test functions in \(\mathrm{Ei}(\Lambda_C)\). To see why, first consider \(g = f+\bar f\). By \eqref{eq: complex eigenspace}, we have that \(\lambda_{\bar f}=\bar{\lambda}_f\) and \(p_{\bar f}=p_f\). Thus, informally we may apply Theorem \ref{cor:crit} to obtain
  \begin{equation}
  \label{eq: limit simple critical eigen}
              \e^{-\lambda_1 nt/2}n^{-\frac{2p_f+p_1-2}{2}}\F(X_{nt}[g]-\sum_{i=1}^{m_L}\sum_{j=1}^{p_i}\sum_{k=1}^{k_{i,j}}\tilde\varphi_{i,j}^{(k)}[\e^{(\lambda_i+\mathcal{N})n t}g]W_{i,j}^{(k)}(x) \R)\conindis \mathrm{Re}\F(\e^{(\lambda_f-\lambda_1/2)nt}Z_C^f(t)\R)
  \end{equation}
  in the sense of finite dimensional distributions, where indeed this only holds informally as the proposed limit oscillates with \(n\). However, by considering a single time point and using that \(Z_C^f(1)\) is radially symmetric in law, Theorem \ref{cor:crit} implies
  \begin{equation}
  \label{eq: real conv in dis critical}
       \e^{-\lambda_1 n/2}n^{-\frac{2p_f+p_1-2}{2}}\F(X_{n}[g]-\sum_{i=1}^{m_L}\sum_{j=1}^{p_i}\sum_{k=1}^{k_{i,j}}\tilde\varphi_{i,j}^{(k)}[\e^{(\lambda_i+\mathcal{N})n }g]W_{i,j}^{(k)}(x) \R)\conindis \mathrm{Re}\F(Z_C^f(1)\R)
  \end{equation}
  as \(n\rightarrow \infty\). This clearly cannot be extended to the finite dimensional distribution setting, since the covariance function of \(\mathrm{Re}\F(\e^{(\lambda_f-\lambda_1/2)nt}Z_C^f(t)\R)\) oscillates in \(n\) when different time points are considered. Thus, in the functional setting, the ``large \(n\) limit" behaves as a Gaussian process with sinusoidal oscillations that have period \({2 \pi}({\left|\operatorname{Im}\left(\lambda_f-\lambda_1 / 2\right)\right| n})^{-1}.\)
  
  More generally, since for any \(g \in B(E)\), we have that \(\Phi_{m_L+1}[g],\dots,\Phi_{m_C}[g],(g-\Phi_{m_L+1}[g]-\dots-\Phi_{m_C}[g])\in \mathrm{Ei}(\Lambda_C)\), Theorem \ref{cor:crit} implies
  \begin{equation}
  \label{eq: limit general critical eigen}
      \e^{-\lambda_1 n/2}n^{-\frac{2p_g^*+p_1-2}{2}}\F(X_{n}[g]-\sum_{i=1}^{m_L}\sum_{j=1}^{p_i}\sum_{k=1}^{k_{i,j}}\tilde\varphi_{i,j}^{(k)}[\e^{(\lambda_i+\mathcal{N})n }g]W_{i,j}^{(k)}(x) \R)\conindis \sum_{i\in P_g}Z_C^{\Phi_{i}[g]}(1),
  \end{equation}
  where \(p_g^*:=\mathrm{max}\{p_{\Phi_i[g]},m_{L}+1\leq i \leq m_C\}\), \(P_g=\{m_{L}+1\leq i \leq m_C:p_{\Phi_i[g]}=p_g^*\}\), and we have used that 
  \begin{equation*}
      \sum_{i\in P_g}Z_C^{\Phi_{i}[g]}(1) \eqindis \sum_{i\in P_g}\e^{(\lambda_i-\lambda_1/2)n}Z_C^{\Phi_{i}[g]}(1), \quad n\geq 0,
  \end{equation*}
  by \eqref{eq: critical covariance main 1} and \eqref{eq: critical covariance main 2}. Identically to \eqref{eq: limit simple critical eigen}, unless \(\lambda_1/2\) is the only critical eigenvalue this cannot be formally extended to a functional result due to \(n\)-dependent oscillations of the limit process. This highlights the importance of functional convergence in the critical regime, since convergence at a single time point masks the oscillatory behaviour of the limit due to radial symmetry of the complex Gaussian.
    \begin{remark}~
    \label{remark:small and critical regimes}
    \begin{itemize}
        \item Using a similar approach to that of \cite{bJanson}, it should be straightforward to show Theorems \ref{thm: main result small}-\ref{cor:small} and Theorems \ref{thm: main result crit}-\ref{cor:crit} hold jointly, where conditionally on \(\mathcal{W}(x)\), the limit processes \(Z_S^f\) and \(Z_C^g\) are independent. Since usually only the asymptotically largest fluctuations are of interest, we do not give proof here for the sake of brevity.
        \item We give some technical comments on the different time scales seen in the small and critical regimes. Recall that the heuristics as to why these time scales appear were given at the end of Section 1. By \eqref{eq: main result crit}, we have on the timescale \((n+t)_{t\geq 0}\), the critical regime converges to \(Z_C^f(1)\). By exponential decay of the covariance function in Theorem \ref{thm: main result small}, on the timescale of \((nt)_{t\geq 0}\), the finite dimensional distributions of \eqref{eq: main result small} are i.i.d.\ Gaussian. This can be shown using the same proof techniques of Theorem \ref{cor:small}. Thus, we do not obtain a functional limit theorem on this timescale. However, note that if there are no critical eigenvalues, then \(\mathrm{Ei}(\Lambda_C)=B(E)\), and the limit in \eqref{eq: main result crit} is equal to 0 for all \(f \in B(E)\). Thus, under the additional polynomial rescaling of \(n^{-1/2}\), \eqref{eq: main result small} converges to 0 on the timescale \((nt)_{t\geq 0}\). 
        \end{itemize}
    \end{remark}}

\section{Discussion}\label{sec:discussion}
The purpose of this section is threefold. We first discuss assumption \ref{H1b} in further detail and provide sufficient conditions for it to hold. Next we discuss the restriction $f \in C_\psi(E)$ and explain why this condition is natural for functional convergence. Finally, we present several examples that demonstrate different behaviours covered by the theory in the previous section. 

\subsubsection*{Assumption \ref{H1b}}
In practice, showing Assumption \ref{H1b} is satisfied is a non-trivial task. We therefore give sufficient conditions for this assumption to hold in terms of the underlying Markov process and associated many-to-one Markov process that are easier to check.

To this end, let us introduce the aforementioned many-to-one formula. By \cite[Lemma 8.2]{bNTEbook}, under (M\(1\)) and \ref{H4}, for any \(f \in B(E)\), \(x \in E\), \(t\geq 0\), we have
\begin{equation}
    \psi_t[f](x) = \hat {\bf E}_x\F[\mathrm{exp}\F(\int_0^t\gamma(\hat\xi_s)\F(\mathcal{E}_{\hat\xi_s}[N]-1\R)\mathrm{d}s\R)f(\hat\xi_t)\R],
    \label{eq: sufficient for semigroup in L2 1}
\end{equation}   
{where $(\hat\xi, \hat{\mathbf P})$ is a Markov process on $E \cup \{\dagger\}$ which evolves according to $(\xi, {\mathbf P})$ and independently with rate \(\gamma(x)\mathcal{E}_x[N]\) is sent to a new location in \(E\), such that for each Borel \(A \subseteq E\), the probability of the new position being in \(A\) is \(\mathcal{E}_x[\mathcal{Z}[\bs 1_A]]/\mathcal{E}_x[N]\).} We also let, for \(f\in B(E)\), \(t\geq 0\), \(x\in E\),
\[
    \hat P_t[f](x) = \hat {\bf E}_x[f(\hat \xi_t)]
\]
denote the semigroup associated with $\hat\xi$. 

As we will shortly see, the process \(\hat\xi\) will be used as a tool for verifying \ref{H1b}. Indeed, in applications it is often easier to check compactness, smoothing, or density assumptions for the Markov semigroup \(\hat P_t\) than directly for \(\psi_t\). Proposition \ref{prop: conditions for Linfinity spectrum} below formalises this idea by giving sufficient conditions, stated in terms of \(\hat P_t\), under which the \(L^\infty\) spectral decomposition required in \ref{H1b} holds.

\begin{proposition}
\label{prop: conditions for Linfinity spectrum}
   Let \(E\) be a locally compact Hausdorff space and \(\mu\) a finite measure on \(E\). Assume (M\(1\)), \ref{H4} and the following two conditions.
   \begin{enumerate}
       \item There exists \(T>0\), such that, for all \(t\geq T\), there exists \((p_{t}(x,\cdot),x\in E)\), such that 
       \begin{equation}
           \hat P_t[f](x) = \int_E p_t(x,y)f(y)\mu(\mathrm{d}y), \quad f \in B(E), \quad x \in E, \label{eq: semigroup in L2 1}
       \end{equation}
       and
       \begin{equation}
             \sup_{x\in E}\|p_t(x,\cdot)\|_{L^2(E,\mu)}<\infty. \label{eq: semigroup in L2 2}
       \end{equation}
       \item There exists a dense subset, \(C(E,\mu)\), of \(L^2(E,\mu)\), such that, \(C(E,\mu) \subseteq B(E)\), and for each \(f\in C(E,\mu)\), \[\lim_{t\rightarrow 0}\|\hat P_t[f]-f\|_{L^2(E,\mu)}=0.\]
   \end{enumerate} 
  {Then, for any \(K\in \mathbb{R}\), there exist\footnote{It is possible for \(m=0\), in which case the list \(\lambda_1,\dots,\lambda_m\) is empty} \(\lambda_1,\dots,\lambda_m\in \mathbb{C}\) with \(\mathrm{Re}\lambda_1\geq \dots \geq \mathrm{Re}\lambda_m\geq K\), such that \(\lambda_1,\dots,\lambda_m\) are eigenvalues of \(\psi\) with finite dimensional eigenspaces \(A_{\lambda_1},\dots,A_{\lambda_m}\). Furthermore, using the notation of \ref{H1b} we have}
\begin{equation}
    \label{eq: h1b assum example}
  \sup_{x \in E, f \in B_1(E)}\e^{-K t}t\bigg|\psi_t[f](x) - \sum_{i=1}^{m}\e^{(\lambda_i+\mathcal{N})t}\Phi_{i}[f](x)\bigg|
  \to 0, \quad \text{as } t \to \infty.
\end{equation}
\end{proposition}

The proof of the above proposition is given in Appendix \ref{sec:H1+}.  

\begin{remark}
Note that the above proposition does not imply that $\lambda_1$ is real. Thus, for \ref{H1b} to hold for a general supercritical MBP, we implicitly assume that this is true. This is satisfied if \ref{Ha} holds, which can be shown to hold for a large class of processes, see \cite{bNTEbook} and references therein.
\end{remark}

Let us now compare the assumptions given in Proposition \ref{prop: conditions for Linfinity spectrum} to those made in the literature. First note that in the case that \(E\) is finite, Proposition \ref{prop: conditions for Linfinity spectrum} holds with \(\mu\) as the counting measure. Thus, we extend the classical results of \cite{bJanson} to the case where \(\lambda_1\) is non-simple. 

Next, consider the setting of \cite{CLTNS}, where the authors consider (local) branching Markov processes. Here the authors assume conditions that allow them to derive a spectral decomposition of $\psi$ in \(L^2(E,\mu)\). In particular, their assumptions satisfy those of Proposition \ref{prop: conditions for Linfinity spectrum} except for the assumptions that \(\mu\) is finite and \eqref{eq: semigroup in L2 2} holds, where \(C(E,\mu)\) is taken to be the space of continuous compactly supported functions. Indeed, \cite{CLTNS} only assumes \(p_t(\cdot,\cdot)\) is in \(L^2(E\times E,\mu \times \mu)\). In our setting, these two additional assumptions are used to transfer the spectral decomposition of \(\psi\) in \((L^2(E,\mu),\|\cdot\|_{L^2(E,\mu)})\) to a spectral decomposition in \((B(E),\|\cdot\|_{\infty})\). At least heuristically, an assumption of the form \eqref{eq: semigroup in L2 2} seems necessary to obtain this stronger version of the spectral decomposition. Indeed, in \ref{H1b}, we have that the asymptotic behaviour of \(\psi_t\) is characterised by pairs of the form \((\varphi_{i,j}^{(k)}, \tilde\varphi_{i,j}^{(k)})\), where \(\tilde\varphi_{i,j}^{(k)}\) does not depend on \(x\). If \eqref{eq: semigroup in L2 2} did not hold, then this would suggest that \(\varphi_{i,j}^{(k)}\notin B(E)\), and hence \ref{H1b} would not hold.

We see two main benefits of Proposition~\ref{prop: conditions for Linfinity spectrum} compared to the setting of \cite{CLTNS}. Firstly, having a spectral decomposition of the form given in \ref{H1b} is a natural extension of the assumptions typically made for non-local branching processes in recent literature; see \cite{bNTEbook} and references therein. In particular, \ref{H1b} allows us to exploit recent results on the asymptotic behaviour of the moments of non-local BMPs, originally developed in \cite{bmoments} and extended in \cite{Moments_Chris}, which are crucial to proving the functional CLT results given in Section~\ref{sec: main results general}. These results are formulated in the \(L^\infty\) norm, and this is natural for the fluctuation theory developed here. Indeed, although spectral decompositions are often most naturally obtained in an \(L^2(E,\mu)\) setting, the branching process may be initiated from a point mass \(\delta_x\), and the quantities of interest are the random variables \(X_t[f]\), their martingale limits and their covariances, all of which depend pointwise on the initial state \(x\). The moment estimates used below are therefore uniform in \(x\in E\) and \(f\in B_1(E)\), rather than merely \(L^2(E,\mu)\) estimates. Moreover, since \(X_t\) is a random finite measure, \(X_t[f]\) is defined for bounded measurable test functions independently of any reference measure, whereas an \(L^2(E,\mu)\) formulation identifies functions only up to \(\mu\)-null sets and does not directly control \(\psi_t[f](x)\) uniformly in \(x\).

This distinction is also relevant from the point of view of simulation and numerical approximation. An \(L^2(E,\mu)\) spectral expansion controls the approximation error only after averaging over the reference measure \(\mu\), while Monte Carlo approximations of branching systems are typically built from particles started at specified locations and observed along realised trajectories. Regions of small \(\mu\)-measure may nevertheless contribute significantly to the growth or variance of such estimators, for example through high reproductive weight, rare but influential offspring displacements, or boundary effects. The \(B(E)\)-formulation of \ref{H1b} gives precisely the uniform control of \(\psi_t[f](x)\), martingale limits and covariance terms needed for such point-mass and simulation-based interpretations.

Secondly, Proposition~\ref{prop: conditions for Linfinity spectrum} only asks for \(\hat P_t\) to have an \(L^2(E,\mu)\) density for \(t\) sufficiently large. This is particularly useful in the non-local setting, since it allows the underlying Markov motion to be singular, provided that the non-local branching mechanism produces sufficient smoothing at positive times. A typical example of such dynamics is the neutron branching process, \cite{bSNTE2,bNTEbook,bSNTE1}. Thus, our approach does not preclude the use of \(L^2\) spectral theory; rather, Proposition~\ref{prop: conditions for Linfinity spectrum} shows how \(L^2\) methods can be used to verify the \(L^\infty\)-type assumptions required for the branching fluctuation theory in concrete examples.

\subsubsection*{Restriction of $f$ to $C_{\psi}(E)$}
As briefly discussed in the previous section, the restriction to $C_{\psi}(E)$ can be thought of as a regularity assumption on the branching process, since it is a condition on the Markov process, $\xi$, the state space, $E$, and the function, $f$. We start by discussing this class in specific settings and then show that the regularity imposed by this restriction transfers to similar regularity on the branching process.

{Since (M\(4\)) is required for the vast majority of our results, we consider $C_{\psi}(E)$ under the minimal assumption (M\(4)\). Let \(d_E\) be a metric on \(E\). Then $C_{\psi}(E)$ contains all bounded Lipchitz functions with respect to \(d_E\) whenever there exists a constant \(C>0\), such that
\begin{align}
                &\sup_{x \in E,0\leq t \leq {\color{black}1}} t^{-1}\mathbf E_x[d_E(\xi_{t},x)^{2}]\leq C,\nonumber \\
                &\sup_{x \in E,0\leq t \leq {\color{black}1}} t^{-1}\mathbf E_x[d_E(\xi_{t},x)^{4}]\leq C\label{eq: mild conditions}. 
                \end{align}
Taking \(d_E\) to be the standard Euclidean norm, \eqref{eq: mild conditions} holds for piecewise deterministic Markov processes with uniformly bounded velocities and uniformly bounded jump rates, as well as multitype branching processes, see Examples \ref{example: mtbp} and \ref{example: nbp} below. Other examples include L\'evy processes in bounded domains with absorbing boundary conditions and whose L\'evy measure has finite fourth moments. Similarly, Itô diffusions (in such domains) which take the form 
\begin{equation*}
    d\xi_t = b(\xi_t)dt + \sigma(\xi_t)dW_t,
\end{equation*}
with \(b,\sigma\) bounded.}

\subsubsection*{Examples}
In this section we discuss several examples where \ref{H1b} holds. The examples are chosen to illustrate different aspects of the theory including a non-simple leading eigenvalue, Jordan blocks and polynomial corrections arising from the existence of generalised eigenfunctions and a non-local spatial model where the many-to-one process regularises the otherwise singular motion. We will assume throughout that \ref{H4} and (M4) are in place.

\begin{example}\label{eq:example BBM}
Our first example demonstrates a case when the leading eigenvalue is not simple by introducing an immigration mechanism into the branching process. 

Let $D \subset \mathbb R^d$ be an open, bounded and connected, and let $D_1$ and $D_2$ be two disjoint copies of $D$. We introduce an additional state $I$, which should be interpreted as a persistent source of immigrants. The full state space is therefore $E = I \cup D_1 \cup D_2$. 

Particles in $D_1$ and $D_2$ evolve independently according the dynamics of (supercritical) branching Brownian motion with constant branching rate and absorbing boundaries. On the other hand, a particle in $I$ does not move and branches at a constant rate $\gamma_I$. When a branching event occurs in $I$ the parent particle is replaced by one offspring in $I$, and $N_I$ additional offspring are sent to $D_1 \cup D_2$, where $N_I \ge 1$ almost surely and has finite fourth moments. The locations of the additional offspring in $D_1 \cup D_2$ are distributed uniformly (with respect to Lebesgue measure).  We note that while the process may be viewed as a branching process with immigration into $D_1 \cup D_2$, the immigration is endogenous  and hence the process remains a branching Markov process on $E$. 

Due to the dynamics of the particles in $D_1 \cup D_2$ and the offspring distribution from branching events in $I$, it is straightforward to check that properties $1$ and $2$ of Proposition \ref{prop: conditions for Linfinity spectrum} are satisfied with \(\mu\) Lebesgue on \(D_1\cup D_2\) and the counting measure on \(I\). 

The purpose of this example is to demonstrate a situation when the leading eigenvalue of the system is not simple, which arises due to the reducibility of the system. Indeed, initiating the process with one particle in \(D_1\) (resp.\ \(D_2\)) restricts all particles to \(D_1\) (resp.\ \(D_2\)). As such, there are two independent supercritical regions, each with the same leading eigenvalue and thus the first order LLN is no longer described by a single martingale limit but by two martingale limits (one for each region). In the context of \ref{H1b}, we have $p_1 = 1$ and $k_{1, 1} = 2$ where \(\tilde\varphi_{1,1}^{(1)}\) (resp.\ \(\tilde\varphi_{1,1}^{(2)}\)) is a probability measures with support \(D_1\) (resp.\ \(D_2\)). In particular, these are exactly the probability measures seen in \ref{Ha} if the process is restricted to \(D_1\) (resp.\ \(D_2\)). In the large regime, Corollary \ref{thm: SLLNs} gives 
\[
{\rm e}^{-\lambda_1 t}X_t[f] \to \tilde\varphi_{1,1}^{(1)}[f] W_{1,1}^{(1)}(x) + \tilde\varphi_{1,1}^{(2)}[f] W_{1,1}^{(2)}(x).
\]

One can consider many natural extensions of Example \eqref{eq:example BBM} that still satisfy the conditions of Proposition \ref{prop: conditions for Linfinity spectrum}. For example, one can take the underlying Markov motion on \(D\) to satisfy the assumptions of either \cite{bCLT} or \cite{CLTNS} along with \(\mu\) finite and \eqref{eq: semigroup in L2 2}. Note that since we assume branching is local on \(D\), the many-to-one Markov process restricted to \(D_1 \cup D_2\) is exactly the underlying Markov process. Then, under the additional assumption that the immigration offspring law has an \(L^2\) density with respect to \(\mu\), Proposition \ref{prop: conditions for Linfinity spectrum} holds for the immigration process as defined in Example \ref{eq:example BBM}. Under these assumptions, the only non-trivial component of Proposition \ref{prop: conditions for Linfinity spectrum} to check is that, for \(x \in I\), \(\|p_t(x,\cdot)\|_{L^2(D_1\cup D_2,\mu^*)}<\infty\), where the restriction of \(\mu^*\) on \(D_1\) or \(D_2\) is equal to \(\mu\). By our assumption on the immigration offspring law, there exists \(q \in L^2(D_1\cup D_2,\mu^*)\), such that, for \(x \in I\),
    \begin{equation}
\mathcal{E}_x[\mathcal{Z}[f]]=f(x)+\int_{D_1\cup D_2}q(y)f(y)\mu^*(\mathrm{d}y), \quad f \in B(E). \label{eq: imigration behaviour example}
    \end{equation}
Therefore,
        \begin{equation}
        p_t(x,y) =\frac{1}{\mathcal{E}_x[\mathcal{Z}[\bs1_{D_1\cup D_2}]]}\int_0^t\gamma_I^*\e^{-\gamma_I^* s}\int_{D_1\cup D_2}p_{t-s}(z,y)q(z)\mu^*(dz)\mathrm{d}s, \quad y \in D_1\cup D_2, \label{eq: control on density}
    \end{equation}
  where \(\gamma_I^*=\gamma_I\mathcal{E}_{x}[\mathcal{Z}[\bs1_{D_1\cup D_2}]]\). By the assumptions of \cite{bCLT} or \cite{CLTNS}, we have that \((\hat{P}_t)_{t\geq 0}\) as a semigroup on \(L^2(D_1\cup D_2,\mu^*)\) is a contraction semigroup. This, Jensen's inequality, Fubini's theorem, Cauchy-Schwarz and \eqref{eq: control on density} imply, for \(x \in I\),
    \begin{align*}
       \|p_t(x,\cdot)\|_{L^2(D_1\cup D_2,\mu)}^2&=\frac{1}{\mathcal{E}_x[\mathcal{Z}[\bs1_{D_1\cup D_2}]]^2}\int_{D_1\cup D_2} \F(\int_0^t\gamma_I^*\e^{-\gamma_I^* s}\int_{D_1\cup D_2}p_{t-s}(z,y)q(z)\mu(dz)\mathrm{d}s\R)^2\mu(dy)\\
       &\leq \frac{1}{\mathcal{E}_x[\mathcal{Z}[\bs1_{D_1\cup D_2}]]^2}\int_{D_1\cup D_2} \int_0^t\gamma_I^*\e^{-\gamma_I^* s}\F(\int_{D_1\cup D_2}p_{t-s}(z,y)q(z)\mu(dz)\R)^2\mathrm{d}s\mu(dy)\\
       & \leq \frac{1}{\mathcal{E}_x[\mathcal{Z}[\bs1_{D_1\cup D_2}]]^2}\int_{D_1\cup D_2} \int_0^t\gamma_I^*\e^{-\gamma_I^* s}\int_{D_1\cup D_2}p_{t-s}(z,y)q(z)^2\mu(dz)\mathrm{d}s\mu(dy)\\
       &\le \frac{\|q\|_{L^2(D_1\cup D_2,\mu^*)}^2}{\mathcal{E}_x[\mathcal{Z}[\bs1_{D_1\cup D_2}]]^2}<\infty
    \end{align*}
    as required, where in the penultimate line we have used that \(\|p_t(\cdot,x)\|_{L^1(D_1\cup D_2,\mu^*)}\leq 1\) for all \(x \in D_1\cup D_2\) and all \(t\geq0\) (this holds by the assumptions of \cite{bCLT} or \cite{CLTNS}). In particular, these assumptions are satisfied for the class of Markov processes given in Example 1.5 of \cite{bCLT}, and the Markov processes in Examples 1.1-1.7 of \cite{CLTNS}.

In addition, the same construction can be made with finitely or countably many communicating classes, provided the reference measure is finite and the $L^2$-density bounds in Proposition \ref{prop: conditions for Linfinity spectrum} hold uniformly over the classes. {In particular, this allows for the immigration set to be as general as \(D\).}

\end{example}

\begin{example}\label{example: mtbp}
While the previous example provides a relevant setting where the leading eigenvalue is non-simple, the leading eigenvalue only has eigenfunctions of rank 1. This example provides a setting where the leading eigenvalue has eigenfunctions of rank 2.

Here we consider the setting of multi-type branching processes. In this case, take \(E=\{1,\dots,d\}\) for some $d \ge 2$. It is well-known (see \cite{bJanson}) that, for any \(f\in B(E)\), \(x \in E\), \(t\geq 0\),
\begin{equation*}
    \psi_t[f](x)=\bs f^T\e^{At}\bs e_x, 
\end{equation*}
where \(\bs e_x\) is the \(x\)th canonical basis vector, \(\bs f = (f(1),\dots,f(d))\), and \(A = (\gamma(j)\mathcal{E}_j[\bs 1_i]-\delta_{ij})_{i,j=1}^d\). Then, (M\(1\)) and \ref{H4} imply \ref{H1b} immediately by taking the Jordan-Normal form of \(A\). Thus, the results of Section \ref{sec: main results general} extend the results of \cite{bJanson} when the leading eigenvalue of \(A\) is non-simple. 

An example of \(A\) with non-simple $\lambda_1$ with generalised eigenfunctions is given by

$$
A=\left(\begin{array}{ll}
1 & 1 \\
0 & 1
\end{array}\right),
$$
so that the semigroup is given by
\[
{\rm e}^{At} = {\rm e}^t
\left(\begin{array}{ll}
1 & t \\
0 & 1
\end{array}\right)
\]

As such, we can see that the type $2$ population grows at rate ${\rm e}^t$ but it also continuously feeds the type $1$ population. Since type $1$ particles themselves grow at the same exponential rate, the cumulative contribution to the type $1$ population from the type $2$ population is of order $t{\rm e}^t$. This interplay between the two growth mechanisms is precisely what produces a Jordan block structure. Indeed, $A$ has eigenvalue 1 with corresponding eigenvector $(1,0)^T$ and generalised eigenvector of rank 2 given by $(0,1)^T$. 

\end{example}

\begin{example} \label{example: nbp}
Our final example is the aforementioned neutron branching process. This is an important model in nuclear physics as it is used to model the evolution of neutrons in fissile material \cite{bNTEbook}. Notably, this is an example that intrinsically contains non-local branching due to how neutrons interact with matter. Furthermore, since the underlying Markov motion is singular, this example demonstrates the need for an $L^\infty$ fluctuation theory. We shall subsequently see that the density of the many-to-one semigroup only exists for sufficiently large times as a consequence of the singular motion and non-local branching, highlighting the importance of Proposition \ref{prop: conditions for Linfinity spectrum}.

The neutron branching process is a stochastic representation of the evolution of neutrons undergoing nuclear fission. In this setting, the state space $E$ is a six dimensional space $D \times V$, where $D \subset \mathbb R^3$ is a open and bounded, representing the set of spatial locations a neutron can take, and $V = \{\upsilon \in \mathbb R^3 : \upsilon_{\min} \le |\upsilon| \le \upsilon_{\max}\}$, $0 < \upsilon_{\min} \le \upsilon_{\max} < \infty$, is the set of possible velocities. The underlying Markov process, $\xi$, is a piecewise deterministic Markov process that, when at $(r, \upsilon) \in D \times V$, changes velocity at rate $\sigma_s(r, \upsilon)$, continuing its motion but now with velocity $\upsilon' \in V$ with probability $\pi_s(r, \upsilon, \upsilon'){\rm d}\upsilon'$. The branching rate is denoted $\sigma_f$, which is also a function of position and velocity. Finally, the offspring distribution is supported on $\{0, \dots, N_{max}\}$, for some $N_{max} \ge 1$ and satisfies
\[
\mathcal E_{(r, v)}[\mathcal Z[g]] = \int_V g(r, \upsilon') \pi_f(r, \upsilon, \upsilon'){\rm d}\upsilon',
\]
for some kernel $\pi_f(\cdot, \cdot, \cdot) : D\times V\times V \to [0, \infty)$. Thus, the branching mechanism is local in the spatial variable, $r$, but non-local in the velocity variable, $\upsilon$.

In practice, the so-called {\it cross-sections}, $\sigma_s$, $\sigma_f$, $\pi_s$ and $\pi_f$, are uniformly bounded above and we shall assume this is the case here. Thus, \ref{H4} and $(M1)$ are satisfied and so it can be shown, see e.g. \cite{bSNTE1, bNTEbook}, that a many-to-one formula in the form of \eqref{eq: sufficient for semigroup in L2 1} holds. In this case, $\hat\xi$ is also a piecewise deterministic Markov process with jump rate $\alpha(r, \upsilon) = \sigma_s(r, \upsilon) + \sigma_f(r, \upsilon)\int_V\pi_f(r, \upsilon, \upsilon'){\rm d}\upsilon'$, and jump kernel $\pi(r, \upsilon, \upsilon') = (\sigma_s(r, \upsilon)\pi_s(r, \upsilon, \upsilon') + \sigma_f(r, \upsilon)\pi_f(r, \upsilon, \upsilon'))/\alpha(r, \upsilon)$. Under the additional assumptions 
\begin{itemize}
\item $\sup_{r, \in D, \upsilon, \upsilon' \in V}\alpha(r, \upsilon)\pi(r, \upsilon, \upsilon') > 0$,
\item there exists $\varepsilon > 0$ such that $D_\varepsilon := \{r \in D : \inf_{y \in \partial D}|y - r| > \varepsilon\}$ is non-empty and connected,
\item there exist $0 < s_\varepsilon < t_\varepsilon$ and $\gamma >0$ such that, for all $r\in D\backslash D_\varepsilon$,there exists $K_r \subset V$ measurable such that ${\rm Vol}(K_r) \ge \gamma > 0$ and for all $\upsilon \in K_r$, $r+\upsilon s \in  D_\varepsilon$ for every $s \in [s_\varepsilon, t_\varepsilon]$ and $r + \upsilon s \notin \partial D$ for all $s \in [0, s_\varepsilon]$,
\end{itemize}
it was shown in \cite{bSNTE1} that \ref{Ha} holds. Moreover, a consequence of the steps in proving this result show that $\hat P_t$ has a density with respect to Lebesgue measure for sufficiently large $t$, and that this density is bounded (uniformly in $x = (r_x, \upsilon_x), y = (r_y, \upsilon_y)$) by a constant, $C_t$ say, and hence \eqref{eq: semigroup in L2 1} and \eqref{eq: semigroup in L2 2} hold. Finally, since $P_t$ is the semigroup of a piecewise deterministic Markov process with intensity $\sigma_s$, which is uniformly bounded from above, Assumption 2 of Proposition \ref{prop: conditions for Linfinity spectrum} is satisfied with \(C(E,\mu)\) the space of Lipschitz continuous functions with respect to the Euclidean norm. Thus, Proposition \ref{prop: conditions for Linfinity spectrum} holds in this setting. Moreover, as mentioned above, \ref{Ha} holds so that in the supercritical setting \ref{H1b} holds for the neutron branching process.
\end{example}

\section{Proofs of the main results}\label{sec:proofs2}
 In this section, we give the proofs of our main theorems. Recall the notation $\Cov_x$ and $\Var_x$, $x \in E$, to denote, respectively, the covariance and variance with respect to $\mathbb P_{\delta_x}$. Also, we use \((\mathcal{F}_t)_{t\geq 0}\) to denote the natural filtration of \((X_t)_{t \ge 0}\). As stated in the main results, we shall assume \ref{H1b} and \ref{H4} are in force throughout the proofs.

\subsection{Large Regime}
Throughout this section we assume (M\(2\)) is in force. We will state explicitly when (M\(4\)) and (M\(k\)) are invoked for the functional and \(L^k\) convergence in Theorems \ref{thm: main result large} and \ref{cor:large} respectively, where whenever (M\(k\)) is mentioned we assume that \(k\geq 2\) is even. We first prove the convergence in Theorem \ref{thm: main result large} along integer times and then upgrade this to continuous time. During the proof of integer time convergence, we shall obtain the \(L^k\) convergence in Theorem \ref{cor:large}. For the proof along integer times, we decompose $f \in B(E)$ as $f = f_1 + f_2$, where $f_1$ and $f_2$ were defined in \eqref{eq:decomp}. 

We first prove the result (along integer times) for $f = f_1$ by studying the martingales defined in \eqref{eq: general martingales}. In particular, we show they have almost sure and $L^2$ limits (\(L^k\) limits under (M\(k\))). We then use the fact that \(f_2 \in B(E)\setminus \bigoplus_{i=1}^{m_L}A_{\lambda_i}\) to study the asymptotics of $X_t[f_2$].


For \(1\leq i  \leq m_L\), \(1\leq j \leq p_i\), \(1\leq k \leq k_{i,j}\), recall the martingales defined in \eqref{eq: general martingales},
    \begin{equation*}
        M_{i,j,t}^{(k)} :=  X_t[\e^{-(\lambda_i+\mathcal{N}) t}\varphi_{i,j}^{(k)}], \quad t\geq0.
    \end{equation*}
By Lemma 3.1 of \cite{Moments_Chris} and (M\(2\)), we have that
\begin{equation*}
   \sup_{t\geq0, x \in E} \E_{\delta_x}[|M_{i,j,t}^{(k)}|^2]<\infty.
\end{equation*}
Therefore, by Doob's martingale convergence theorem, the \(M_{i,j}^{(k)}\) have \(\mathbb{P}_{\delta_x}\)-almost sure and \(L^2\) limits, denoted by \(W_{i,j}^{(k)}(x)\). Furthermore, again by Lemma 3.1 of \cite{Moments_Chris}, we have convergence in $L^k$ if (M$k$) holds.

The next lemma studies the fluctuations of these martingales around their limit. {\color{black} Before introducing this lemma, we first define a type of shift operator for the martingale limits \(W_{i,j}^{(k)}(x)\). For \(x \in E\), we let
\begin{align}
    &\mathcal{N}^*W_{i,1}^{(k)}(x) = 0,\quad 1\leq i \leq m, \quad  1\leq k \leq k_{i,1}, \nonumber\\
    &\mathcal{N}^ *W_{i,j}^{(k)}(x) = W_{i,j-1}^{(k)}(x), \quad 1\leq i \leq m, \quad 2 \leq j \leq p_i, \quad 1\leq k \leq k_{i,j}. \label{eq: shift operator}
\end{align}
Note that this definition is consistent with applying the operator \(\mathcal{N}\) to the eigenfunction in \eqref{eq: general martingales}.}
\begin{lemma}
\label{lemma: conv of phi for fdd}
   For \(1\leq i  \leq m_L\), \(1\leq j \leq p_i\), \(1\leq k \leq k_{i,j}\), and \(x \in E\),
    \begin{equation}
    \label{eq: fluct large a.s. 2}
        n^{-\frac{p_1+1}{2}}\e^{-\frac{\lambda_1 n}{2}}\F(X_n[\varphi_{i,j}^{(k)}]-\e^{(\lambda_i +\mathcal{N}^*)n}W_{i,j}^{(k)}(x)\R) \rightarrow 0, \qquad n \to \infty,
    \end{equation}
    \(\mathbb{P}_{\delta_x}\)-almost surely, for \(n\in \mathbb{N}_{0}\), and in \(L^2\). Moreover, if (M\(k\)) holds, then this convergence holds in \(L^k\).
\end{lemma}
\begin{proof}
    For the \(\mathbb{P}_{\delta_x}\)-a.s.\ convergence, we appeal to a Borel-Cantelli argument and show that for any \(\varepsilon>0\), \(x \in E\),
    \begin{equation*}
        \sum_{n=1}^{\infty}\mathbb{P}_{\delta_x}\F(n^{-\frac{p_1+1}{2}}\e^{-\frac{\lambda_1 n}{2}}\F|X_n[\varphi_{i,j}^{(k)}]-\e^{(\lambda_i +\mathcal{N}^*)n}W_{i,j}^{(k)}(x)\R|>\varepsilon\R)<\infty.
    \end{equation*}
    By Markov's inequality, for any \(n\geq 1\), we have that
    \begin{equation}
    \label{eq: borelli phi bound}
        \mathbb{P}_{\delta_x}\F(n^{-\frac{p_1+1}{2}}\e^{-\frac{\lambda_1 n}{2}}\F|X_n[\varphi_{i,j}^{(k)}]-\e^{(\lambda_i +\mathcal{N}^*)n}W_{i,j}^{(k)}(x)\R|>\varepsilon\R) \leq \frac{\E_{\delta_x}\F[\F|X_n[\varphi_{i,j}^{(k)}]-\e^{(\lambda_i +\mathcal{N}^*)n}W_{i,j}^{(k)}(x)\R|^2\R]}{n^{p_1+1}\e^{\lambda_1 n}\varepsilon^2}.
    \end{equation}
Next, by the branching property and the martingale limit, for any \(s\geq 0\), we have that
\begin{equation*}
    X_{s+t}[\e^{-(\lambda_i+\mathcal{N})(s+t)}\varphi_{i,j}^{(k)}]=\sum_{\ell=1}^{N_s}X_t^{(\ell)}[\e^{-(\lambda_i+\mathcal{N}) s}\e^{-(\lambda_i+\mathcal{N}) t}\varphi_{i,j}^{(k)}] \rightarrow \e^{-(\lambda_i+\mathcal{N}^*) s}\sum_{\ell =1}^{N_s}W_{i,j}^{(k)}(X_s^{\ell}),
\end{equation*}
\(\mathbb{P}_{\delta_x}\)-almost surely as \(t\rightarrow \infty\), where conditionally on \(\mathcal{F}_{s}\), the \(X_t^{(\ell)}\) are independent copies of \(X_t\) initiated from \(X_s^{\ell}\), the \(\ell\)-th particle alive at time \(s\). Furthermore, conditionally on \(\mathcal{F}_{s}\), the \(W_{i,j}^{(k)}(X_s^{\ell})\) are independent. This implies, for any $s \ge 0$, 
\begin{equation}
\label{eq: Winf trick}
    W_{i,j}^{(k)}(x)=\e^{-(\lambda_i+\mathcal{N}^*) s}\sum_{i=1}^{N_s}W_{i,j}^{(k)}(X_s^{\ell}).
\end{equation}
Therefore, for any $n \ge 0$,
\begin{equation}
   X_n[\varphi_{i,j}^{(k)}]-\e^{(\lambda_i +\mathcal{N}^*)n}W_{i,j}^{(k)}(x)=\sum_{\ell=1}^{N_{n}}(\varphi_{i,j}^{(k)}(X_s^{\ell})-W_{i,j}^{(k)}(X_s^{\ell})).
   \label{eq: sum it phi}
\end{equation}
Conditionally on \(\mathcal{F}_{n}\), the right-hand side is the sum of independent random variables. Furthermore, by \ref{H1b} and \(L^2\) convergence of \(M_{i,j}^{(k)}\), we have, for any \(x \in E\), \(\E[W_{i,j}^{(k)}(x)]=\varphi_{i,j}^{(k)}(x)\), therefore these random variables are mean-zero. Moreover, by Lemma 3.1 of \cite{Moments_Chris}, for \(k\geq 2\), under (M\(k\)),
\begin{equation}
\label{eq: added for lk conv}
    \sup_{x \in E}\E_{\delta_x}\F[\F|\varphi_{i,j}^{(k)}(x)-W_{i,j}^{(k)}(x)\R|^{k}\R] < \infty.
\end{equation}
This, \eqref{eq: sum it phi}, and Lemma \ref{lemma: combin} imply
\begin{align}
\label{eq: need for markov ineq}
  \E_{\delta_x}\F[\F|X_n[\varphi_{i,j}^{(k)}]-\e^{(\lambda_i +\mathcal{N}^*)n}W
  _{i,j}^{(k)}(x)\R|^k\R] &\leq C\E_{\delta_x}[N_n^{k/2}]\sup_{x \in E}\E\F[\F|\varphi_{i,j}^{(k)}(x)-W_{i,j}^{(k)}(x)\R|^k\R] \nonumber\\
  &\leq Cn^{k(p_1-1)/2}\e^{\lambda_1 kn/2}, 
\end{align}
for some constant \(C\), where in the final inequality we have used Theorem 2.1 of \cite{Moments_Chris} to bound \(\E_{\delta_x}[N_n^{k/2}]\). This implies \(L^k\) convergence of \eqref{eq: fluct large a.s. 2} under (M\(k\)). 
To complete the proof of \(\mathbb{P}_{\delta_x}\)-a.s. convergence, using \eqref{eq: need for markov ineq} in \eqref{eq: borelli phi bound} gives 
\begin{equation*}
   \sum_{n=1}^{\infty} \mathbb{P}_{\delta_x}\F(n^{-\frac{p_1+1}{2}}\e^{-\frac{\lambda_1 n}{2}}\F|X_n[\varphi_{i,j}^{(k)}]-\e^{(\lambda_i +\mathcal{N}^*)n}W_{i,j}^{(k)}(x)\R|>\varepsilon\R)\leq C\varepsilon^{-2}\sum_{n=1}^{\infty}n^{-2}<\infty. 
\end{equation*}
\end{proof}
{Next, we present a corollary of Lemma \ref{lemma: conv of phi for fdd} which will be used to prove Theorem \ref{thm: main result crit}.
\begin{corollary}
\label{corollary: conv of phi for fdd}
   For \(1\leq i  \leq m_L\), \(1\leq j \leq p_i\), \(1\leq k \leq k_{i,j}\), and \(x \in E\), under (M\(4\)),
    \begin{equation}
    \label{eq: cor fluct large a.s. 2}
        n^{-\frac{p_1}{2}}\e^{-\frac{\lambda_1 n}{2}}\F(X_n[\varphi_{i,j}^{(k)}]-\e^{(\lambda_i +\mathcal{N}^*)n}W_{i,j}^{(k)}(x)\R) \rightarrow 0, \qquad n \to \infty,
    \end{equation}
    \(\mathbb{P}_{\delta_x}\)-almost surely, for \(n\in \mathbb{N}_{0}\).
    \end{corollary}
    \begin{proof}
        The argument is identical to the proof of Lemma \ref{lemma: conv of phi for fdd}, where a fourth moment Markov inequality is taken in \eqref{eq: borelli phi bound}.
    \end{proof}}
    {\color{black}Lemma \ref{lemma: conv of phi for fdd} gives us Theorem \ref{cor:large} and integer time \(\mathbb{P}_{\delta_x}\)-a.s.\ convergence in Theorem \ref{thm: main result large} for \(f = f_1\) defined in \eqref{eq:decomp}. We next show the same holds for \(f \in B(E)\setminus \bigoplus_{i=1}^{m_L}A_{\lambda_i}\), which, in particular, implies the result for \(f_2 \), also given in \eqref{eq:decomp}.

\begin{lemma}
\label{lemma: large interim 1}
    For \(f \in B(E)\setminus \bigoplus_{i=1}^{m_L}A_{\lambda_i}\), define
    \begin{equation}
    \label{eq: big u mcb}
        U_t = \e^{-\mathrm{Re}\lambda_{m_L} t}X_t[f].
    \end{equation}
    Fix $\delta > 0$. Then, for any non-decreasing sequence \((q_n)_{n\geq 0}\) with \(q_0>0\), and any \(x \in E\),
    \begin{equation*}
        \lim_{n\rightarrow \infty}|U_{(q_n+n)\delta}-\E[U_{(q_n+n)\delta}| \mathcal{F}_{n\delta}]| = 0, \quad \mathbb{P}_{\delta_x}-a.s.\ 
    \end{equation*}
\end{lemma}
    \begin{proof}
    The proof of this lemma is a line-by-line copy of Lemma 12.2 of \cite{bNTEbook}, thus we omit it.
    \end{proof}

The next lemma provides control over the moments of \(X[f]\), for \(f \in B(E)\setminus \bigoplus_{i=1}^{m_L}A_{\lambda_i}\).
    \medskip
\begin{lemma}
        \label{lemma: interim large 1}
            For \(k\geq 1\), assume that (M\(k\)) holds. Then, for \(1\leq \ell \leq k\), there exists a constant \(C_{\ell}\), such that
            \begin{equation}
            \label{eq: proof by induct result}
                \sup_{x \in E,f^{(1)},\dots,f^{(\ell)} \in B_1(E)\setminus \bigoplus_{i=1}^{m_L}A_{\lambda_i},t\geq 0}{\rm e}^{-\ell\mathrm{Re}\lambda_{m_L} t}(1+t)^{\ell}\F|\E_{\delta_x}\F[\prod_{i=1}^{\ell}X_t[f^{(i)}]\R]\R| \leq C_{\ell}.
            \end{equation}
        \end{lemma}
        \begin{proof}
            The proof technique is similar to Theorem 2.1 of \cite{Moments_Chris}, thus we only sketch the arguments. The case of \(k=1\) follows from \ref{H1b}. Now fix \(2\leq q \leq k\) and assume that \eqref{eq: proof by induct result} holds for all \(1\leq \ell <q\). By Lemma \ref{lem:evo-2}, for \(f^{(1)},\dots,f^{(q)} \in B(E)\) and \(x \in E\),
            \begin{equation}
            \label{eq: bound for second 1}
 \E_{\delta_x}\F[\prod_{i=1}^{q}X_t[f^{(i)}]\R]=\psi_t\F[f^{(1)}\cdots f^{(q)}\R](x)+\int_0^t \psi_s\big[\gamma \eta_{t-s}^{(q)}[f^{(1)},\dots,f^{(q)}]\big](x) \mathrm{d} s, \quad t \geq 0.
            \end{equation}
	 For the first term on the right-hand side of \eqref{eq: bound for second 1}, by \ref{H1b}, there exists a constant \(C\), such that
            \begin{equation}
                \sup_{x \in E,f^{(1)},..,f^{(q)} \in B_1(E)\setminus \bigoplus_{i=1}^{m_L}A_{\lambda_i},t\geq 0}{\rm e}^{-q\mathrm{Re}\lambda_{m_L} t}(1+t)^{q}\F|\psi_t\F[f^{(1)}\cdots f^{(q)}\R](x)\R| \leq C,
                \label{eq: control large 1}
            \end{equation}
            where we have used that \(2\mathrm{Re}\lambda_{m_L}>\lambda_1\) and \(q\geq 2\). 

            For the integral term in \eqref{eq: bound for second 1}, first consider the case of \(0\leq t \leq 2\). By the induction hypothesis, \ref{H1b}, (M\(k\)), and \ref{H2}, there exists a constant \(C\), such that
            \begin{equation}
               \sup_{x \in E,f^{(1)},\dots,f^{(q)} \in B_1(E)\setminus \bigoplus_{i=1}^{m_L}A_{\lambda_i},0\leq t \leq 2}{\rm e}^{-q\mathrm{Re}\lambda_{m_L} t}(1+t)^{q} \F|\int_0^{t} \psi_s\big[\gamma \eta_{t-s}^{(q)}[f^{(1)},\dots,f^{(q)}]\big](x) \mathrm{d} s \R|\leq C, \label{eq: control large 2}
            \end{equation}
            where we have used that the induction hypothesis implies
            \begin{equation}
             \label{eq: moment bound combintor}
                |\eta_{t}^{(q)}[f^{(1)},\dots,f^{(q)}](x)| \leq C_1\mathcal E_x \F[\sum_{\sigma \in \mathcal{P}^*([q])}\sum_{\bs i \in B_{|\sigma|,N}}1\R] \leq C_1\sup_{x \in E}\mathcal{E}_x[N^q]\leq C_2
            \end{equation} 
            for some constants \(C_1\) and \(C_2\) which are independent of the choice in the supremum of \eqref{eq: control large 2}. Thus, the integrand in \eqref{eq: control large 2} is uniformly bounded for \(0\leq t \leq 2\), \(x \in E\), \(f^{(1)},\dots,f^{(q)} \in B_1(E)\setminus \bigoplus_{i=1}^{m_L}A_{\lambda_i}\).  
            
            Now assume that \(t>2\). We split the integral over the intervals \([0,t-1]\) and \([t-1,t]\). 
            For the first interval, the induction hypothesis, \ref{H2} and (M\(k\)) allow one to control the term involving $\gamma \eta_{t-s}^{(k)}[\cdot]$ to obtain a bound of the form \eqref{eq: moment bound combintor}, and \ref{H1b} provides control of $\psi_s$ to yield the existence of a constant $C_1$ such that
         \begin{align}
               &\sup_{x \in E,f^{(1)},\dots,f^{(q)} \in B_1(E)\setminus \bigoplus_{i=1}^{m_L}A_{\lambda_i},t> 2}{\rm e}^{-q\mathrm{Re}\lambda_{m_L} t}(1+t)^{q}\F|\int_{0}^{t-1} \psi_s\left[\gamma \eta_{t-s}^{(q)}[f^{(1)},\dots,f^{(q)}] \right](x){\rm d} s\R| \nonumber \\
               &\leq C_1\sup_{t> 2} (1+t)^{q}\int_{0}^{t-1}s^{p_1-1}(t-s)^{-q}\e^{(\lambda_1 -q\mathrm{Re}\lambda_{m_L})s}\mathrm{d}s.
               \label{eq: control large 3}
	\end{align}
    Furthermore, there exists a constant \(C_2\), such that, for any  \(t> 2\), 
            \begin{align*}
                \int_0^{t-1} s^{p_1-1}(t-s)^{-q}\e^{(\lambda_1-q\mathrm{Re}\lambda_{m_L})s}\mathrm{d}s & \leq \F(\frac{t}{2}\R)^{-q} \int_0^{t/2} s^{p_1-1}\e^{(\lambda_1-q\mathrm{Re}\lambda_{m_L})s}\mathrm{d}s \\
                &\quad+\e^{(\lambda_1-q\mathrm{Re}\lambda_{m_L})\frac{t}{2}}t^{p_1-1}\int_{t/2}^{t-1} \e^{(\lambda_1-q\mathrm{Re}\lambda_{m_L})(s-\frac{t}{2})}\mathrm{d}s\\
                & \leq C_2 t^{-q}.
            \end{align*}
	Along with \eqref{eq: control large 3}, this implies that there exists a constant $C_3$ such that
            \begin{align}
       \sup_{x \in E,f^{(1)},\dots,f^{(q)} \in B_1(E)\setminus \bigoplus_{i=1}^{m_L}A_{\lambda_i},t> 2}{\rm e}^{-q\mathrm{Re}\lambda_{m_L} t}(1+t)^{q}\F|\int_{0}^{t-1} \psi_s\left[\gamma \eta_{t-s}^{(q)}[f^{(1)},\dots,f^{(q)}] \right](x){\rm d} s\R| \leq C_3. 
               \label{eq: control large 4}
            \end{align}
The interval \([t-1,t]\) is handled almost identically. The only difference being that, since \((t-s)\in [0,1]\), the lower order moments in \(\eta_{t-s}^{(q)}[f_1,\dots,f_k]\) are controlled uniformly by a constant. Therefore, there exist constants \(C_1\) and \(C_2\), such that
            \begin{align}
             &  \sup_{x \in E,f^{(1)},\dots,f^{(q)} \in B_1(E)\setminus \bigoplus_{i=1}^{m_L}A_{\lambda_i},t> 2}{\rm e}^{-q\mathrm{Re}\lambda_{m_L} t}(1+t)^{q}\F|\int_{t-1}^{t} \psi_s\left[\gamma \eta_{t-s}^{(q)}[f^{(1)},\dots,f^{(q)}] \right](x){\rm d} s\R|  \nonumber \\
              &\leq C_1\sup_{t> 2} {\rm e}^{-q\mathrm{Re}\lambda_{m_L} t}(1+t)^{q}\int_{t-1}^{t}s^{p_1-1}\e^{\lambda_1 s}\mathrm{d} s  \nonumber \\
              &\leq C_2. \label{eq: control large 5}
            \end{align} 
            The lemma follows from \eqref{eq: control large 1}, \eqref{eq: control large 2}, \eqref{eq: control large 4}, and \eqref{eq: control large 5}.
        \end{proof}

\begin{lemma}
\label{lemma: interim large 1 new}
For any \({\color{black}f \in B(E)\setminus \bigoplus_{i=1}^{m_L}A_{\lambda_i}}\), and any initial position \(x \in E\),
\begin{equation}
   \e^{-\mathrm{Re}\lambda_{m_L}n}X_{n}[f] \to 0,
\end{equation}
\(\mathbb{P}_{\delta_x}\)-almost surely as $n \to \infty$, for \(n\in \mathbb{N}_{0}\).
\end{lemma}

\medskip

\begin{proof}
Using the notation from Lemma \ref{lemma: large interim 1}, let us write
\begin{equation}\label{eq: Un}
    U_{n}=U_{n}-\E[U_{n}|\mathcal{F}_{n/2}]+\E[U_{n}|\mathcal{F}_{n/2}].
\end{equation}
By Lemma \ref{lemma: large interim 1}, $U_{n}-\E[U_{n}|\mathcal{F}_{n/2}]$ tends to zero \(\mathbb{P}_{\delta_x}\)-almost surely as \(n\rightarrow \infty\). For the third term, by the branching property,
    \begin{align}
        \E[U_{n}|\mathcal{F}_{n/2}]&= \e^{-\frac{\mathrm{Re}\lambda_{m_L}n}{2}}\sum_{i=1}^{N_{n/2}}\e^{-\frac{\mathrm{Re}\lambda_{m_L}n}{2}}\psi_s[f](X_{n/2}^i), 
        \label{eq: large conv 1}
    \end{align}
where $\{X_{n/2}^i, i = 1, \dots, N_{n/2}\}$ denotes the collection of particles alive at time $n/2 $. We show this tends to zero almost surely using a standard Borel-Cantelli argument. First note that since \(f \in B(E)\setminus \bigoplus_{i=1}^{m_L}A_{\lambda_i}\), \eqref{eq:gen-eigenfunctions} implies
\begin{equation*}
    \tilde\varphi_{i,j}^{(k)}[\psi_t[f]]=0, \quad t\geq 0, \quad 1\leq i \leq m_L, \quad 1\leq j \leq p_i, \quad 1\leq k \leq k_{i,j}.
\end{equation*}Thus, \(\psi_t[f]\) and its complex conjugate, \(\psi_t[\bar{f}]\), are in \(B(E)\setminus \bigoplus_{i=1}^{m_L}A_{\lambda_i}\). This and Lemma \ref{lemma: interim large 1} imply that there exists a constant \(C>0\) that is independent of \(n\), such that
   \begin{align}
   \nonumber
       \E_{\delta_x}\Bigg[\Big|\e^{-\frac{\mathrm{Re}\lambda_{m_L}n}{2}}\sum_{i=1}^{N_{n/2}}\e^{-\frac{\mathrm{Re}\lambda_{m_L}n}{2}}\psi_{n/2}[f](X_{n/2}^i)\Big|^{2}\Bigg] &=\e^{-2\mathrm{Re}\lambda_{m_L} n}\E_{\delta_x}\F[X_{n/2}[\psi_{n/2}[f]]X_{n/2}[\psi_{n/2}[\bar{f}]]\R]\\
       &\leq C {n}^{-4},   \label{eq: Borel-Cant large}
   \end{align}
where \ref{H1b} is used to control \(\psi_{n/2}[f]\) and \(\psi_{n/2}[\bar{f}]\). From here \eqref{eq: Borel-Cant large} and Markov's inequality imply the right-hand side of \eqref{eq: large conv 1} tends to 0 almost surely as \(n\rightarrow \infty\) by Borel-Cantelli.
\end{proof}
We are now in a position to prove Theorem \ref{thm: main result large} \(\mathbb{P}_{\delta_x}\)-a.s.\ along integer times and Theorem \ref{cor:large}. For \(f \in B(E)\), recall the decomposition given in \eqref{eq:decomp}. Lemma \ref{lemma: conv of phi for fdd} implies that, as \(t\rightarrow \infty\),
\begin{equation}
\label{eq: large result new 1 part 1}
   \e^{-\mathrm{Re}\lambda_{m_L}t}\F(X_{t}[f_1]-\sum_{i=1}^{m_L}\sum_{j=1}^{p_i}\sum_{k=1}^{k_{i,j}}\tilde\varphi_{i,j}^{(k)}[f_1]\e^{(\lambda_i+\mathcal{N}^*) t}W_{i,j}^{(k)}(x) \R)\to 0,
\end{equation}
\(\mathbb{P}_{\delta_x}\)-almost surely for \(t \in \mathbb{N}_{0}\) and in \(L^k\) under (M\(k\)). Lemmas \ref{lemma: interim large 1} and \ref{lemma: interim large 1 new} imply as \(t\rightarrow \infty\),
\begin{equation}
\label{eq: large result new 1 part 2}
      \e^{-\mathrm{Re}\lambda_{m_L}t}\F(X_{t}[f_2]-\sum_{i=1}^{m_L}\sum_{j=1}^{p_i}\sum_{k=1}^{k_{i,j}}\tilde\varphi_{i,j}^{(k)}[f_2]\e^{(\lambda_i+\mathcal{N}^*) t}W_{i,j}^{(k)}(x) \R)\to 0,  
\end{equation}
in \(L^k\) under (M\(k\)) and \(\mathbb{P}_{\delta_x}\)-almost surely for \(t \in \mathbb{N}_{0}\) respectively, where we have used that the second term is $0$ since \(f_2 \in B(E)\setminus \bigoplus_{i=1}^{m_L}A_{\lambda_i}\). We note that for the \(L^k\) convergence, in Lemma \ref{lemma: interim large 1}, we take \(f^{(1)},\dots,f^{(k/2)}=f_2\) and \(f^{(k/2)+1},\dots,f^{(k)}=\bar{f}_2\). Combining \eqref{eq: large result new 1 part 1} and \eqref{eq: large result new 1 part 2}, then performing some simple algebra gives us as \(t\rightarrow \infty\)
\begin{equation}
\label{eq: large result new 1 part 3}
   \e^{-\mathrm{Re}\lambda_{m_L}t}\F(X_{t}[f]-\sum_{i=1}^{m_L}\sum_{j=1}^{p_i}\sum_{k=1}^{k_{i,j}}\tilde\varphi_{i,j}^{(k)}[\e^{(\lambda_i+\mathcal{N}) t}f]W_{i,j}^{(k)}(x) \R)\to 0,
\end{equation}
\(\mathbb{P}_{\delta_x}\)-almost surely for \(t \in \mathbb{N}_{0}\) and in \(L^k\) under (M\(k\)).

We now extend this result to almost sure continuous time convergence. From here, we assume (M\(4\)) holds and that \(f\in C_{\psi}(E)\). For {\color{black}\(f\in B(E)\)} and \(t\geq 0\), define
\begin{align*}
    &Y_{t}^1[f]=\e^{-\mathrm{Re}\lambda_{m_L} t}X_{t}[f],\\
    &Y_{t}^2[f] = \e^{-\mathrm{Re}\lambda_{m_L} t}\sum_{i=1}^{m_L}\sum_{j=1}^{p_i}\sum_{k=1}^{k_{i,j}}\tilde\varphi_{i,j}^{(k)}[\e^{(\lambda_i+\mathcal{N)}t}f]W_{i,j}^{(k)}(x), \\
    &Y_{t}[f] = Y_{t}^1[f]-Y_{t}^2[f].
\end{align*}

We first prove the following lemma. 

\begin{lemma}\label{lem:tight}
For \(f \in C_{\psi}(E)\), there exists \(C\geq 0\), such that, for any \(n\geq 1\) and \(0\leq r\leq s\leq t\leq 1\),
    \begin{equation}
    \label{eq: large tight control 2}
        \mathbb{E}_{\delta_x}[|Y_{n+t}[f]-Y_{n+s}[f]|^{k_f}|Y_{n+s}[f]-Y_{n+r}[f]|^{k_f}] \leq C n^{-2}(t-r)^{3/2},
    \end{equation}
    where \(k_f\) is as in \eqref{eq: tight assum 1112}-\eqref{eq: tight assum 3112}.
\end{lemma}

\begin{proof}
   Let \(f_1\) and \(f_2\) be defined as in \eqref{eq:decomp}. We have that
    \begin{align*}
      &Y_{n+t}[f] - Y_{n+s}[f] = \alpha_{n,s,t} + \beta_{n,s,t},
      \end{align*}
      where
      \begin{align*}
      \alpha_{n,s,t} &= Y_{n+t}^1[f] - Y_{n+s}^1[f_2] - \E_{\delta_x}[Y_{n+t}^1[f_1]|\mathcal{F}_{n+s}],\\
      \beta_{n,s,t} &= \F(\E_{\delta_x}[Y_{n+t}^1[f_1]|\mathcal{F}_{n+s}]-Y_{n+s}^1[f_1]\R)-\F(Y_{n+t}^2[f]-Y_{n+s}^2[f]\R).
    \end{align*}
    
    Using this decomposition, it follows that there exists a constant {\color{black}\(C\geq 0\)}, such that
    \begin{align}
        \mathbb{E}_{\delta_x}[|Y_{n+t}[f]-Y_{n+s}[f]|^{k_f}|Y_{n+s}[f]-Y_{n+r}[f]|^{k_f}] &\leq C\E_{\delta_x}[|\alpha_{n,s,t}|^{k_f}|\alpha_{n,r,s}|^{k_f}]+C\E_{\delta_x}[|\beta_{n,s,t}|^{k_f}|\beta_{n,r,s}|^{k_f}]\nonumber\\
        &\quad +C\E_{\delta_x}[|\alpha_{n,s,t}|^{k_f}|\beta_{n,r,s}|^{k_f}] +C\E_{\delta_x}[|\beta_{n,s,t}|^{k_f}|\alpha_{n,r,s}|^{k_f}].\label{eq: large tight control 1}
    \end{align}   
Since \(\alpha_{n,r,s}\) is \(\mathcal{F}_{n+s}\)-measurable, we have
\begin{equation*}
    \E_{\delta_x}[|\alpha_{n,s,t}|^{k_f}|\alpha_{n,r,s}|^{k_f}] = \E_{\delta_x}[\E[|\alpha_{n,s,t}|^{k_f}|\mathcal{F}_{n+s}]|\alpha_{n,r,s}|^{k_f}].
\end{equation*}
This, \eqref{eq: large tight control 1}, and the Cauchy-Schwarz inequality imply that to show \eqref{eq: large tight control 2} it is sufficient to find a constant \(C\) such that, for any \(n\geq 1\) and $0 \le s \le t \le 1$,
\begin{align}
&\E_{\delta_x}\big[\E\big[|\alpha_{n,s,t}|^{k_f}|\mathcal{F}_{n+s}]^2]\leq Cn^{-2}(t-s)^2,\label{eq: tight control 1l}\\
&\E_{\delta_x}\big[|\alpha_{n,s,t}|^{2{k_f}}]\leq Cn^{-2}(t-s),\label{eq: tight control 2l}\\
&\E_{\delta_x}\big[|\beta_{n,s,t}|^{2{k_f}}] \leq Cn^{-2}(t-s)^2. \label{eq: tight control 3l}       
\end{align}

 We start by proving \eqref{eq: tight control 1l}. Firstly, there exists a constant \(C>0\) such that
    \begin{align}
       \E_{\delta_x}\big[|\alpha_{n,s,t}|^{k_f}|\mathcal{F}_{n+s}]^2
       &\leq C\E_{\delta_x}[|Y_{n+t}^1[f]-\E[Y_{n+t}^1[f]|\mathcal{F}_{n+s}]|^{k_f}|\mathcal{F}_{n+s}]^2\nonumber\\
       &\quad+C\F|\E_{\delta_x}[Y_{n+t}^1[f]|\mathcal{F}_{n+s}]-Y_{n+s}^1[f_2]-\E_{\delta_x}[Y_{n+t}^1[f_1]|\mathcal{F}_{n+s}]\R|^{2{k_f}},   \label{eq: cond exp large}
    \end{align}
   where we have used that the final line is \(\mathcal{F}_{n+s}\) measurable. For the final line, taking expectations, we have that
    \begin{align*}
        &\E_{\delta_x}\F[\F|\E[Y_{n+t}^1[f]|\mathcal{F}_{n+s}]-Y_{n+s}^1[f_2]-\E_{\delta_x}[Y_{n+t}^1[f_1]|\mathcal{F}_{n+s}]\R|^{2{k_f}}\R]
        =\e^{-2{k_f}\mathrm{Re}\lambda_{m_L}(n+s)}\E_{\delta_x}[|X_{n+s}\F[h_{t-s}\R]|^{2{k_f}}], 
    \end{align*}
    where 
    \begin{equation*}
        h_{t} := \e^{-\mathrm{Re}\lambda_{m_L} t}\psi_{t}[f_2]-f_2 = (\e^{-\mathrm{Re}\lambda_{m_L} t}\psi_{t}[f]-f) - (\e^{-\mathrm{Re}\lambda_{m_L} t}\psi_{t}[f_1]-f_1).
    \end{equation*}
Using the definition of $f_1$ along with \ref{H1b}, it follows that there exists a constant \(C_1\), such that
  \begin{equation}
      \sup_{0\leq t \leq 1,x \in E}t^{-1/{k_f}}\F|\e^{-\mathrm{Re}\lambda_{m}t}\psi_{t}[f_1](x)-f_1(x)\R| \leq C_1.\label{eq: result of mart expect}
  \end{equation}
  This, the fact that \(f\in C_{\psi}(E)\), and \eqref{eq: tight assum 111} imply that there exists a constant \(C_2\), such that
  \begin{equation}
  \label{eq: assum ext}
      \sup_{0  \leq t \leq 1,x \in E}t^{-1/{k_f}}|h_{t}(x)| \leq   \sup_{0\leq  t \leq 1,x \in E}t^{-1/{k_f}}|\psi_{t}[f](x)-f(x)|+C_1\leq C_2.
  \end{equation}
{Furthermore, since \(f_2 \in B_1(E)\setminus \bigoplus_{i=1}^{m_L}A_{\lambda_i}\), \ref{H1b} implies that, for all \(0\leq t \leq 1\), \(h_{t} \in B_1(E)\setminus \bigoplus_{i=1}^{m_L}A_{\lambda_i}\).} Therefore, we can apply Lemma \ref{lemma: interim large 1} to \(h_{t-s}\) which gives the existence of a constant \(C\) independent of \(s,t\) and \(n\), such that     
  \begin{equation}
  \label{eq: large assumption exp bound 12}
       \e^{-2{k_f}\mathrm{Re}\lambda_{m_L}(n+s)}\E_{\delta_x}[|X_{n+s}\F[h_{t-s}\R]|^{2{k_f}}]
        \leq Cn^{-2}(t-s)^2,
    \end{equation}
    where we have used that \(k\geq 2\) to bound by \(n^{-2}\). For the first term on the right-hand side of \eqref{eq: cond exp large}, letting $X_{t-s}^{(i)}$ denote an independent copy of $X_{t-s}$ initiated from $X_{n+s}^i$, the $i$-th particle alive at time $n+s$, we have 
    \begin{align}
      \E_{\delta_x}&[|Y_{n+t}^1[f]-\E[Y_{n+t}^1[f]|\mathcal{F}_{n+s}]|^{k_f}|\mathcal{F}_{n+s}]\nonumber \\
       &= \E_{\delta_x}\F[\F|\e^{-\mathrm{Re}\lambda_{m_L}(n+t)}\sum_{i=1}^{N_{n+s}}X_{t-s}^{(i)}[f]-\psi_{t-s}[f](X_{n+s}^i)\R|^{k_f}\M|\mathcal{F}_{n+s}\R]\nonumber \\
       &\leq C_1\e^{-{k_f}\mathrm{Re}\lambda_{m_L}(n+t)}N_{n+s}^{ {k_f}/2 }\max_{i=1,\dots,N_{n+s}}\E_{\delta_{X_{n+s}^i}}\big[\big|X_{t-s}[f]-\psi_{t-s}[f](X_{n+s}^i)\big|^{k_f}\big] \nonumber\\
       &\leq C_2\e^{-{k_f}\mathrm{Re}\lambda_{m_L}(n+t)}N_{n+s}^{ {k_f}/2 }(t-s), \label{eq: large part 2 exp bound}
    \end{align}
   where \(C_1,C_2\) are independent of \(s,t\) and \(n\). In the second line, we have used Lemma \ref{lemma: combin}, and in the final line we have used \eqref{eq: tight assum 211}. {Theorem 2.1 of \cite{Moments_Chris}} implies 
    \begin{equation}
    \label{eq: control on pop}
      \sup_{x \in E, t\geq 0}\e^{-{k_f}\lambda_1 t}(1+t)^{-{k_f}(p_1-1)}\E_{\delta_x}[N_{t}^{{k_f}}] \leq C
    \end{equation}
    for some constant \(C\). Since \(2\mathrm{Re}\lambda_{m_L}>\lambda_1\), this and \eqref{eq: large part 2 exp bound} imply the existence of constants \(C_1\) and \(C_2\), independent of \(n,s\) and \(t\), such that
   \begin{align}
        \E_{\delta_x}[\E[|Y_{n+t}^1[f]-\E[Y_{n+t}^1[f]|\mathcal{F}_{n+s}]|^{k_f}|\mathcal{F}_{n+s}]^2]
        &\leq  C_1n^{{k_f}(p_1-1)}\e^{{k_f}(\lambda-2\mathrm{Re}\lambda_{m_L})n}(t-s)^2\notag\\
        &\leq C_2 n^{-2}(t-s)^2. \label{eq: large bound 33}
   \end{align}
Combining \eqref{eq: cond exp large}, \eqref{eq: large assumption exp bound 12}, and \eqref{eq: large bound 33}, gives \eqref{eq: tight control 1l}. 
   
   The proof of \eqref{eq: tight control 2l} follows an almost identical argument. The notable difference being that in \eqref{eq: cond exp large}, the power of two in the first line appears inside the expectation instead. From this \eqref{eq: large assumption exp bound 12} is identical, and an analogous version of \eqref{eq: large bound 33} follows. However, we note that in the penultimate line of \eqref{eq: large part 2 exp bound}, the \(k_f\)-th moment is replaced by the \(2k_f\)-th moment. To control this term, we use \eqref{eq: tight assum 311} in place of \eqref{eq: tight assum 211}. 
   
   Lastly, we prove \eqref{eq: tight control 3l}. First note that, by \ref{H1b},
   \begin{equation*}
       \beta_{n,s,t} = \sum_{i=1}^{m_L}\sum_{j=1}^{p_i}\sum_{k=1}^{k_{i,j}}\beta_{n,s,t}(i,j,k),
   \end{equation*}
   where 
   \begin{equation*}
       \beta_{n,s,t}(i,j,k) = Y_{n+s}\F[\tilde\varphi_{i,j}^{(k)}\F[f\R]\F(\e^{(\lambda_i+\mathcal{N}-\mathrm{Re}\lambda_{m_L})(t-s)}-1\R)\varphi_{i,j}^{(k)}\R].
   \end{equation*} By an identical argument to \eqref{eq: Winf trick}-\eqref{eq: need for markov ineq}, we have that
   \begin{align*}
       \E_{\delta_x}\big[|\beta_{n,s,t}(i,j,k)|^{2{k_f}}]
       &\leq C^{(1)}_{i,j,k}\F|\e^{(\lambda_i+\mathcal{N}-\mathrm{Re}\lambda_{m_L})(t-s)}-1\R|^{2{k_f}}n^{{k_f}(p_1-1)}\e^{{k_f}(\lambda_1-2\mathrm{Re}\lambda_{m_L}) n}\\
       & \leq C^{(2)}_{i,j,k}n^{-2}(t-s)^2,
   \end{align*}
where \(C^{(1)}_{i,j,k}, C^{(2)}_{i,j,k}\) are constants that only depend on \(i,j,k\). Therefore, \eqref{eq: tight control 3l} holds.
\end{proof}

\begin{proof}[Proof of Theorem \ref{thm: main result large}]
Let \(f\in C_{\psi}(E)\). To complete the proof of Theorem \ref{thm: main result large}, it remains to extend the \(\mathbb{P}_{\delta_x}\)-almost sure convergence from integer times to continuous time. For this, we will show that, for any \(n\geq 1\) and \(\varepsilon>0\),
    \begin{equation}
    \label{eq: large tight proof main result 1}
        \mathbb{P}_{\delta_x}(\sup_{0\leq s,t \leq 1}|Y_{n+t}[f]-Y_{n+s}[f]|>\varepsilon) \leq C_{\varepsilon}n^{-3/2},
    \end{equation}
    for some constant \(C_{\varepsilon}\) which is independent of \(n\). Then, the result follows from Borel-Cantelli. 
    
    \medskip
    
    To this end, fix \(n\geq 1\) and \(\varepsilon>0\). We first show that
    \begin{align}
        \mathbb{P}_{\delta_x}&(\sup_{0\leq s,t \leq 1}|Y_{n+t}[f]-Y_{n+s}[f]|>\varepsilon) \notag \\
        &\leq \mathbb{P}_{\delta_x}(\sup_{0\leq r \leq s \leq t \leq 1}|Y_{n+t}[f]-Y_{n+s}[f]|\wedge |Y_{n+s}[f]-Y_{n+r}[f]| >\varepsilon/4)\notag \\
        &\hspace{4cm}+ \PP_{\delta_x}(|Y_{n+1}[f]-Y_{n}[f]|> \varepsilon/4). \label{eq: number 12}
    \end{align}
    Then we find appropriate bounds for the two terms on the right-hand side that will allow us to conclude.
    
    \medskip
    
  To prove \eqref{eq: number 12}, for \(\delta>0\), choose \(t_{\inf}^{\delta}\) and \(t_{\sup}^{\delta}\) such that 
  \begin{align*}
      \sup_{0\leq s,t \leq 1}|Y_{n+t}[f]-Y_{n+s}[f]|\leq |Y_{n+t_{\sup}^{\delta}}[f]-Y_{n+t_{\inf}^{\delta}}[f]|+\delta
  \end{align*}
  Thus,
    \begin{align}
        \sup_{0\leq s,t \leq 1}|Y_{n+t}[f]-Y_{n+s}[f]| &\leq   |Y_{n+t_{\sup}^{\delta}}[f]-Y_{n}[f]| + |Y_{n}[f]- Y_{n+t_{\inf}^{\delta}}[f]|+\delta. \label{eq: sup inq 1}
    \end{align}
    Next, we have that
    \begin{align*}
        &2\sup_{0\leq r \leq s \leq t \leq 1}|Y_{n+t}[f]-Y_{n+s}[f]|\wedge |Y_{n+s}[f]-Y_{n+r}[f]|\geq |Y_{n}[f]-Y_{n+t_{\sup}^{\delta}}[f]|\wedge |Y_{n+t_{\sup}^{\delta}}[f]-Y_{n+1}[f]|\\
        &\phantom{2\sup_{0\leq r \leq s \leq t \leq 1}|Y_{n+t}[f]-Y_{n+s}[f]|\wedge |Y_{n+s}[f]-Y_{n+r}[f]|\geq}+|Y_{n}[f]-Y_{n+t_{\inf}^{\delta}}[f]|\wedge |Y_{n+t_{\inf}^{\delta}}[f]-Y_{n+1}[f]|,
    \end{align*}
    where on the right-hand side we have taken $r = 0, t = 1$ and \(s\) equal to either \(t_{\inf}^{\delta}\) or \(t_{\sup}^{\delta}\). 
    {\color{black}
    Therefore, on the event \(\{|Y_{n+1}[f]-Y_{n}[f]|\leq \varepsilon/4\}\), we have that
       \begin{align*}
        &2\sup_{0\leq r \leq s \leq t \leq 1}|Y_{n+t}[f]-Y_{n+s}[f]|\wedge |Y_{n+s}[f]-Y_{n+r}[f]|\\
        &\qquad\geq |Y_{n}[f]-Y_{n+t_{\sup}^{\delta}}[f]|+|Y_{n}[f]-Y_{n+t_{\inf}^{\delta}}[f]|-\varepsilon/2.
    \end{align*}
    }
  Combining this with \eqref{eq: sup inq 1} implies that, on the event \(\{|Y_{n+1}[f]-Y_{n}[f]|\leq \varepsilon/4\}\), 
    \begin{equation*}
         \sup_{0\leq s,t \leq 1}|Y_{n+t}[f]-Y_{n+s}[f]| \leq 2\sup_{0\leq r \leq s \leq t \leq 1}|Y_{n+t}[f]-Y_{n+s}[f]|\wedge |Y_{n+s}[f]-Y_{n+r}[f]| + \varepsilon/2+\delta.
    \end{equation*}
    Taking \(\delta \rightarrow 0\) gives
    \begin{multline*}
        \{\sup_{0\leq s,t \leq 1}|Y_{n+t}[f]-Y_{n+s}[f]| > \varepsilon,|Y_{n+1}[f]-Y_{n}[f]|\leq \varepsilon/4\} \\
        \subseteq \{\sup_{0\leq r \leq s \leq t \leq 1}|Y_{n+t}[f]-Y_{n+s}[f]|\wedge |Y_{n+s}[f]-Y_{n+r}[f]| >\varepsilon/4,|Y_{n+1}[f]-Y_{n}[f]|\leq \varepsilon/4\}.
    \end{multline*}
   This implies \eqref{eq: number 12}. 
   
   We now find appropriate bounds for the terms on the right-hand side of \eqref{eq: number 12}. For the first term, we will use \cite[Theorem 10.3]{bBill}, which says that if there exist \(\eta>0\), \(k\geq 0\), and \(C_1\geq 0\), such that, for all \(\varepsilon>0\) and \(0\leq r \leq s \leq t \leq 1\),
    \begin{equation}
    \label{eq: cond need large t}
        \mathbb{P}_{\delta_x}(|Y_{n+t}[f]-Y_{n+s}[f]|\wedge |Y_{n+s}[f]-Y_{n+r}[f]|>\varepsilon)\leq C_1\varepsilon^{-2k}n^{-(1+\eta)}(t-r)^{1+\eta},
    \end{equation}
    then, for all \(\varepsilon>0\),
    \begin{equation}
        \mathbb{P}_{\delta_x}(\sup_{0\leq r \leq s \leq t \leq 1}|Y_{n+t}[f]-Y_{n+s}[f]|\wedge |Y_{n+s}[f]-Y_{n+r}[f]| >\varepsilon)\leq C_2\varepsilon^{-2k}n^{-(1+\eta)}, \label{eq: bill part 1}
    \end{equation}
    where \(C_2\) only depends on \(\eta\), \(k\), and \(C_1\). To show \eqref{eq: cond need large t}, first note by Markov's inequality, for any $k > 0$,
    \begin{multline*}
        \mathbb{P}_{\delta_x}(|Y_{n+t}[f]-Y_{n+s}[f]|\wedge |Y_{n+s}[f]-Y_{n+r}[f]|>\varepsilon) \\
        \leq \varepsilon^{-2k}\mathbb{E}[|Y_{n+t}[f]-Y_{n+s}[f]|^k|Y_{n+s}[f]-Y_{n+r}[f]|^k].
    \end{multline*}

By Lemma \ref{lem:tight}, taking \(k=k_f\) there exists a constant \(C_1\), such that the expectation on the right-hand side is bounded above by $C_1 n^{-3/2}(t-r)^{3/2}$. Thus, \eqref{eq: bill part 1} holds with \(\eta =1/2\).

The bound for the second term on the right-hand side of \eqref{eq: number 12} is similar. By Markov's inequality,
    \begin{align}
        \PP_{\delta_x}(|Y_{n+1}[f]-Y_{n}[f]|> \varepsilon/4) &\leq (\varepsilon/4)^{-2{k_f}}\E_{\delta_x}[|Y_{n+1}[f]-Y_{n}[f]|^{{2k_f}}]\notag \\
        & =(\varepsilon/4)^{-2{k_f}}\E_{\delta_x}[|\alpha_{n,0,1}+\beta_{n,0,1}|^{{2k_f}}]\notag \\
        & \leq   C(\varepsilon/4)^{-2{k_f}}n^{-2},\label{eq: large tight control 3}
    \end{align}
   where in the final line we have applied \eqref{eq: tight control 2l} and \eqref{eq: tight control 3l}. Combining \eqref{eq: large tight control 3}, \eqref{eq: bill part 1} and \eqref{eq: number 12} gives us \eqref{eq: large tight proof main result 1}. Thus, the theorem follows by Borel-Cantelli.
\end{proof}
\begin{proof}[Proof of Corollary \ref{cor:large}]
The triple sum in \eqref{eq: large result new 1} is dominated by terms with \(i=1\). Furthermore, by definition of \(\mathcal{N}\), for each \(f \in B(E)\), we have that
\begin{align*}
    \tilde\varphi_{i,j}^{(k)}[\mathcal{N}f]=\tilde\varphi_{i,{j+1}}^{(k)}[f], \quad 1\leq i \leq m, \quad 1\leq j \leq p_i-1,\quad 1\leq k\leq k_{i,j},\\
    \tilde\varphi_{i,p_i}^{(k)}[\mathcal{N}f]=0, \quad 1\leq i \leq m, \quad 1\leq k\leq k_{p_i,j}.
\end{align*}
By iterating this formula, one obtains that \(\varphi_{1,j}^{(k)}[\mathcal{N}^{p_1-1}f]\neq 0\) only for \(j=1\), \(1\leq k \leq k_{1,1}\), and that \(\Phi_1[\mathcal{N}^{p_1+q}f]=0\) for all \(q\geq 0\). Using this when expanding out the exponential series in the aforementioned triple sum gives the corollary.
\end{proof}
\subsection{Small regime}


To prove Theorem \ref{thm: main result small} we will use Theorem \ref{Theorem:Bill con}. Over a series of lemmas, we will prove that the three conditions of Theorem \ref{Theorem:Bill con} hold. We start with convergence of the finite dimensional distributions which is precisely Theorem \ref{cor:small}. To prove Theorem \ref{cor:small}, we use an intermediate lemma which we state after introducing some additional notation. For {\color{black}\(f\in B(E)\)}, \(t\geq 0\), \(n\geq 1\), and initial position \(x \in E\), define
\begin{align}
    &Y_{n,t}^1[f]=n^{-\frac{p_1-1}{2}}\e^{-\frac{\lambda_1(n+t)}{2}}X_{n+t}[f],\nonumber \\
    &Y_{n,t}^2[f] = n^{-\frac{p_1-1}{2}}\e^{-\frac{\lambda_1(n+t)}{2}}\sum_{i=1}^{m_L}\sum_{j=1}^{p_i}\sum_{k=1}^{k_{i,j}}\tilde\varphi_{i,j}^{(k)}[\e^{(\lambda_i+\mathcal{N)}(n+t)}f]W_{i,j}^{(k)}(x), \nonumber\\
    & Y_{n,t}[f] = Y_{n,t}^1[f]-Y_{n,t}^2[f] \label{eq: compact ver of small},
\end{align}
and for \(0\leq s \leq t <\infty\), \(n\geq 1\), define
\begin{equation}\label{eq:Y*}
    Y^*_{n,s,t}[f] = Y_{n, t}^1[f] - \mathbb E_{\delta_x}[Y_{n, t}^1[f] | \mathcal{F}_{n+s}].
\end{equation}

\medskip

    \begin{lemma}
    \label{lemma: fdd small inter}
        Assume the setting of Theorem \ref{cor:small}. Let \(k\geq 1\) and \(0\leq t_1 < \dots< t_k<\infty\). Furthermore, for each \(1\leq i \leq k\), let \(r_i\geq 1\), let \(f_1^{(1)},\dots,f_{r_1}^{(1)}\in {\color{black}B(E)\setminus \bigoplus_{i=1}^{m_L}A_{\lambda_i}}\), and \(f_1^{(i)},\dots,f_{r_i}^{(i)}\in {\color{black}B(E)}\), \(2\leq i \leq k\). Then, for any initial position \(x \in E\), jointly as \(n\rightarrow \infty\),
        \begin{align*}
            &Y^1_{n,t_1}[f_j^{(1)}]\conindis Z_j^{(1)}, \quad 1 \leq j \leq r_1,\\
            &Y^*_{n,t_{i-1},t_i}[f_j^{(i)}] \conindis Z_j^{(i)}, \quad 2 \leq i \leq k, \quad  1 \leq j \leq r_i,\\
            &Y_{n,t_k}[\varphi_{a,b}^{(c)}] \conindis Z_{a,b}^{(c)}, \quad 1\leq a \leq m_L, \quad 1 \leq b \leq p_a, \quad 1 \leq c \leq k_{a,b},
        \end{align*}
        where conditionally on \(\mathcal{W}(x)\), the \(Z_j^{(i)},Z_{a,b}^{(c)}\) are jointly complex Gaussian, the \(Z_j^{(i)}\) are {independent for different \(i\)'s}, independent of the \(Z_{a,b}^{(c)}\), and otherwise satisfy
        \begin{alignat}{3}
           \Cov(Z_{j_1}^{(1)}, Z_{j_2}^{(1)}|\mathcal{W}(x)) &= \mathcal{C}^{(3)}_{0}(f^{(1)}_{j_1}, f^{(1)}_{j_2},\mathcal{W}(x)), &&\quad1 \leq j_1, j_2\leq r_1, \nonumber\\
            \Cov(Z_{j_1}^{(i)}, Z_{j_2}^{(i)}|\mathcal{W}(x)) &= \mathcal{C}^{(2)}_{t_i-t_{i-1}}(f_{j_1}^{(i)}, f_{j_2}^{(j)},\mathcal{W}(x)), &&\quad2\leq i \leq k, \quad  1 \leq j_1, j_2\leq r_i,\notag\\
            \Cov(Z_{a_1,b_1}^{(c_1)},Z_{a_2,b_2}^{(c_2)}|\mathcal{W}(x)) &= \mathcal{C}^{(1)}_{0}(\varphi_{a_1,b_1}^{(c_1)}, \varphi_{a_2,b_2}^{(c_2)},\mathcal{W}(x)), &&\quad 1\leq a_1,a_2\leq m_L, \quad 1\leq b_1\leq p_{a_1}, \notag \\
            &&&\qquad1\leq b_2\leq p_{a_2}, \quad 1\leq c_1 \leq k_{a_1,b_1},\notag\\
            &&&\qquad 1\leq c_2\leq k_{a_2,b_2}\label{eq: cov function 1 small 1}
        \end{alignat}
 
        where we recall \(\mathcal{C}^{(1)},\mathcal{C}^{(2)},\) and \(\mathcal{C}^{(3)}\) from Section \ref{sec: main results general} 
    \end{lemma}

    {Since the proof of this lemma is rather long and technical, we first show how Theorem \ref{cor:small} follows from Lemma \ref{lemma: fdd small inter}.}

        \begin{proof}[Proof of Theorem \ref{cor:small}]
            Fix \(k\geq 1\) and \(0\leq t_1< \dots < t_k<\infty\) and note that the left-hand side of \eqref{eq: main result small} can be written as $Y_{n, t}[f] = Y_{n, t}[f_1] + Y_{n, t}[f_2] = Y_{n, t}[f_1] + Y_{n, t}^1[f_2]$, where \(f_1\) and \(f_2\) were defined in \eqref{eq:decomp}, and where the second equality follows from the fact that $f_2 \in {\color{black}B(E)\setminus \bigoplus_{i=1}^{m_L}A_{\lambda_i}}$. Next note that for each \(1 \leq j \leq k\), the branching and semigroup properties allow us to decompose \(Y_{n,t_j}[f_2]\) into the following telescoping sum
            \begin{align}
               Y_{n,t_j}[f_2] = Y_{n,t_j}^1[f_2]&=  n^{-\frac{p_1-1}{2}}\sum_{i=2}^j \e^{-\frac{\lambda_1(n+t_i)}{2}}\Big(X_{n+t_i}[\e^{-{\lambda_1(t_j-t_i)}/2}\psi_{t_j - t_i}[f_2]]\notag\\
               &\hspace{4cm}-\E[X_{n+t_i}[\e^{-{\lambda_1(t_j-t_i)}/2}\psi_{t_j - t_i}[f_2]]|\mathcal{F}_{n+t_{i-1}}]\Big)\nonumber\\
               &\qquad +n^{-\frac{p_1-1}{2}}\e^{-\frac{\lambda_1(n+t_1)}{2}}X_{n+t_1}[\e^{-\lambda_1(t_j-t_1)/2}\psi_{t_j - t_1}[f_2]]\nonumber\\
               &=\sum_{i=2}^j Y^*_{n,t_{i-1},t_i}[\e^{-{\lambda_1(t_j-t_i)}/2}\psi_{t_j - t_i}[f_2]] + Y^1_{n,t_1}[\e^{-\lambda_1(t_j-t_1)/2}\psi_{t_j - t_1}[f_2]].     \label{eq: lem eq 1}       
            \end{align}
           Similarly, by \eqref{eq:gen-eigenfunctions}, for any \(1\leq a \leq m_L\), \(1\leq b \leq p_a \), \(1\leq c \leq k_{a,b}\), we have
            \begin{align}
              Y_{n,t_j}[\tilde\varphi_{a,b}^{(c)}[f]\varphi_{a,b}^{(c)}]&=\notag-\sum_{i=1}^{k-j} Y_{n,t_{j+i-1},t_{j+i}}^*[\e^{(\frac{\lambda_1}{2}-\lambda_a-\mathcal{N})(t_{j+i}-t_j)}\tilde\varphi_{a,b}^{(c)}[f]\varphi_{a,b}^{(c)}]\\
               &\qquad +Y_{n,t_k}[\e^{(\frac{\lambda_1}{2}-\lambda_a-\mathcal{N})(t_k-t_j)}\tilde\varphi_{a,b}^{(c)}[f]\varphi_{a,b}^{(c)}].\label{eq: lem eq 1111}       
            \end{align}
            {\color{black}The summands in the final equations of both \eqref{eq: lem eq 1} and \eqref{eq: lem eq 1111} are random variables to which we may apply Lemma \ref{lemma: fdd small inter}.} Furthermore, by definition we have that \[Y_{n,t_j}[f]=Y_{n,t_j}[f_2]+\sum_{a=1}^{m_L}\sum_{b=1}^{p_a}\sum_{c=1}^{k_{a,b}}Y_{n,t_j}[\tilde\varphi_{a,b}^{(c)}[f]\varphi_{a,b}^{(c)}].\] Thus, by Lemma \ref{lemma: fdd small inter}, it follows that, for \(1\leq j \leq k\), jointly as \(n\rightarrow \infty\), 
            \begin{equation*}
                Y_{n,t_j}[f] \conindis Z^{(j)},
            \end{equation*}
            where conditionally on \(\mathcal{W}(x)\), the \(Z^{(j)}\) are jointly Gaussian. Therefore, to prove Theorem \ref{cor:small} it remains to check that conditionally on \(\mathcal{W}(x)\) the covariance matrix of \((Z^{(1)},\dots,Z^{(k)})\) is equal to the covariance matrix of \((Z_S^f(t_1),\dots,Z_S^f(t_k))\) almost surely. Using \eqref{eq: lem eq 1} and \eqref{eq: lem eq 1111}, we may write, for \(1\leq i\leq j \leq k\),
 \begin{align}
&Y_{n,t_j}[f]= Y^1_{n,t_i}[\e^{-{\lambda_1(t_j-t_i)}/2}\psi_{t_j-t_i}[f_2]]+Y_{n,t_i,t_j}^*[f_2]+Y_{n,t_j}[f_1], \nonumber\\
  &Y_{n,t_i}[f]=Y_{n,t_i}^{1}[f_2]-Y_{n,t_i,t_j}^*[f_1^{(t_j-t_i)}]+Y_{n,t_j}[f_1^{(t_j-t_i)}],\label{eq: crit identically and on}
\end{align}
where \(f_1^{(t)}\) is as in \eqref{eq:decomp}. Using this decomposition, we again can apply Lemma \ref{lemma: fdd small inter} to give jointly as \(n\rightarrow \infty\)
\begin{align*}
    &Y_{n,t_j}[f] \conindis Z_*^{(j)},\\
    &Y_{n,t_i}[f] \conindis Z_*^{(i)},
\end{align*}
where conditionally on \(\mathcal{W}(x)\), \(Z_*^{(i)}\) and \(Z_*^{(j)}\) are jointly Gaussian. By the consistency of limits, we necessarily have that
\begin{equation*}
    \Cov(Z^{(i)},Z^{(j)}|\mathcal{W}(x))=\Cov(Z_*^{(i)},Z_*^{(j)}|\mathcal{W}(x)), \quad x \in E.
\end{equation*}
From this, \eqref{eq: crit identically and on} and the covariances given in Lemma \ref{lemma: fdd small inter}, we obtain that conditionally on \(\mathcal{W}(x)\) the covariance matrix of \((Z^{(1)},\dots,Z^{(k)})\) is equal to the covariance matrix of \((Z_S^f(t_1),\dots,Z_S^f(t_k))\) almost surely. Indeed, the components involving \(\mathcal{C}^{(2)}\) and \(\mathcal{C}^{(3)}\) follow immediately and for the component involving \(\mathcal{C}^{(1)}\), it is easy to check that \(\mathcal{C}^{(1)}_0(\cdot,\cdot,\mathcal W(x))\) is bilinear, which gives us
\begin{equation*}
    \mathcal{C}^{(1)}_0(f_1^{(t_j-t_i)},f_1, \mathcal W(x)) = \mathcal{C}^{(1)}_{t_j-t_i}(f_1,f_1, \mathcal W(x)) 
\end{equation*}
as required. The joint convergence for \(f,g \in B(E)\) follows by replacing \(f\) with \(g\) in the first line of \eqref{eq: crit identically and on}.
\end{proof}
Before giving the proof of Lemma \ref{lemma: fdd small inter}, we first recall some standard facts on multivariate Gaussian random variables that follow from linear algebra. For \(\bs U\sim\mathcal{N}(0,\Sigma) \in \mathbb{R}^q\) with \(\Sigma \neq 0\), there exists a subset \(Q\subseteq \{1, \dots, q\}\) such that the covariance matrix, \(\Sigma_Q\), of \(\bs U_Q:=(U_i)_{i \in Q}\) is non-degenerate, and for \(i \notin Q\),
\begin{equation}
   \Sigma_{i,i}=\bs \Sigma_{Q,i}^T \Sigma_Q^{-1}\bs \Sigma_{Q,i}, \text{ and }\bs \Sigma_{Q,i}^T\Sigma_Q^{-1}\bs U_Q=U_i, \quad a.s, \label{eq: non-degen result small fdd}
\end{equation} 
where \(\bs \Sigma_{Q,i} = (\Sigma_{i,j})_{j\in Q}\). To make this subset unique, we say that \(\bs U\) is \(Q\)-ND, if \(Q\) satisfies the above and, for any other such \(Q^*\), \(Q_i\leq Q^*_i, \) \(1\leq i \leq |Q|\). To extend to the case where \(\Sigma=0\), we use the convention that \(\bs U_Q=0\), \(\Sigma_Q^{-1}=0\), and \(\bs \Sigma_{Q,i}=0\) for \(i \in \{1, \dots, q\}\), which ensures that \eqref{eq: non-degen result small fdd} still holds. The motivation for this definition is that the multivariate analogue of the Berry-Esseen Theorem (Theorem \ref{theorem: Berry-Esseen}) requires a non-degenerate covariance matrix, thus convergence of the non-degenerate and degenerate components in Lemma \ref{lemma: fdd small inter} must be handled separately.

\begin{proof}[Proof of Lemma \ref{lemma: fdd small inter}]
        For brevity, we give the proof in the case that the random variables and limits in Lemma \ref{lemma: fdd small inter} are real. In the complex setting, the proof is structured identically since \(\mathbb{C}\) is isomorphic to \(\mathbb{R}^2\). 
        
        Let us first introduce some notation. Fix some initial position \(x \in E\). For each \(n\geq 1\) and \(t\geq 0\), define
 \begin{align*}
	&\bs X_t^{(i)} := (X_{t}[f_1^{(i)}],\dots,X_{t}[f_{r_i}^{(i)}]), \quad  1\leq i \leq k, \\
        &\bs Y_n^{(1)} := (Y_{n,t_1}^1[f_1^{(1)}],\dots,Y_{n,t_1}^1[f_{r_1}^{(1)}]),  \\
        &\bs Y_n^{(i)} := (Y^*_{n,t_{i-1},t_{i}}[f_1^{(i)}],\dots,Y_{n,t_{i-1},t_{i}}^*[f_{r_i}^{(i)}]), \quad 2 \leq i \leq k,\\
        & \bs Y_n^{(k+1)}:=(Y_{n,t_k}[\varphi_{1,1}^{(1)}],\dots,Y_{n,t_k}[\varphi_{m_L,p_{m_L}}^{(k_{m_L,p_{m_L}})}]), \\
        & \bs Z^{(i)} := (Z_1^{(i)}, \dots, Z_{r_i}^{(i)}), \quad 1 \leq i \leq k, \\
        &\bs Z^{(k+1)}:=(Z_{1,1}^{(1)},\dots,Z_{{m_L},p_{m_L}}^{(k_{m_L,p_{m_L}})}).
    \end{align*}
    Also, let 
    \begin{equation*}
       \kappa_{n,1} := \log(n), \quad \kappa_{n,2}=n-\log(n),
    \end{equation*}
    and
            \begin{align}
      \label{eq: renorm term}
          &{\bs S}_n^{(1,1)}:=\bs Y_n^{(1)}-\E_{\delta_x}[\bs Y_n^{(1)}|\mathcal{F}_{\kappa_{n,1}}],\\
          &{\bs S}_n^{(1,2)} := \E_{\delta_x}[\bs Y_n^{(1)}|\mathcal{F}_{\kappa_{n,1}}].
      \end{align}  
Furthermore, define, for each \(n\geq 1\), the random vectors \((\bs Z^{(i)}_n,1\leq i \leq k+1)\) as follows. We extend the space in which \(X\) lives to house \(k+1\) additional vector-valued càdlàg jump processes \((\bs {\hat{Z}}^{(i)}_{t},t\geq 0,1\leq i \leq k+1)\). For each \(1\leq i \leq k+1\), \((\bs {\hat{Z}}^{(i)}_{t})_{t\geq 0}\) has jumps at the times \(n+t_{i}\), with \(t_{k+1}=t_k+1\). We let \(\boldsymbol{Z}_{n}^{(i)}:=\boldsymbol{{\hat{Z}}}_{n+t_i}^{(i)}\) which has the following law. Conditionally on \(\mathcal{F}_{\infty}\), which remains as the filtration generated only by \(X\), the \(\boldsymbol{Z}_{n}^{(i)}\) are mean-zero Gaussian random vectors, independent for different \(1\leq i\leq k+1\), different \(n\geq 2\), and otherwise have covariance matrices 
        \begin{alignat}{3}
            \Cov(Z_{n,j_1}^{(1)},Z_{n,j_2}^{(1)}|\mathcal{F}_{\infty}) 
            &= \mathcal{C}^{(3)}_{0}(f^{(1)}_{j_1},f^{(1)}_{j_2},(M_{1,1,\kappa_{n,1}}^{(q)})_{q=1}^{k_{1,1}}), &&\quad 1 \leq j_1, j_2\leq r_1, \nonumber\\
            \Cov(Z_{n,j_1}^{(i)},Z_{n,j_2}^{(i)}|\mathcal{F}_{\infty}) 
            &= \mathcal{C}^{(2)}_{t_i-t_{i-1}}(f_{j_1}^{(i)},f_{j_2}^{(i)},(M_{1,1,\kappa_{n,1}}^{(q)})_{q=1}^{k_{1,1}}), &&\quad2\leq i \leq k, \quad  1 \leq j_1, j_2\leq r_i,\nonumber\\
            \Cov(Z_{n,a_1,b_1}^{(c_1)},Z_{n,a_2,b_2}^{(c_2)}|\mathcal{F}_{\infty}) &= \mathcal{C}^{(1)}_{0}(\varphi_{a_1,b_1}^{(c_1)},\varphi_{a_2,b_2}^{(c_2)},(M_{1,1,\kappa_{n,1}}^{(q)})_{q=1}^{k_{1,1}}), &&\quad 1\leq a_1,a_2\leq m_L, \quad 1\leq b_1\leq p_{a_1}, \notag \\
            &&& \qquad 1\leq b_2\leq p_{a_2}, \quad 1\leq c_1 \leq k_{a_1,b_1},\notag\\
            &&&\qquad 1\leq c_2 \leq k_{a_2,b_2}, \label{eq: cov function 1 small}
           \end{alignat}
  where the right-hand side should be seen as the usual definitions of \(\mathcal{C}^{(1)},\mathcal{C}^{(2)}\) and \(\mathcal{C}^{(3)}\) with \(M_{1,1,\kappa_{n,1}}^{(q)}\) in place of \(W_{1,1}^{(q)}\), \(1\leq q \leq k_{1,1}\). Moreover, let \(F^{Z^{(i)}}\), \(F_n^{Z^{(i)}}\), \(F_n^{Y^{(i)}}\) and \(F_n^{S^{(1,1)}}\) be the CDFs of \(\bs Z^{(i)}\), \(\bs Z_n^{(i)}\), \(\bs Y_n^{(i)}\), and \(\bs S_n^{(1,1)}\) respectively, and let \(F^Z\), \(F_n^Z\), \(F_n^{Y_*}\) be the CDFs of 
   \begin{align*}
       &(\bs Z^{(1)},\dots ,\bs Z^{(k+1)}), \quad (\bs Z_n^{(1)},\dots ,\bs Z_n^{(k+1)}), \quad (\bs S_n^{(1,1)}, \bs Y_n^{(2)}, \dots ,\bs Y_n^{(k+1)}),    
   \end{align*}
   respectively. Finally, let \(r_{k+1}=\sum_{i=1}^{m_L}\sum_{j=1}^{p_i}k_{i,j}\), \(r = \sum_{i=1}^{k+1}r_i\). 
   
   The majority of the proof is spent showing the following convergence results. Firstly, for any \(\bs y \in (\mathbb{R}\setminus \{0\})^r\),
   \begin{equation}
      \label{eq: need to show small fdd 1}
       F_n^{Z}(\bs y) \rightarrow F^Z(\bs y), \qquad n \to \infty.
   \end{equation}
   Secondly, for \(\bs y \in (\mathbb{R}\setminus \{0\})^{r_1}\),
   \begin{equation}
      \label{eq: need to show small fdd 2}
       F_n^{S^{(1,1)}}(\bs y|\mathcal{F}_{\kappa_{n,1}}) - F^{Z^{(1)}}_n(\bs y|\mathcal{F}_{\kappa_{n,1}}) \coninprob 0,
   \end{equation}
   as \(n\rightarrow \infty\). Next, for \(2\leq i \leq k+1\), \(\bs y \in (\mathbb{R}\setminus \{0\})^{r_i}\),
   \begin{equation}
   \label{eq: need to show small fdd 3}
       F_n^{Y^{(i)}}(\bs y|\mathcal{F}_{n+t_{i-1}}) - F^{Z^{(i)}}_n(\bs y|\mathcal{F}_{n+t_{i-1}}) \coninprob 0,
   \end{equation}
  as \(n\rightarrow \infty\). Finally, 
  \begin{equation}
     \label{eq: need to show small fdd 4}
      {\bs S}_n^{(1,2)} \coninprob 0, \qquad n \to \infty.
  \end{equation}

  Assuming \eqref{eq: need to show small fdd 1}--\eqref{eq: need to show small fdd 4} hold, the proof concludes as follows. Firstly, by definition of the \(\bs Z^{(i)}_n\), \(1\leq i \leq k+1\), we have, for \(\bs y \in (\mathbb{R}\setminus \{0\})^{r_i}\), 
   \begin{equation*}
       F^{Z^{(i)}}_n(\bs y|\mathcal{F}_{n+t_{i-1}})=F^{Z^{(i)}}_n(\bs y|\mathcal{F}_{\kappa_{n,1}}).
   \end{equation*}
   Furthermore, by definition, for \(2\leq i \leq k+1\), the \(\bs Y_n^{(i)}\) are \(\mathcal{F}_{n+t_{i}}\) measurable and \(\boldsymbol{S}_n^{(1,1)}\) is \(\mathcal{F}_{n+t_{1}}\) measurable. This, \eqref{eq: need to show small fdd 2} and \eqref{eq: need to show small fdd 3} allow for repeated use of the tower property to obtain, for any \(\bs y_1 \in (\mathbb{R}\setminus \{0\})^{r_1}\),\dots, \( \bs y_{k+1} \in (\mathbb{R}\setminus \{0\})^{r_{k+1}}\),
  \begin{align}
  \label{eq: combine 113}
    &\mathbb{P}_{\delta_x}(\bs S_n^{(1,1)}\leq \bs y_1,\dots, \bs Y_n^{(k+1)}\leq  \bs y_{k+1}) - \E_{\delta_x}\big[\PP(\bs Z_n^{(1)}\leq \bs y_1|\mathcal{F}_{\kappa_{n,1}})\cdots \PP(\bs Z_n^{(k+1)}\leq  \bs y_{k+1}|\mathcal{F}_{\kappa_{n,1}})\big] \rightarrow 0,
\end{align}
as \(n\rightarrow \infty\), where to separate the probabilities associated to different \((\bs Z_n^{(i)}, 1\leq i \leq k+1)\), we use that conditionally on \(\mathcal{F}_{\kappa_{n,1}}\), the \((\bs Z^{(i)}_n,1\leq i \leq k+1)\) are independent. This and \eqref{eq: combine 113} imply that for any \(\bs y \in (\mathbb{R}\setminus \{0\})^{r}\),
\begin{equation*}
     F_n^{Y_*}(\bs y) - F^{Z}_n(\bs y) \rightarrow 0,
\end{equation*}
as \(n\rightarrow \infty\). This and \eqref{eq: need to show small fdd 1} give, for any \(\bs y \in (\mathbb{R}\setminus \{0\})^{r}\),
\begin{equation}
\label{eq: conv of CDF fdd 1}
    F_n^{Y_*}(\bs y) \rightarrow F^{Z}(\bs y),
\end{equation}
as \(n\rightarrow \infty\). In order to complete the proof, we need to extend this convergence to all $\bs y \in \mathcal B$, the set of continuity points of $F^Z$. Thus, it remains to show \eqref{eq: conv of CDF fdd 1} for \(\bs y \in \mathcal{B}\setminus (\mathbb{R}\setminus \{0\})^{r}\). To this end, let \(\varepsilon>0\) and let \(\bs y(\varepsilon)\) denote the vector whose entries are given by $\bs y(\varepsilon)_i = \varepsilon\mathbf 1_{\{y_i = 0\}}$. Then, 
\begin{equation*}
 F_n^{Y_*}(\bs y-\bs y(\varepsilon))  \leq   F_n^{Y_*}(\bs y) \leq  F_n^{Y_*}(\bs y+\bs y(\varepsilon)).
\end{equation*}
Furthermore, we can apply \eqref{eq: conv of CDF fdd 1} to the left and right terms. This and continuity of \(F^{Z}\) at \(\bs y\) imply \eqref{eq: conv of CDF fdd 1} for \(\bs y\) by sending \(\varepsilon\rightarrow 0\). This and \eqref{eq: need to show small fdd 4} imply the lemma. 

Thus, it remains to show \eqref{eq: need to show small fdd 1}-\eqref{eq: need to show small fdd 4}. The rest of the proof is split into five steps. We start by showing \eqref{eq: need to show small fdd 3} when \(i=k+1\). Then, similar ideas will allow us to prove \eqref{eq: need to show small fdd 3} for \(2\leq i \leq k\). Steps 3, 4 and 5 are then dedicated to the proofs of \eqref{eq: need to show small fdd 2}, \eqref{eq: need to show small fdd 1} and \eqref{eq: need to show small fdd 4}, respectively.

{\it Step 1: Proof of \eqref{eq: need to show small fdd 3} for $i = k+1$.} {The strategy for this step in the proof is to find upper and lower bounds for $F_n^{Y^{(k+1)}}(\bs y|\mathcal{F}_{n+t_k})$ in terms of $F_n^{Z^{(k+1)}}(\bs y|\mathcal{F}_{n+t_k}) \pm C_n$, where $C_n$ is a non-negative \(n\)-dependent random variable that tends to $0$ in probability as $n \to \infty$.}

To obtain these bounds, we will use a combination of Markov's inequality and Theorem \ref{theorem: Berry-Esseen}, which requires us to work with the covariance matrices associated with $\bs Y_n^{(k+1)}$ and $\bs Z_n^{(k+1)}$. As such, we first need to understand the behaviour of these matrices as $n \to \infty$. To this end, define \(\Sigma^{(k+1)}:=\Cov(\bs Z^{(k+1)}|\mathcal{W}(x))\). We first show that, as \(n\rightarrow \infty\),
\begin{align}
    &\Sigma_{n}^{Z^{(k+1)}}:=\Cov_x(\bs Z^{(k+1)}_n|\mathcal{F}_{n+t_k})\coninprob \Sigma^{(k+1)}, \label{eq: small fdd bound 31}\\
   &\Sigma_{n}^{Y^{(k+1)}}:=\Cov_x(\bs Y^{(k+1)}_n|\mathcal{F}_{n+t_k}) \coninprob \Sigma^{(k+1)}. \label{eq: small fdd bound 13}
\end{align}

First note that \eqref{eq: small fdd bound 31} follows from almost sure convergence of the martingales \((M_{1,1,\kappa_{n,1}}^{(q)})_{q=1}^{k_{1,1}}\) introduced in Theorem \ref{thm: main result large}, along with the continuous mapping theorem. Next, to show \eqref{eq: small fdd bound 13}, by \eqref{eq: sum it phi}, for \(1\leq a \leq m_L\), \(1 \leq b \leq p_a\), \(1 \leq c \leq k_{a,b}\),
\begin{align*}
    Y_{n,t_k}[\varphi_{a,b}^{(c)}] &=n^{-\frac{p_1-1}{2}}\e^{-\frac{\lambda_1(n+t_k)}{2}}( X_{n+t_k}[\varphi_{a,b}^{(c)}]-\e^{(\lambda_a +\mathcal{N}^*)(n+t_k)}W_{a,b}^{(c)}(x))\\
    &=n^{-\frac{p_1-1}{2}}\e^{-\frac{\lambda_1(n+t_k)}{2}}\sum_{\ell=1}^{N_{n+t_k}}(\varphi_{a,b}^{(c)}(X^{\ell}_{n+t_k})-W_{a,b}^{(c)}(X^{\ell}_{n+t_k})),
\end{align*} 
where in the first line we have used \eqref{eq: shift operator} and in the second line we note that conditionally on \(\mathcal{F}_{n+t_k}\), the right-hand side is the sum of independent random variables. Furthermore, by \ref{H1b} and the \(L^2\) convergence of the martingale \(M_{a,b}^{(c)}\), these random variables are mean-zero. This implies, for \(1\leq a_1,a_2\leq m_L\), \(1\leq b_1\leq p_{a_1}\), \(1\leq b_2\leq p_{a_2}\), \(1\leq c_1 \leq k_{a_1,b_1}\), \(1\leq c_2 \leq k_{a_2,b_2}\),
\begin{align}
    &\nonumber\Cov_x(Y_{n,t_k}[\varphi_{a_1,b_1}^{(c_1)}]Y_{n,t_k}[\varphi_{a_2,b_2}^{(c_2)}]|\mathcal{F}_{n+t_{k}}) \\
    &= n^{-(p_1-1)}\e^{-\lambda_1(n+t_k)}\sum_{\ell=1}^{N_{n+t_k}}(\E[W_{a_1,b_1}^{(c_1)}(X^{\ell}_{n+t_k})W_{a_2,b_2}^{(c_2)}(X^{\ell}_{n+t_k})]-\tilde\varphi_{a_1,b_1}^{(c_1)}(X^{\ell}_{n+t_k})\tilde\varphi_{a_2,b_2}^{(c_2)}(X^{\ell}_{n+t_k}))\nonumber\\
   &\coninprob \mathcal{C}_0^{(1)}(\varphi_{a_1,b_1}^{(c_1)},\varphi_{a_2,b_2}^{(c_2)},\mathcal{W}(x)),\label{eq: Dn conv phi}
\end{align}
as \(n\rightarrow \infty\), where the second moments of the martingale limits follow from Lemma 3.1 of \cite{Moments_Chris}, and the convergence follows from Corollary \ref{thm: SLLNs}. This is precisely \eqref{eq: small fdd bound 13}. 

The proof of the remainder of step 1 is laid out as follows. We have that $\Sigma^{(k+1)}$ is $Q$-ND for some \(Q\subseteq \{1, \dots, r_{k+1}\}\) which depends on \(\mathcal{W}(x)\). Using such \(Q\), we are going to control
\begin{align*}
&|\mathbb P_x(Y_{n, i}^{(k+1)} \le y_i, i \in Q|\mathcal{F}_{n+t_k}) - \mathbb P_x(Z_{n, i}^{(k+1)} \le y_i, i \in Q |\mathcal{F}_{n+t_k})| \\
&\qquad \le |\mathbb P_x(Y_{n, i}^{(k+1)} \le y_i, i \in Q|\mathcal{F}_{n+t_k}) - \mathbb P_x(U_{n,Q, i} \le y_i, i \in Q|\mathcal{F}_{n+t_k})|\\
& \qquad \qquad+ |\mathbb P_x(U_{n,Q, i} \le y_i, i \in Q|\mathcal{F}_{n+t_k}) - \mathbb P_x(Z_{n, i}^{(k+1)} \le y_i, i \in Q|\mathcal{F}_{n+t_k})|,
\end{align*}
 where conditionally on \(\mathcal{F}_{n+t_k}\), \(\bs U_{n, Q} \sim \mathcal{N}(0,\Sigma_{n,Q}^{Y^{(k+1)}})\) (such random variables can formally be defined in the same way as the \(\bs{Z}_n^{(i)}\)). Now, the first term on the right-hand side is of the form for use in Theorem \ref{theorem: Berry-Esseen}, where \((\Sigma_{n,Q}^{Y^{(k+1)}})^{-1/2}\) will be controlled by \(Q\)-ND of $\Sigma^{(k+1)}$ and \eqref{eq: small fdd bound 13}. By \eqref{eq: small fdd bound 31} and \eqref{eq: small fdd bound 13}, the second term is simply the difference of two Gaussian CDFs with covariance matrices converging in probability, this will be controlled using standard Gaussian bounds. The fact that \(Q\) is random due to the dependence on \(\mathcal{W}(x)\) will be overcome by simply taking an infimum over all possible \(Q\).  Finally, for the CDF evaluated at $y_i, i \notin Q$, motivated by \eqref{eq: non-degen result small fdd}, we will use $(\bs \Sigma_{n,Q,i}^{Y^{(k+1)}})^T(\Sigma_{n,Q}^{Y^{(k+1)}})^{-1}\bs Y_{n,Q}^{(k+1)}$ as a proxy for $Y_{n ,i}^{(k+1)}$ (and similarly for $ Z_{n,i}^{(k+1)}$), then control the difference between these proxies.

With the latter in mind, fix \(\bs y \in (\mathbb{R}\setminus0)^{r_{k+1}}\), \(\varepsilon>0\) and \(Q\subseteq \{1, \dots, r_{k+1}\}\). Up to but excluding equation \eqref{eq: need for fdd small 25} we shall assume \(\Sigma_{n,Q}^{Z^{(k+1)}}\) and \(\Sigma_{n,Q}^{Y^{(k+1)}}\) are non-degenerate. Note that if this is not the case, then the bounds stated up to but excluding \eqref{eq: need for fdd small 25} are invalid. This is not a concern, since we shall use the aforementioned infimum argument and the fact that \(\Sigma^{(k+1)}\) is \(Q\)-ND for some \(Q\). Second moment Markov inequalities give, for \(i \notin Q\),
\begin{align}
    &\mathbb{P}_x(|(\bs \Sigma_{n,Q,i}^{Y^{(k+1)}})^T(\Sigma_{n,Q}^{Y^{(k+1)}})^{-1}\bs Y_{n,Q}^{(k+1)}-Y^{(k+1)}_{n,i}|\geq \varepsilon |\mathcal{F}_{n+t_{k}})\nonumber\\
    &\qquad\qquad\qquad\qquad \leq \varepsilon^{-2}\Big(\Sigma_{n,(i,i)}^{Y^{(k+1)}}- (\bs \Sigma_{n,Q,i}^{Y^{(k+1)}})^T(\Sigma_{n,Q}^{Y^{(k+1)}})^{-1}\bs \Sigma_{n,Q,i}^{Y^{(k+1)}}\Big):=C_{n,Q,i,\varepsilon}^{(1)}, \label{eq: small fdd bound 1}\\
    &\mathbb{P}_x(|(\bs \Sigma_{n,Q,i}^{Z^{(k+1)}})^T(\Sigma_{n,Q}^{Z^{(k+1)}})^{-1}\bs Z_{n,Q}^{(k+1)}-Z^{(k+1)}_{n,i}|\geq \varepsilon|\mathcal{F}_{n+t_{k}})\nonumber \\
    &\qquad\qquad\qquad\qquad \leq \varepsilon^{-2}\Big(\Sigma_{n,(i,i)}^{Z^{(k+1)}}- (\bs \Sigma_{n,Q,i}^{Z^{(k+1)}})^T(\Sigma_{n,Q}^{Z^{(k+1)}})^{-1}\bs \Sigma_{n,Q,i}^{Z^{(k+1)}}\Big):=C_{n,Q,i,\varepsilon}^{(2)}, \label{eq: small fdd bound 2}\\
    &\mathbb{P}_x\Big(\Big|\Big((\bs \Sigma_{n,Q,i}^{Z^{(k+1)}})^T(\Sigma_{n,Q}^{Z^{(k+1)}})^{-1}-(\bs \Sigma_{n,Q,i}^{Y^{(k+1)}})^T(\Sigma_{n,Q}^{Y^{(k+1)}})^{-1}\Big)\bs Z_{n,Q}^{(k+1)}\Big|\geq \varepsilon\Big|\mathcal{F}_{n+t_{k}}\Big) \nonumber\\
    &\leq \varepsilon^{-2}\Big((\bs \Sigma_{n,Q,i}^{Z^{(k+1)}})^T(\Sigma_{n,Q}^{Z^{(k+1)}})^{-1}-(\bs \Sigma_{n,Q,i}^{Y^{(k+1)}})^T(\Sigma_{n,Q}^{Y^{(k+1)}})^{-1}\Big)\Sigma_{n,Q}^{Z^{(k+1)}} \nonumber\\
    &\quad \quad \quad \times\Big((\bs \Sigma_{n,Q,i}^{Z^{(k+1)}})^T(\Sigma_{n,Q}^{Z^{(k+1)}})^{-1}-(\bs \Sigma_{n,Q,i}^{Y^{(k+1)}})^T(\Sigma_{n,Q}^{Y^{(k+1)}})^{-1}\Big)^T:=C_{n,Q,i,\varepsilon}^{(3)}. \label{eq: small fdd bound 3}
\end{align}
{Then, by splitting on whether $\big|(\bs \Sigma_{n,Q,i}^{Y^{(k+1)}})^T(\Sigma_{n,Q}^{Y^{(k+1)}})^{-1}\bs Y_{n,Q}^{(k+1)}-Y^{(k+1)}_{n,i}\big|$ is smaller than or greater than $\varepsilon$ for \(i \notin Q\) and applying \eqref{eq: small fdd bound 1}, we obtain,}
\begin{align}
    &F_n^{Y^{(k+1)}}(\bs y|\mathcal{F}_{n+t_k}) \nonumber \\
    &\quad \leq \mathbb{P}_x((Y^{(k+1)}_{n,i}\leq y_i)_{i\in Q},((\bs \Sigma_{n,Q,i}^{Y^{(k+1)}})^T(\Sigma_{n,Q}^{Y^{(k+1)}})^{-1}\bs Y_{n,Q}^{(k+1)}\leq y_{i}+\varepsilon)_{i\notin Q}|\mathcal{F}_{n+t_k})+\sum_{i\notin Q}C_{n,Q,i,\varepsilon}^{(1)},\label{eq: small fdd bound 4}\\
    &F_n^{Y^{(k+1)}}(\bs y|\mathcal{F}_{n+t_k})\nonumber\\
    &\quad \geq \mathbb{P}_x((Y^{(k+1)}_{n,i}\leq y_i)_{i\in Q},((\bs \Sigma_{n,Q,i}^{Y^{(k+1)}})^T(\Sigma_{n,Q}^{Y^{(k+1)}})^{-1}\bs Y_{n,Q}^{(k+1)}\leq y_{i}-\varepsilon)_{i\notin Q}|\mathcal{F}_{n+t_k})-\sum_{i\notin Q}C_{n,Q,i,\varepsilon}^{(1)}. \label{eq: small fdd bound 5}
\end{align}
Next, we show that the right-hand sides of \eqref{eq: small fdd bound 4} and \eqref{eq: small fdd bound 5} can be upper and lower bounded respectively with the same probabilities but with $\bs U_{n,Q}$ in place of \(\bs Y_{n,Q}^{(k+1)}\), at the expense of a constant. For each \(i \notin Q\), \[\{(\bs \Sigma_{n,Q,i}^{Y^{(k+1)}})^T(\Sigma_{n,Q}^{Y^{(k+1)}})^{-1}\bs Y_{n,Q}^{(k+1)}\leq y_{i}+\varepsilon\}\] is a convex set. Therefore, by Theorem \ref{theorem: Berry-Esseen}, we have 
\begin{align}
  &\Big|\mathbb{P}_x((Y^{(k+1)}_{n,i}\leq y_i)_{i\in Q},((\bs \Sigma_{n,Q,i}^{Y^{(k+1)}})^T(\Sigma_{n,Q}^{Y^{(k+1)}})^{-1}\bs Y_{n,Q}^{(k+1)}\leq y_{i}+\varepsilon)_{i\notin Q}|\mathcal{F}_{n+t_k})\nonumber\\
  &\quad -\mathbb{P}_x((U_{n,Q,i}\leq y_i)_{i\in Q},((\bs \Sigma_{n,Q,i}^{Y^{(k+1)}})^T(\Sigma_{n,Q}^{Y^{(k+1)}})^{-1} \bs U_{n,Q}\leq y_{i}+\varepsilon)_{i\notin Q}|\mathcal{F}_{n+t_k})\Big|\nonumber\\
 &\leq C_1n^{-3(p_1-1)/2}\e^{-3\lambda_1(n+t_k)/2}\|(\Sigma_{n,Q}^{Y^{(k+1)}})^{-1/2}\|^3_2\sum_{i=1}^{N_{n+t_k}}\sup_{\substack{1\leq a \leq m_L\\ 1\leq b \leq p_a\\1\leq c \leq k_{a,b}}}\E[|\varphi_{a,b}^{(c)}(x_{\ell})-W_{a,b}^{(c)}(x_{\ell})|^3]\nonumber\\
            &\leq C_2n^{-3(p_1-1)/2}\e^{-3\lambda_1(n+t_k)/2}\|(\Sigma_{n,Q}^{Y^{(k+1)}})^{-1/2}\|^3_2 N_{n+t_k}:=C^{(4)}_{n,Q}, \label{eq: berry esseen phi}
\end{align}
 where conditionally on \(\mathcal{F}_{n+t_k}\), \(\bs U_{n, Q} \sim \mathcal{N}(0,\Sigma_{n,Q}^{Y^{(k+1)}})\). In the final line, we have used Lemma 3.1 of \cite{Moments_Chris} which gives us uniform control of the third moments over \(x\in E\) by (M\(4\)). An analogous inequality holds for the first term on the right-hand side of \eqref{eq: small fdd bound 5}. 
 
 The next step is to replace $\bs U_{n, Q}$ with $\bs Z_{n, Q}^{(k+1)}$, again at the expense of a constant. We have that
 {\begin{align}
     &\Big|\mathbb{P}_x((U_{n,Q,i}\leq y_i)_{i\in Q},((\bs \Sigma_{n,Q,i}^{Y^{(k+1)}})^T(\Sigma_{n,Q}^{Y^{(k+1)}})^{-1} \bs U_{n,Q}\leq y_{i}+\varepsilon)_{i\notin Q}|\mathcal{F}_{n+t_k})\nonumber\\
     &\quad -\mathbb{P}_x((Z_{n,i}^{(k+1)}\leq y_i)_{i\in Q},((\bs \Sigma_{n,Q,i}^{Y^{(k+1)}})^T(\Sigma_{n,Q}^{Y^{(k+1)}})^{-1} \bs Z_{n,Q}^{(k+1)}\leq y_{i}+\varepsilon)_{i\notin Q}|\mathcal{F}_{n+t_k})\Big|\nonumber\\
     &\le \Big|\mathbb{P}_x((U_{n,Q,i}\leq y_i)_{i\in Q},((\bs \Sigma_{n,Q,i}^{Y^{(k+1)}})^T(\Sigma_{n,Q}^{Y^{(k+1)}})^{-1} \bs U_{n,Q}\leq y_{i}+\varepsilon)_{i\notin Q}|\mathcal{F}_{n+t_k})\nonumber\\
     &\quad -\frac{\mathrm{det}(\Sigma_{n,Q}^{Y^{(k+1)}})^{1/2}}{\mathrm{det}(\Sigma_{n,Q}^{Z^{(k+1)}})^{1/2}}\mathbb{P}_x((U_{n,Q,i}\leq y_i)_{i\in Q},((\bs \Sigma_{n,Q,i}^{Y^{(k+1)}})^T(\Sigma_{n,Q}^{Y^{(k+1)}})^{-1} \bs U_{n,Q}\leq y_{i}+\varepsilon)_{i\notin Q}|\mathcal{F}_{n+t_k})\Big|\nonumber\\
     & + \Bigg| \frac{\mathrm{det}(\Sigma_{n,Q}^{Y^{(k+1)}})^{1/2}}{\mathrm{det}(\Sigma_{n,Q}^{Z^{(k+1)}})^{1/2}}\mathbb{P}_x((U_{n,Q,i}\leq y_i)_{i\in Q},((\bs \Sigma_{n,Q,i}^{Y^{(k+1)}})^T(\Sigma_{n,Q}^{Y^{(k+1)}})^{-1} \bs U_{n,Q}\leq y_{i}+\varepsilon)_{i\notin Q}|\mathcal{F}_{n+t_k}) \nonumber \\
     & \quad - \mathbb{P}_x((Z_{n,i}^{(k+1)}\leq y_i)_{i\in Q},((\bs \Sigma_{n,Q,i}^{Y^{(k+1)}})^T(\Sigma_{n,Q}^{Y^{(k+1)}})^{-1} \bs Z_{n,Q}^{(k+1)}\leq y_{i}+\varepsilon)_{i\notin Q}|\mathcal{F}_{n+t_k})\Bigg| \notag \\
     &\leq \Bigg |1-\frac{\mathrm{det}(\Sigma_{n,Q}^{Y^{(k+1)}})^{1/2}}{\mathrm{det}(\Sigma_{n,Q}^{Z^{(k+1)}})^{1/2}}\Bigg|\nonumber\\
     &\quad +\frac{1}{(2\pi)^{|Q|/2}\mathrm{det}(\Sigma_{n,Q}^{Z^{(k+1)}})^{1/2}}\int_{\bs v \in \mathbb{R}^{|Q|}}|\e^{-\frac{1}{2}\bs v^T (\Sigma_{n,Q}^{Z^{(k+1)}})^{-1}\bs v}-\e^{-\frac{1}{2}\bs v^T (\Sigma_{n,Q}^{Y^{(k+1)}})^{-1}\bs v}|\mathrm{d}\bs v. \label{eq: small fdd bound 8}
 \end{align}}
 To deal with the second term on the right-hand side, fix \(K\geq 0\). Then, we have
 \begin{align}
     &\frac{1}{(2\pi)^{|Q|/2}\mathrm{det}(\Sigma_{n,Q}^{Z^{(k+1)}})^{1/2}}\int_{\bs v \in \mathbb{R}^{|Q|}}|\e^{-\frac{1}{2}\bs v^T (\Sigma_{n,Q}^{Z^{(k+1)}})^{-1}\bs v}-\e^{-\frac{1}{2}\bs v^T (\Sigma_{n,Q}^{Y^{(k+1)}})^{-1}\bs v}|\mathrm{d}\bs v\nonumber\\
     &\leq \frac{1}{(2\pi)^{|Q|/2}\mathrm{det}(\Sigma_{n,Q}^{Z^{(k+1)}})^{1/2}}\int_{\bs v \in\mathbb{R}^{|Q|}\setminus \mathrm{B}(0,K)}(\e^{-\frac{1}{2}\bs v^T (\Sigma_{n,Q}^{Z^{(k+1)}})^{-1}\bs v}+\e^{-\frac{1}{2}\bs v^T (\Sigma_{n,Q}^{Y^{(k+1)}})^{-1}\bs v})\mathrm{d}\bs v\nonumber\\
     &\quad +\Big|1-\e^{K\|(\Sigma_{n,Q}^{Y^{(k+1)}})^{-1}-(\Sigma_{n,Q}^{Z^{(k+1)}})^{-1}\|_2}\Big|\nonumber\\
     &\leq \sum_{i\in Q}\F(\e^{-{K^2}/{2|Q|\Sigma_{n,(i,i)}^{Z^{(k+1)}}}}+\frac{\mathrm{det}(\Sigma_{n,Q}^{Y^{(k+1)}})^{1/2}}{\mathrm{det}(\Sigma_{n,Q}^{Z^{(k+1)}})^{1/2}}\e^{-{K^2}/{2|Q|\Sigma_{n,(i,i)}^{Y^{(k+1)}}}}\R)\nonumber\\
     &\quad +\Big|1-\e^{K\|(\Sigma_{n,Q}^{Y^{(k+1)}})^{-1}-(\Sigma_{n,Q}^{Z^{(k+1)}})^{-1}\|_2}\Big|:=C^{(5)}_{n,Q,K}, \label{eq: small fdd bound 6}
 \end{align}
where \(\mathrm{B}(0,K) :=\{\bs v \in \mathbb{R}^{|Q|}:\|\bs v\|_2\leq K\}\), and in the first line of the second inequality we have used Chernoff bounds for \(\mathbb{P}_x( U_{n,Q,i}>K|\mathcal{F}_{n+t_k})\) and \(\mathbb{P}_x( Z^{(k+1)}_{n,Q,i}>K|\mathcal{F}_{n+t_k})\), \(i \in Q\), along with that to be in the set \(\mathbb{R}^{|Q|}\setminus \mathrm{B}(0,K)\) requires at least one component to be greater than \(K/|Q|^{1/2}\). Next, using \eqref{eq: small fdd bound 6} in \eqref{eq: small fdd bound 8}, and then using \eqref{eq: small fdd bound 8} and \eqref{eq: berry esseen phi} in \eqref{eq: small fdd bound 4} and \eqref{eq: small fdd bound 5}, we obtain
\begin{align}
    &F_n^{Y^{(k+1)}}(\bs y|\mathcal{F}_{n+t_k})\nonumber\leq \mathbb{P}_x((Z^{(k+1)}_{n,i}\leq y_i)_{i\in Q},((\bs \Sigma_{n,Q,i}^{Y^{(k+1)}})^T(\Sigma_{n,Q}^{Y^{(k+1)}})^{-1}\bs Z_{n,Q}^{(k+1)}\leq y_{i}+\varepsilon)_{i\notin Q}|\mathcal{F}_{n+t_k})\\
   &\qquad \qquad\qquad\qquad +\sum_{i\notin Q}C_{n,Q,i,\varepsilon}^{(1)}+C^{(4)}_{n,Q}+C^{(6)}_{n,Q,K},\label{eq: small fdd bound 9}\\
    &F_n^{Y^{(k+1)}}(\bs y|\mathcal{F}_{n+t_k})\geq \nonumber\mathbb{P}_x((Z^{(k+1)}_{n,i}\leq y_i)_{i\in Q},((\bs \Sigma_{n,Q,i}^{Y^{(k+1)}})^T(\Sigma_{n,Q}^{Y^{(k+1)}})^{-1}\bs Z_{n,Q}^{(k+1)}\leq y_{i}-\varepsilon)_{i\notin Q}|\mathcal{F}_{n+t_k})\\
    &\qquad \qquad\qquad\qquad-\sum_{i\notin Q}C_{n,Q,i,\varepsilon}^{(1)}-C^{(4)}_{n,Q}-C^{(6)}_{n,Q,K}, \label{eq: small fdd bound 10}
\end{align}
where
 \begin{equation*}
     C^{(6)}_{n,Q,K}:= \Bigg|1-\frac{\mathrm{det}(\Sigma_{n,Q}^{Y^{(k+1)}})^{1/2}}{\mathrm{det}(\Sigma_{n,Q}^{Z^{(k+1)}})^{1/2}}\Bigg|+{C^{(5)}_{n,Q,K}}. \label{eq: small fdd bound 7}
 \end{equation*}
Next, using \eqref{eq: small fdd bound 2} and \eqref{eq: small fdd bound 3} in the first term on the right-hand side of both \eqref{eq: small fdd bound 9} and \eqref{eq: small fdd bound 10} gives
\begin{align}
    &F_n^{Y^{(k+1)}}(\bs y|\mathcal{F}_{n+t_k})\nonumber\leq \mathbb{P}_x((Z^{(k+1)}_{n,i}\leq y_i)_{i\in Q},(Z^{(k+1)}_{n,i}\leq y_{i}+3\varepsilon)_{i\notin Q}|\mathcal{F}_{n+t_k})\\
    &\qquad\qquad\qquad\qquad +\sum_{i\notin Q}(C_{n,Q,i,\varepsilon}^{(1)}+C_{n,Q,i,\varepsilon}^{(2)}+C_{n,Q,i,\varepsilon}^{(3)})+C^{(4)}_{n,Q}+C^{(6)}_{n,Q,K},\label{eq: small fdd bound 11}\\
    &F_n^{Y^{(k+1)}}(\bs y|\mathcal{F}_{n+t_k})\nonumber\geq \mathbb{P}_x((Z^{(k+1)}_{n,i}\leq y_i)_{i\in Q},(Z^{(k+1)}_{n,i}\leq y_{i}-3\varepsilon)_{i\notin Q}|\mathcal{F}_{n+t_k})\\
    &\qquad\qquad\qquad\qquad-\sum_{i\notin Q}(C_{n,Q,i,\varepsilon}^{(1)}+C_{n,Q,i,\varepsilon}^{(2)}+C_{n,Q,i,\varepsilon}^{(3)})-C^{(4)}_{n,Q}-C^{(6)}_{n,Q,K}. \label{eq: small fdd bound 12}
\end{align}
Since \(\bs y \in (\mathbb{R}\setminus 0)^{r_{k+1}}\), for any \(\bs \varepsilon:=(\varepsilon,\dots ,\varepsilon)\in \mathbb{R}^{r_{k+1}}\) sufficiently small, \(0\notin [y_i,y_i+3\varepsilon]\), \(1\leq i \leq r_{k+1}\). For such \(\varepsilon\),
\begin{align}
  | F^{Z^{(k+1)}}_n(\bs y+3\bs \varepsilon|\mathcal{F}_{n+t_{k}})- F^{Z^{(k+1)}}_n(\bs y|\mathcal{F}_{n+t_{k}})| &\leq \sum_{i=1}^{r_{k+1}}\frac{1}{(2\pi \Sigma^{Z^{(k+1)}}_{n,(i,i)})^{1/2}}\int_{y_i}^{y_i+3\varepsilon}\e^{-x^2/2\Sigma^{Z^{(k+1)}}_{n,(i,i)}}\mathrm{d}x\nonumber\\
   &\leq \sum_{i=1}^{r_{k+1}}\frac{3(2\Sigma^{Z^{(k+1)}}_{n,(i,i)})^{1/2}\varepsilon}{\pi^{1/2}}\sup_{0\leq x \leq 3\varepsilon}\F\{\frac{1}{(y_i+x)^2}\R\}:=C_{n,\varepsilon}^{(7)}. \label{eq: small fdd bound 16}
\end{align}
An identical argument gives
\begin{align}
     |F^{Z^{(k+1)}}_n(\bs y|\mathcal{F}_{n+t_{k}})- F^{Z^{(k+1)}}_n(\bs y-3\bs \varepsilon|\mathcal{F}_{n+t_{k}})| &\leq \sum_{i=1}^{r_{k+1}}\frac{1}{(2\pi \Sigma^{Z^{(k+1)}}_{n,(i,i)})^{1/2}}\int_{y_i-3\varepsilon}^{y_i}\e^{-x^2/2\Sigma^{Z^{(k+1)}}_{n,(i,i)}}\mathrm{d}x\nonumber\\
   &\leq \sum_{i=1}^{r_{k+1}}\frac{3(2\Sigma^{Z^{(k+1)}}_{n,(i,i)})^{1/2}\varepsilon}{\pi^{1/2}}\sup_{0\leq x \leq 3\varepsilon}\F\{\frac{1}{(y_i-x)^2}\R\}:=C_{n,\varepsilon}^{(8)}. \label{eq: small fdd bound 15}
\end{align}
Placing \eqref{eq: small fdd bound 16} and \eqref{eq: small fdd bound 15} into \eqref{eq: small fdd bound 11} and \eqref{eq: small fdd bound 12} gives us
\begin{align*}
    &|F_n^{Y^{(k+1)}}(\bs y|\mathcal{F}_{n+t_k})-F_n^{Z^{(k+1)}}(\bs y|\mathcal{F}_{n+t_k})|\\
    &\leq \sum_{i\notin Q}(C_{n,Q,i,\varepsilon}^{(1)}+C_{n,Q,i,\varepsilon}^{(2)}+C_{n,Q,i,\varepsilon}^{(3)})+C^{(4)}_{n,Q}+C^{(6)}_{n,Q,K}+C_{n,\varepsilon}^{(7)}+C_{n,\varepsilon}^{(8)}.
\end{align*}
Furthermore, since the choices of \(Q,\varepsilon,\) and \(K\) were arbitrary, we obtain that
\begin{align}
    &|F_n^{Y^{(k+1)}}(\bs y|\mathcal{F}_{n+t_k})-F_n^{Z^{(k+1)}}(\bs y|\mathcal{F}_{n+t_k})|\nonumber\\
    &\leq \inf_{Q\subseteq\{1, \dots, r_{k+1}\},\varepsilon,K>0}\F(\sum_{i\notin Q}(C_{n,Q,i,\varepsilon}^{(1)}+C_{n,Q,i,\varepsilon}^{(2)}+C_{n,Q,i,\varepsilon}^{(3)})+C^{(4)}_{n,Q}+C^{(6)}_{n,Q,K}+C_{n,\varepsilon}^{(7)}+C_{n,\varepsilon}^{(8)}\R). \label{eq: need for fdd small 25}
\end{align}
Finally, to obtain \eqref{eq: need to show small fdd 3} for \(i=k+1\), we show the right-hand side of \eqref{eq: need for fdd small 25} tends to 0 in probability. Let \[A^{(k+1)}=\{Q\subseteq \{1, \dots, r_{k+1}\}:\mathbb{P}_x(\Sigma^{(k+1)}\text{ is \(Q\)-ND})>0\}.\] For \(Q \in A^{(k+1)}\), if \(\Sigma^{(k+1)}\) is \(Q\)-ND, then \((\Sigma^{(k+1)}_Q)^{-1}\) is well-defined, thus the continuous mapping theorem, \eqref{eq: small fdd bound 31}, and \eqref{eq: small fdd bound 13} imply, conditionally on \(\Sigma^{(k+1)}\) being \(Q\)-ND,
\begin{equation}
\label{eq: need for fdd small 20}
    (\Sigma_{n,Q}^{Z^{(k+1)}})^{-1}\coninprob (\Sigma_{Q}^{{(k+1)}})^{-1}, \quad  (\Sigma_{n,Q}^{Y^{(k+1)}})^{-1}\coninprob (\Sigma_{Q}^{{(k+1)}})^{-1}
\end{equation}
as \(n\rightarrow \infty\). This and \eqref{eq: non-degen result small fdd} imply, for \(Q \in A^{(k+1)}\), and \(\varepsilon>0\), conditionally on \(\Sigma^{(k+1)}\) being \(Q\)-ND, 
\begin{equation}
\label{eq: need for fdd small 21}
   \sum_{i\notin Q}(C_{n,Q,i,\varepsilon}^{(1)}+C_{n,Q,i,\varepsilon}^{(2)}+C_{n,Q,i,\varepsilon}^{(3)}) \coninprob 0
\end{equation}
as \(n\rightarrow \infty\). Again by the continuous mapping theorem, conditionally on \(\Sigma^{(k+1)}\) being \(Q\)-ND, the convergence in \eqref{eq: need for fdd small 20} extends to the 2-norm, the determinant, and the inverse square root. This and \eqref{eq: need for fdd small 21} imply, for \(Q\in A^{(k+1)}\) and \(K,\varepsilon> 0\), conditionally on \(\Sigma^{Z^{(k+1)}}\) being \(Q\)-ND,
\begin{align}
  &\sum_{i\notin Q}(C_{n,Q,i,\varepsilon}^{(1)}+C_{n,Q,i,\varepsilon}^{(2)}+C_{n,Q,i,\varepsilon}^{(3)}) +C^{(4)}_{n,Q}+C^{(6)}_{n,Q,K}\nonumber\\
  &\coninprob\sum_{i\in Q}2\e^{-{K^2}/{2|Q|\Sigma_{(i,i)}^{{(k+1)}}}}\label{eq: need for fdd small 22}
\end{align}
as \(n\rightarrow \infty\), where the limit term appears from the penultimate line of \eqref{eq: small fdd bound 6}, and to obtain convergence to 0 for \(C^{(4)}_{n,Q}\), we have used that Corollary \ref{thm: SLLNs} implies
\begin{equation*}
    \e^{-\frac{3\lambda_1(n+t_k)}{2}}N_{n+t_k}\coninprob 0
\end{equation*}
as \(n\rightarrow \infty\). Since the limit of \eqref{eq: need for fdd small 22} decays exponentially in \(K\), this implies, for \(Q \in A^{(k+1)}\), and \(\varepsilon>0\), conditionally on \(\Sigma^{(k+1)}\) being \(Q\)-ND,
\begin{equation}
\label{eq: need for fdd small 23}
  \inf_{K> 0}\F(\sum_{i\notin Q}(C_{n,Q,i,\varepsilon}^{(1)}+C_{n,Q,i,\varepsilon}^{(2)}+C_{n,Q,i,\varepsilon}^{(3)}) +C^{(4)}_{n,Q}+C^{(6)}_{n,Q,K}\R)\coninprob 0  
\end{equation}
as \(n\rightarrow \infty\). Next, by \eqref{eq: small fdd bound 31}, for any \(\varepsilon>0\) sufficiently small, we have that
\begin{equation*}
 C_{n,\varepsilon}^{(7)}+C_{n,\varepsilon}^{(8)} \coninprob \sum_{i=1}^{r_{k+1}}\frac{3(2\Sigma^{{(k+1)}}_{(i,i)})^{1/2}\varepsilon}{\pi^{1/2}}\F(\sup_{0\leq x \leq 3\varepsilon}\F\{\frac{1}{(y_i+x)^2}\R\}+\sup_{0\leq x \leq 3\varepsilon}\F\{\frac{1}{(y_i-x)^2}\R\}\R)
\end{equation*}
as \(n\rightarrow \infty\). In particular, the limit is \(O(\varepsilon)\) in \(\varepsilon\), therefore, combining this with \eqref{eq: need for fdd small 23} gives us for \(Q \in A^{(k+1)}\), conditionally on \(\Sigma^{(k+1)}\) being \(Q\)-ND,
\begin{equation}
\label{eq: need for fdd small 40}
  \inf_{\varepsilon,K>0}\F(\sum_{i\notin Q}(C_{n,Q,i,\varepsilon}^{(1)}+C_{n,Q,i,\varepsilon}^{(2)}+C_{n,Q,i,\varepsilon}^{(3)})+C^{(4)}_{n,Q}+C^{(6)}_{n,Q,K}+C_{n,\varepsilon}^{(7)}+C_{n,\varepsilon}^{(8)}\R) \coninprob 0 
\end{equation} 
as \(n\rightarrow \infty\). Finally, since \(\Sigma^{{(k+1)}}\) is almost surely \(Q\)-ND for some \(Q \in A^{(k+1)}\), \eqref{eq: need for fdd small 40} implies that the right-hand side of \eqref{eq: need for fdd small 25} tends to 0 in probability as \(n\rightarrow \infty\) as required. This completes the proof of \eqref{eq: need to show small fdd 3} for \(i=k+1\).

  \medskip

%
%
%
%

{\it Step 2: Proof of \eqref{eq: need to show small fdd 3} for $2 \le i \le k$.} This follows an identical argument to Step 1, thus we omit most of the details. The only parts to clarify are the analogues of \eqref{eq: small fdd bound 13} and the third moment bound used in \eqref{eq: berry esseen phi}. 

Let us start with the analogue of \eqref{eq: small fdd bound 13}. Fix \(2\leq i \leq k\). We have that,
    \begin{equation}
        \bs Y_n^{(i)} = n^{-\frac{p_1-1}{2}}\e^{-\frac{\lambda_1 (n+t_i)}{2}}\sum_{m=1}^{N_{n+t_{i-1}}}(\bs X^{(m,i)}_{t_i-t_{i-1}}-\E_{\delta_{X_{n + t_{i-1}}^{m}}}[\bs X^{(i)}_{t_i-t_{i-1}}]), \label{eq: berry esseen use}
        \end{equation}
        where, conditionally on \(\mathcal{F}_{n+t_{i-1}}\), the \(\bs X^{(m,i)} :=  (X^{(m)}[f_1^{(i)}],\dots,X^{(m)}[f_{r_i}^{(i)}]) \) are independent, and where \(X^{(m)}\) is an independent copy of \(X\) initiated from the $m$-th particle alive at time $n + t_{i-1}$, \(X_{n + t_{i-1}}^{m}\). Therefore, conditionally on \(\mathcal{F}_{n+t_{i-1}}\), the right-hand side is a sum of mean-zero independent random variables. Let \(\Sigma^{Y^{(i)}}_n=\Cov_x(\bs Y_n^{(i)}|\mathcal{F}_{n+t_{i-1}})\). Then, for \(1\leq j,\ell \leq r_i\), we have
        \begin{align*}
       \Sigma^{Y^{(i)}}_{n,(j,\ell)} 
      &= \Cov_x(Y_{n,t_{i-1}, t_i}[f_{j}^{(i)}],Y_{n,t_{i-1}, t_i}[f_{\ell}^{(i)}]|\mathcal{F}_{n+t_{i-1}}) \\
       &= n^{-{(p_1-1)}}e^{-\lambda_1(n+t_{i-1})}\sum_{m=1}^{N_{n+t_{i-1}}}\e^{-\lambda_1(t_i-t_{i-1})}\Cov_{X_{n + t_{i-1}}^m}(X_{t_i-t_{i-1}}[f_{j}^{(i)}],X_{t_i-t_{i-1}}[f_{\ell}^{(i)}])\\
       &\coninprob \mathcal{C}^{(2)}_{t_i-t_{i-1}}(f_{j}^{(i)},f_{\ell}^{(i)},\mathcal{W}(x))
        \end{align*}
     as \(n\rightarrow \infty\), where the final line follows from Corollary \ref{thm: SLLNs}. This provides the analogue of \eqref{eq: small fdd bound 13}.
        
        The third moment control used in the equivalent of \eqref{eq: berry esseen phi} takes the form
        \begin{equation*}
            \sup_{x \in E}\E_{\delta_x}[\|{\bs X}^{(i)}_{t_i-t_{i-1}}\|_2^3]<\infty.
        \end{equation*}
        This follows from Theorem 2.1 of \cite{Moments_Chris} and (M\(4\)).

        \medskip

        {\it Step 3: Proof of \eqref{eq: need to show small fdd 2}.} This again follows the same argument as \eqref{eq: need to show small fdd 3}. Thus we omit most of the details except again clarifying the analogue of \eqref{eq: small fdd bound 13} and the third moment bound used in \eqref{eq: berry esseen phi}. Let \(\Sigma^{S^{(1,1)}}_n=\Cov_x(\bs S_n^{(1,1)}|\mathcal{F}_{\kappa_{n,1}})\). Then, for \(1\leq j,\ell\leq r_1\), we have that 
      \begin{align}
      \label{eq: cov 2}
          &\Sigma^{S^{(1,1)}}_{n,(j,\ell)}:= n^{-(p_1-1)}\e^{-\lambda_1 \kappa_{n,1}}\sum_{m=1}^{N_{\kappa_{n,1}}} \e^{-\lambda_1(\kappa_{n,2}+t_1)}\Cov_{X_{\kappa_{n,1}}^m}\big(X_{\kappa_{n,2}+t_1}[f_j^{(1)}],X_{\kappa_{n,2}+t_1}[f_{\ell}^{(1)}]\big).
      \end{align}
      Using Lemma 3.2 of \cite{Moments_Chris} to handle the second moment, and \ref{H1b} to handle the product of first moments, we have that 
      \begin{equation}
          \lim_{t\rightarrow \infty}\sup_{x \in E}|\e^{-(\lambda_1+\mathcal{N}) t}\Cov_{x}\big(X_{t}[f_j^{(1)}],X_{t}[f_{\ell}^{(1)}]\big) -L[f_j^{(1)},f_{\ell}^{(1)}](x)|=0,
          \label{eq: small covariance convergence}
      \end{equation}
      where
      \begin{equation}
        L[f,g](x)=  \Phi_{1}[fg](x) + \int_0^{\infty}\e^{-(\lambda_1+\mathcal{N})u}\Phi_{1}[V[\psi_u[f],\psi_u[g]]](x)\mathrm{d}u.
        \label{eq: use of moment limit for small}
      \end{equation}
      The following consequence of Theorem 2.2 of \cite{Moments_Chris} shall also be useful
      \begin{equation}
          \sup_{t\rightarrow \infty}\sup_{x \in E}\e^{-\lambda_1t}(1+t)^{-(p_1-1)}|\Cov_{x}\big(X_{t}[f_j^{(1)}],X_{t}[f_{\ell}^{(1)}]\big)|<\infty.
          \label{eq: small covariance convergence just bound}
      \end{equation} For \(1\leq j,\ell\leq r_1\), let \(h_{t,j,\ell}(x) =\e^{-\lambda_1 t}\Cov_{x}(X_{t}[f_j^{(1)}],X_{t}[f_{\ell}^{(1)}])\). Then, \eqref{eq: small covariance convergence just bound} implies \[\sup_{n\geq 1}\|n^{-(p_1-1)}h_{(\kappa_{n,2}+t_1),j,\ell}\|_{\infty}<\infty.\] The proofs of Lemmas \ref{lemma: large interim 1}-\ref{lemma: interim large 1 new} do not change for sequences of functions \((g_n)_{n\geq 1}\) with \(\sup_{n\geq 1}\|g_n\|<\infty\), therefore we still obtain Theorem \ref{cor:large} for such sequences giving us
      \begin{equation}
      \label{eq: almost usre conv smal}
          n^{-(p_1-1)}\e^{-\lambda_1 \kappa_{n,1}}\F(X_{\kappa_{n,1}}[h_{(\kappa_{n,2}+t_1),j,\ell}]-\sum_{a=1}^{m_L}\sum_{b=1}^{p_a}\sum_{c=1}^{k_{a,b}}\tilde\varphi_{a,b}^{(c)}[\e^{(\lambda_a+\mathcal{N})\kappa_{n,1}}h_{(\kappa_{n,2}+t_1),j,\ell}]W_{a,b}^{(c)}(x)\R)\coninprob 0
      \end{equation}
      as \(n\rightarrow \infty\). Next, we have that
      \begin{align}
      \label{eq: almost usrer conv smal}
          &n^{-(p_1-1)}\Bigg(\sum_{a=1}^{m_L}\sum_{b=1}^{p_a}\sum_{c=1}^{k_{a,b}}\tilde\varphi_{a,b}^{(c)}\F[\e^{\mathcal{N}\kappa_{n,1}-\lambda_1({\kappa_{n,2}}+t_1)+(\lambda_a-\lambda_1)\kappa_{n,1}}\Cov_{\cdot}\big(X_{\kappa_{n,2}+t_1}[f_j^{(1)}],X_{\kappa_{n,2}+t_1}[f_{\ell}^{(1)}]\big)\R]W_{a,b}^{(c)}(x)\nonumber\\
          &\qquad\qquad\qquad\qquad-\sum_{b=1}^{p_1}\sum_{c=1}^{k_{1,b}}\tilde\varphi_{1,b}^{(c)}[\e^{\mathcal{N}(n+t_1)}L[f_j^{(1)},f_{\ell}^{(1)}]]W_{1,b}^{(c)}(x)\Bigg)\coninprob 0,
      \end{align}
      where to justify this convergence we need to show there exists \(\varepsilon>0\), such that, for any \(2 \leq a \leq m_L\),
      \begin{equation}
      \label{eq: specific small bound 1}
          \sup_{t\geq 0}\e^{(\varepsilon-\lambda_1)t}\|\Phi_a[\Cov_{\cdot}\big(X_{t}[f_j^{(1)}],X_{t}[f_{\ell}^{(1)}]\big)]\|_{\infty}<\infty.
      \end{equation}
      Indeed, with this bound we have that all but the projection on \(A_{\lambda_1}\) (\(a=1\)) vanishes in \eqref{eq: almost usrer conv smal}. Therefore, since \(n^{-(p_1-1)}\e^{\mathcal{N}n}\) is a uniformly bounded operator in \(n\) on \(A_{\lambda_1}\), \eqref{eq: almost usrer conv smal} follows from \eqref{eq: small covariance convergence}. To obtain \eqref{eq: specific small bound 1}, first note the product of first moment terms vanish by \ref{H1b} since \(f^{(1)}_j,f_{\ell}^{(1)} \in B(E)\setminus \bigoplus_{i=1}^{m_L}A_{\lambda_i}\). For the second moment term, we apply Lemma \ref{lem:evo-2} and use \ref{H1b}, \ref{H4} and (M\(2\)) to obtain the following bound for \(2\leq a \leq m_L\)
      \begin{align*}
          &\|\Phi_a[\E_{\delta_{\cdot}}[X_{t}[f_j^{(1)}]X_{t}[f_{\ell}^{(1)}]]]\|_{\infty}\leq \|\Phi_a[\psi_{t}[f_j^{(1)}f_{\ell}^{(1)}]]\|_{\infty}+\int_0^{t}\|\Phi_a[\psi_s[\gamma\eta_{t-s}^{(2)}[f_j^{(1)}f_{\ell}^{(1)}]]]\|_{\infty} \mathrm{ds}\\
          &\leq C_1\F(\e^{\mathrm{Re}\lambda_2}t^{p_{\max}}+\int_0^{t}\e^{\mathrm{Re}\lambda_2s}s^{p_{\max}}\e^{(t-s)2\mathrm{Re}\lambda_{m_L+1}}(t-s)^{2p_{\max}}\mathrm{d}s\R)\\
          &\leq C_2\e^{(\mathrm{Re}\lambda_2\vee 2\mathrm{Re}\lambda_{m_L+1} )t}t^{3p_{\max}},
      \end{align*}
      where \(p_{\max}=\max\{p_1,\dots,p_m\}\). Since \((\mathrm{Re}\lambda_2\vee 2\mathrm{Re}\lambda_{m_L+1} )<\lambda_1\), \eqref{eq: specific small bound 1} holds. Finally, we have that
      \begin{equation}
      \label{eq: for the p_1}
           n^{-(p_1-1)}\F(\sum_{b=1}^{p_1}\sum_{c=1}^{k_{1,b}}\tilde\varphi_{1,b}^{(c)}[\e^{\mathcal{N}n}L[f_j^{(1)},f_{\ell}^{(1)}]]W_{1,b}^{(c)}(x)\R)\rightarrow \frac{1}{(p_1-1)!}\sum_{c=1}^{k_{1,1}}\tilde\varphi_{1,1}^{(c)}[\mathcal{N}^{p_1-1}L[f_j^{(1)},f_{\ell}^{(1)}]]W_{1,b}^{(c)}(x)
      \end{equation}
      This, \eqref{eq: almost usre conv smal} and \eqref{eq: almost usrer conv smal} give
      \begin{equation*}
          \Sigma^{S^{(1,1)}}_{n,(j,\ell)}\coninprob \mathcal{C}_0^{(3)}(f_{j}^{(1)},f_{\ell}^{(1)},\mathcal{W}(x))
      \end{equation*}
      as \(n\rightarrow \infty\). This completes the analogue of \eqref{eq: small fdd bound 13}. 
      
      The third moment control used in the analogue of \eqref{eq: berry esseen phi} takes the form
      \begin{equation*}
          \sup_{n\geq 1}n^{-3(p_1-1)/2}\e^{-3\lambda_1(\kappa_{n,2}+t_1)/2}\sup_{x\in E}\E_{\delta_x}[\|\bs X^{(1)}_{\kappa_{n,2}+t_1}\|_2^3] < \infty,
      \end{equation*}
      which follows from Theorem 2.2 of \cite{Moments_Chris} and (M\(4\)).

     {\it Step 4: Proof of \eqref{eq: need to show small fdd 1}.} It is sufficient to show, for any \(\bs y \in (\mathbb{R}\setminus\{0\})^r\),
     \begin{equation*}
         F_n^Z(\bs y|\mathcal{F}_{\kappa_{n,1}})-F^Z(\bs y|\mathcal{W}(x)) \coninprob 0
     \end{equation*} 
     as \(n\rightarrow \infty\). Again, the proof of this result follows an identical argument to the proof of \eqref{eq: need to show small fdd 3}. The convergence of the conditional covariance matrix of \(\bs Z_n\) is given by \eqref{eq: small fdd bound 31}, and the left-hand side of the analogue of \eqref{eq: berry esseen phi} is equal to zero, since \(\bs Z_n\) is already conditionally Gaussian. We omit the details.  
     
     {\it Step 5: Proof of \eqref{eq: need to show small fdd 4}.} Since the functions \(f^{(1)}_1,\dots,f_{r_1}^{(1)} \in B(E)\setminus \bigoplus_{i=1}^{m_L}A_{\lambda_i}\), and thus by \eqref{eq:gen-eigenfunctions} \(\psi_{\kappa_{n,2}+t_1}[f_{r_1}^{(1)}],\dots,\psi_{\kappa_{n,2}+t_1}[f_1^{(1)}]\in  B(E)\setminus \bigoplus_{i=1}^{m_L}A_{\lambda_i}\), we can apply Theorem 2.2 of \cite{Moments_Chris} and \ref{H1b} to obtain 
     \begin{align}
        \E_{\delta_x}\F[\F\|\boldsymbol{S}_{n}^{(1,2)}\R\|^2_2\R]  &=n^{-(p_1-1)}\e^{-\lambda_1(n+t)}\sum_{j=1}^{r_1}\E_{\delta_x}\F[|X_{\kappa_{n,1}}[\psi_{\kappa_{n,2}+t_1}[f_j^{(1)}]]|^2\R]\nonumber\\
         &\leq C_2 n^{-(p_1-1)}\e^{\lambda_1\kappa_{n,1}+2\mathrm{Re}\lambda_{m_L+1} \kappa_{n,2}-\lambda_1n}\kappa_{n,2}^{2(p_{\max}-1)} \kappa_{n,1}^{(p_1-1)}\rightarrow 0 \label{eq: small regime last condition for fdd}
     \end{align}
    as \(n\rightarrow \infty\), where the convergence follows from the fact that \(2\mathrm{Re}\lambda_{m_L+1}<\lambda_1\). This along with a second moment Markov inequality imply \eqref{eq: need to show small fdd 4}.
\end{proof}
        \begin{lemma}

        \label{lemma: tightness small}
             Under the assumptions of Theorem \ref{thm: main result small}, for each \(f\in C_{\psi}(E)\) and \(T\geq 0\), there exists a constant \(C_{f,T}>0\), such that, for any \(0\leq r \leq s \leq t \leq T\), and \(n\geq 1\), 
            \begin{equation}
              \E_{\delta_x}\big[|Y_{n,t}[f]-Y_{n,s}[f]|^{k_f}|Y_{n,s}[f]-Y_{n,r}[f]|^{k_f}\big] \leq C_{f,T}(t-r)^{3/2}. \label{eq: eq bound 1}
            \end{equation}
        \end{lemma}
        
        \smallskip
        \begin{proof}
            The proof of Lemma \ref{lemma: tightness small} is almost a line-by-line copy of Lemma \ref{lem:tight}. Thus, we omit the majority of the proof, but clarify two minor details for which the proofs differ. Firstly, the sufficient bounds presented in \eqref{eq: tight control 1l}-\eqref{eq: tight control 3l} now take the following form. There exists constant \(C_1>0\), such that, for any \(n\geq 1\), \(0\leq s\leq t \leq T\),
            \begin{align*}
&\E_{\delta_x}\big[\E\big[|\alpha_{n,s,t}|^{k_f}|\mathcal{F}_{n+s}]^2]\leq C_1(t-s)^2,\\
&\E_{\delta_x}\big[|\alpha_{n,s,t}|^{2{k_f}}]\leq C_1(t-s),\\
&\E_{\delta_x}\big[|\beta_{n,s,t}|^{2{k_f}}] \leq C_1(t-s)^2,     
\end{align*}
where
      \begin{align*}
      \alpha_{n,s,t} &= Y_{n,t}^1[f] - Y_{n,s}^1[f_2] - \E_{\delta_x}[Y_{n,t}^1[f_1]|\mathcal{F}_{n+s}],\\
      \beta_{n,s,t} &= \F(\E_{\delta_x}[Y_{n,t}^1[f_1]|\mathcal{F}_{n+s}]-Y_{n,s}^1[f_1]\R)-\F(Y_{n,t}^2[f]-Y_{n,s}^2[f]\R).
    \end{align*}
Secondly, the analogue of \eqref{eq: large assumption exp bound 12} is now a consequence of Theorem 2.2 of \cite{Moments_Chris}, where we use that \(h_{t-s} \in B(E)\setminus \bigoplus_{i=1}^{m_L}A_{\lambda_i}\) with \(2\mathrm{Re}\lambda_{m_L+1}<\lambda_1\). This gives the analogue bound
\begin{equation*}
    n^{-{k_f}(p_1-1)}\e^{-{k_f}\lambda_1(n+s)}\E_{\delta_x}[|X_{n+s}\F[h_{t-s}\R]|^{2{k_f}}]
        \leq C(t-s)^2.
\end{equation*}
The remainder of the proof can be constructed by following the proof of Lemma \ref{lem:tight} line-by-line. 
        \end{proof}

    \medskip
    
\begin{proof}[Proof of Theorem \ref{thm: main result small}]
Let \(f \in C_{\psi}(E)\), and, for \(n\geq 1\), \(t\geq 0\), let \(Y_{n,t}[f]\) be as in \eqref{eq: compact ver of small}. We will apply Theorem \ref{Theorem:Bill con} to \((Y_{n,\cdot}[f])_{n\geq 1}\) with limit \(Z_S^f\). Conditions 1 (convergence of finite dimensional distributions) and 3 (tightness) of the aforementioned theorem have already been shown in Theorem \ref{cor:small} and Lemma \ref{lemma: tightness small} respectively. Therefore, it remains to check the second condition. We have, for \(x \in E\), and \(t\geq\delta \geq 0\), 
\begin{align}
   \Var(Z_S^f(t)-Z_S^f(t-\delta)|\mathcal{W}(x)) &= 2\mathcal{C}^{(1)}_{0}(f_1,f_1,\mathcal{W}(x))-2\mathcal{C}^{(1)}_{\delta}(f_1,f_1,\mathcal{W}(x))\nonumber\\
   &+2\mathcal{C}^{(2)}_{\delta}(f_1^{(\delta)},f_2,\mathcal{W}(x))\nonumber\\
   &+2\mathcal{C}^{(3)}_{0}(f_2,f_2,\mathcal{W}(x))-2\mathcal{C}^{(3)}_{\delta}(f_2,f_2,\mathcal{W}(x)). \label{eq: continuity req}
\end{align}
To prove the second condition of Theorem \ref{Theorem:Bill con}, it is sufficient to show that the right-hand side is continuous in \(\delta\) at 0 almost surely. Starting with the third line, we have that
\begin{align*}
    &\mathcal{C}^{(3)}_{\delta}(f_2,f_2,\mathcal{W}(x))-\mathcal{C}^{(3)}_0(f_2,f_2,\mathcal{W}(x)) \\
    &\quad = \frac{1}{(p_1-1)!}\sum_{i=1}^{k_{1,1}}W_{1,1}^{(i)}(x)\Bigg(\tilde\varphi_{1,1}^{(i)}[\mathcal{N}^{p_1-1}f_2(\e^{-\lambda_1\delta/2} \psi_{\delta}[f_2]-f_2)]\\
    & \quad\quad \quad \qquad \qquad \qquad \qquad +\int_0^{\infty}\e^{-\lambda_1 u}\tilde\varphi_{1,1}^{(i)}\F[\mathcal{N}^{p_1-1}V[\psi_u[f_2],
    \psi_u[\e^{-\lambda_1\delta/2}\psi_\delta[f_2] - f_2]
     \R]\mathrm{d}u\Bigg).
\end{align*} 
{Therefore, almost sure continuity at \(\delta=0\) is satisfied if we can show
\begin{align}
   &\lim_{\delta\rightarrow 0}\sup_{x \in E}|\e^{-\lambda_1 \delta/2}\psi_\delta[f_2](x)-f_2(x)|=0.\label{eq: small billing prop 2 1}
\end{align}}
To see this, note that
\begin{multline}
\label{eq: argument similar 12}
    \lim_{\delta\rightarrow 0}\sup_{x \in E}|{\e^{-\lambda_1 \delta/2}}\psi_\delta[f_2](x)-f_2(x)|\\ \leq \lim_{\delta\rightarrow 0}\sup_{x \in E}(|{\e^{-\lambda_1 \delta/2}}\psi_\delta[f](x)-f(x)|+|{\e^{-\lambda_1 \delta/2}}\psi_\delta[f_1](x)-f_1(x)|)=0.
\end{multline} 
The first term converges to $0$ due to \eqref{eq: tight assum 111}. The second term converges to $0$ due to \eqref{eq:gen-eigenfunctions} since \(f_1\) is a linear combination of eigenfunctions.

For the first line of \eqref{eq: continuity req}, a similar argument shows almost sure continuity at \(\delta=0\) is satisfied if
\begin{equation*}
    \lim_{\delta \rightarrow 0}\sup_{1\leq i \leq m_L}\|(\e^{(\frac{\lambda_1}{2}-\lambda_i-\mathcal{N})\delta}-1)\Phi_i[f]\|_{\infty } =0
\end{equation*}
which trivially holds by the exponential series. 

Finally, for the second line of \eqref{eq: continuity req}, we first note that by Lemma \ref{lem:evo-2},
\begin{align}
    \Cov_x(X_{\delta}[f_1^{(\delta)}],X_{\delta}[f_2]) &= \psi_{\delta}[f_1^{(\delta)}f_2](x)+ \int_0^{\delta}\psi_s[V[\psi_{\delta-s}[f_1^{(\delta)}],\psi_{\delta-s}[f_2]]](x)\mathrm{d}s \nonumber \\
    &\quad -\psi_{\delta}[f_1^{(\delta)}](x)\psi_{\delta}[f_2](x). \label{eq: tightness bill cond 2 1}
\end{align}
By \ref{H1b} and (M\(2\)), the integrand is uniformly bounded in \([0,\delta]\), thus the integral term is continuous in \(\delta\). For the first and third terms, first note that, by definition, \(\lim_{\delta\rightarrow 0}\|f_1-f^{(\delta)}_1\|_{\infty}=0\). This and \eqref{eq: argument similar 12} imply 
\begin{align}
\label{eq: conv of both 1}
    &\lim_{\delta\rightarrow 0}\sup_{x\in E}|\psi_{\delta}[f_1^{(\delta)}](x)\psi_{\delta}[f_2](x) - f_1(x)f_2(x)|=0 \\
    &\lim_{\delta\rightarrow 0}\sup_{x\in E}|\psi_{\delta}[f_1^{(\delta)}f_2](x) - \psi_{\delta}[f_1f_2](x)|=0
\end{align}
These two limits along with  \eqref{eq: tightness bill cond 2 1} imply 
\begin{equation*}
    \mathcal{C}^{(2)}_{\delta}(f_1^{(\delta)},f_2,\mathcal{W}(x))= \frac{\e^{-{\lambda_1}\delta}}{(p_1-1)!}\sum_{i=1}^{k_{1,1}}W_{1,1}^{(i)}(x)\tilde\varphi_{1,1}^{(i)}[\mathcal{N}^{p_1-1}(\psi_{\delta}[f_1f_2]-f_1f_2+R_{\delta})],
\end{equation*}
where \(\lim_{\delta\rightarrow 0}\|R_{\delta}\|_{\infty}=0\). Since the \(\tilde\varphi_{i,j}^{(k)}\) are bounded with respect to \(\|\cdot \|_{\infty}\), we have that
\begin{equation*}
    \lim_{\delta\rightarrow 0}\sup_{1\leq i\leq k_{1,1}}\tilde\varphi_{1,1}^{(i)}[\mathcal{N}^{p_1-1}R_{\delta}]=\lim_{\delta\rightarrow 0}\sup_{1\leq i\leq k_{1,p_1}}\tilde\varphi_{1,p_1}^{(i)}[R_{\delta}]=0.
\end{equation*}
Furthermore, by \eqref{eq:gen-eigenfunctions}, we have that
\begin{equation*}
    \lim_{\delta\rightarrow 0}\sup_{1\leq i\leq k_{1,1}}\tilde\varphi_{1,1}^{(i)}[\mathcal{N}^{p_1-1}(\psi_{\delta}[f_1f_2]-f_1f_2)]=\lim_{\delta\rightarrow 0}\sup_{1\leq i\leq k_{1,p_1}}\tilde\varphi_{1,p_1}^{(i)}[\psi_{\delta}[f_1f_2]-f_1f_2]=0.
\end{equation*}Therefore, we have almost sure continuity of the second line of \eqref{eq: continuity req} at \(\delta =0\). Hence, the second condition of Theorem \ref{Theorem:Bill con} holds and so Theorem \ref{Theorem:Bill con} gives the required functional convergence.
\end{proof}

\subsection{Critical regime}
Finally, we prove Theorem \ref{thm: main result crit}. The results and proofs in this section are almost identical to those presented for the small regime. We thus only sketch the arguments, providing details where they deviate from their analogues in the previous section. Again, to prove Theorem \ref{thm: main result crit}, we show that the three conditions of Theorem \ref{Theorem:Bill con} hold, starting with the convergence of finite dimensional distributions which is precisely Theorem \ref{cor:crit}.

To prove Theorem \ref{cor:crit}, we use an intermediate lemma which we state after introducing some additional notation. For {\color{black}\(f,g\in\mathrm{Ei}(\Lambda_C)\)}, \(t\geq 0\), \(n\geq 1\), and initial position \(x \in E\), define
\begin{align}
    &Y_{n,t}^1[f, g]=\e^{-{\lambda_f nt}}n^{-\frac{2p_f+p_1-2}{2}}X_{nt}[g], \qquad Y_{n,t}^1[f] = Y_{n,t}^1[f, f],\nonumber \\
    &Y_{n,t}^2[f, g] = \e^{-{\lambda_f nt}}n^{-\frac{2p_f+p_1-2}{2}}\sum_{i=1}^{m_L}\sum_{j=1}^{p_i}\sum_{k=1}^{k_{i,j}}\tilde\varphi_{i,j}^{(k)}[\e^{(\lambda_i+\mathcal{N)}nt}g]W_{i,j}^{(k)}(x), \qquad Y_{n,t}^2[f] = Y_{n,t}^2[f, f],\nonumber\\
    & Y_{n,t}[f, g] = Y_{n,t}^1[f, g]-Y_{n,t}^2[f, g], \qquad Y_{n,t}[f] = Y_{n,t}[f, f] \label{eq: compact ver of crit},
\end{align}
and, for \(0\leq s \leq t <\infty\), \(n\geq 1\), define
\begin{equation}\label{eq:Y* crit}
    Y^*_{n,s,t}[f] = Y_{n, t}^1[f] - \mathbb E[Y_{n, t}^1[f] | \mathcal F_{ns}].
\end{equation}
    \begin{lemma}
    \label{lemma: fdd small inter crit}
        Take \(k\geq 1\) and \(0< t_1 <\dots < t_k<\infty\). Furthermore, for each \(1\leq i \leq k\), let \(r_i\geq 1\), and let \(f_1^{(i)},\dots,f_{r_i}^{(i)}\in \mathrm{Ei}(\Lambda_C)\setminus \bigoplus_{i=1}^{m_L}A_{\lambda_i}\). Then, jointly as \(n\rightarrow \infty\),
        \begin{align*}
            &Y_{n,t_1}^1[f_j^{(1)}]\conindis Z_j^{(1)}, \quad 1 \leq j \leq r_1,\\
            &Y_{n,t_{i-1},t_i}^*[f_j^{(i)}] \conindis Z_j^{(i)}, \quad 2 \leq i \leq k, 1 \leq j \leq r_i,
        \end{align*}
        where conditionally on \(\mathcal{W}(x)\), \(Z_{j_1}^{(i_1)},Z_{j_2}^{(i_2)}\) are mean-zero Gaussian random variables, are independent for \(i_1 \neq i_2\), and otherwise satisfy, for $1 \leq j, \ell\leq r_1$,
        \begin{multline}
          \Cov(Z_{j}^{(1)},Z_{\ell}^{(1)}|\mathcal{W}(x))\\
          =\bs{1}_{\lambda_{{1, j}}=\bar\lambda_{{1, \ell}}}B(0,t_1,f_j^{({ 1})},f_{\ell}^{({ 1})})\sum_{q=1}^{k_{1,1}}\tilde\varphi_{1,1}^{(q)}[\mathcal{N}^{p_1-1}V[\mathcal{N}^{p_{1, j}-1}\Phi _{\nu_{1,j}}[f_j^{(1)}],\mathcal{N}^{p_{1, \ell}-1}\Phi _{\nu_{1,\ell}}[f_{\ell}^{(1)}]]]W_{1,1}^{(q)}(x), \label{eq: cov function 1 critical 1}
          \end{multline}
        and for $2\leq i \leq k$, $1 \leq j, \ell\leq r_i$,
        \begin{multline}
            \Cov(Z_{j}^{(i)}, Z_{\ell}^{(i)}|\mathcal{W}(x)) \\
            = \bs{1}_{\lambda_{i, j}=\bar\lambda_{i, \ell}}B(t_{i-1},t_i,f_j^{(i)},f_{\ell}^{(i)})\sum_{q=1}^{k_{1,1}}\tilde\varphi_{1,1}^{(q)}[\mathcal{N}^{p_1-1}V[\mathcal{N}^{p_{i, j}-1}\Phi _{\nu_{i,j}}[f_j^{(i)}],\mathcal{N}^{p_{{i, \ell}}-1}\Phi _{\nu_{i,\ell}}[f_{\ell}^{(i)}]]]W_{1,1}^{(q)}(x),\label{eq: cov function 1 critical 2}
        \end{multline}
        where, 
        \begin{equation*}
          B(t_{i-1},t_i,f_j^{(i)},f_{\ell}^{(i)}) = \int_0^{t_i-t_{i-1}}\frac{(t_i-v)^{p_1-1}v^{p_{i, j}+p_{i, \ell}-2}}{(p_1-1)!(p_{i, j}-1)!(p_{i, \ell}-1)!}\mathrm{d}v,
        \end{equation*}
        and where we have used the notation \(\nu_{i,j}\), $\lambda_{i, j}$ and $p_{i, j}$ in place of \(\nu_{f_j^{(i)}}\), $\lambda_{f_j^{(i)}}$ and $p_{f_j^{(i)}}$, respectively.
    \end{lemma}

    \medskip

    As in the small regime, we first give the proof of Theorem \ref{cor:crit} using Lemma \ref{lemma: fdd small inter crit}.
    
    \medskip
       \begin{proof}[Proof of Theorem \ref{cor:crit}]
            Fix \(k\geq 1\) and \(0\leq t_1< \dots < t_k<\infty\). 
            Using the notation defined in \eqref{eq:decomp}, we can write the left-hand side of \eqref{eq: main result crit} as 
            \begin{equation*}
                Y_{n, t}^1[f, f_1] + Y_{n, t}^1[f, \Phi_{\nu_f}[f]] + Y_{n, t}^1[f, f_3] - Y_{n, t}^2[f, f_1],
            \end{equation*}
            where $f_3 = f-f_1-\Phi_{\nu_f}[f]$ with $\nu_f$ being defined just before the statement of \eqref{thm: main result crit}, and where we have used that \(Y_{n, t}^2[f, f_2]=0\) since $f_2 \in B(E)\setminus \bigoplus_{i=1}^{m_L}A_{\lambda_i}$.
             By Corollary \ref{corollary: conv of phi for fdd}, we have that
    \begin{equation*}
        Y_{n, t}^1[f, f_1] - Y_{n, t}^2[f, f_1]
        \to 0,
    \end{equation*}
    in the sense of finite dimensional distributions, as \(n\rightarrow \infty\). Next, since \(f \in \mathrm{Ei}(\Lambda_C)\), it follows that \(f_3\in \mathrm{Ei}(\Lambda_C)\setminus \bigoplus_{i=1}^{m_C}A_{\lambda_i}\). Thus, \(f_3\) has zero projection in all critical eigenspaces with \(\Phi_{\nu_{f_3}}[f]=0\), therefore, by Lemma \ref{lemma: fdd small inter crit}, we have that
    \begin{equation*}
        Y_{n, t}^1[f, f_3] \to 0,
    \end{equation*} 
    in the sense of finite dimensional distributions, as \(n\rightarrow \infty\). Thus, it remains to show that \((Y_{n,\cdot}[f, \Phi_{\nu_f}[f]])_{n\geq 1}\) converges to \(Z_C^f\) in the sense of finite dimensional distributions. To this end, following similar steps to the proof of Theorem \ref{cor:small}, 
            for each \(1 \leq j \leq k\), we have that
            \begin{multline}
            X_{nt_j}[\Phi_{\nu_f}[f]]
               =\sum_{i=2}^j \left( X_{nt_i}[\psi_{n(t_j - t_i)}[\Phi_{\nu_f}[f]]] - \mathbb E[X_{nt_i}[\psi_{n(t_j - t_i)}[\Phi_{\nu_f}[f]]] | \mathcal F_{nt_{i-1}}] \right) \\ +  X_{nt_1}[\psi_{n(t_j - t_1)}[\Phi_{\nu_f}[f]]].     \label{eq: lem eq 1 crit}   
            \end{multline}
            The aim is to rewrite the right-hand side of \eqref{eq: lem eq 1 crit}, after an appropriate renormalisation, so that the summands are terms to which we can apply Lemma \ref{lemma: fdd small inter crit}. However, for all summands except \(i=j\), the functions in the argument of the \(X_{nt_i}\) depend on \(n\). To deal with this, we make some observations. Firstly, from the definition of $\nu_f$, it is not too difficult to see that $\nu_{\Phi_{\nu_f}[f]} = \nu_f$ so that $\lambda_{\Phi_{\nu_f}[f]} = \lambda_f$, and similarly for $\Phi_{\nu_f, b}[f]$ in place of $\Phi_{\nu_f[f]}$. Then, combining this with \ref{H1b}, we have, for \(1\leq i \leq j-1\),
            \begin{equation*}
                \psi_{n(t_j-t_i)}[\Phi_{\nu_f}[f]]=\e^{(\lambda_f + \mathcal{N})n(t_j-t_i)}\Phi_{\nu_f}[f].
            \end{equation*}
            Next recall that, each $\varphi_{a, b}^{(c)}$ is of rank $b$, meaning that $p(\varphi_{a, b}^{(c)}):=p_{\varphi_{a, b}^{(c)}} = b$. Due to linearity of $\mathcal N$ and the fact that $\Phi_{a, b}[f]$ is a sum over rank $b$ eigenfunctions, it follows that $p(\Phi_{a, b}[f]) = b$. Moreover, from the definition of $p(\cdot)$, for $f \in B(E)$ and $0 \le c \leq p_f-1$, we have $p(\mathcal N^cf) = p_f - c$. Combining this with the previous observation means that for each \(1\leq b \leq p_f\), \(0\leq c \leq b-1\), we have that
            \begin{equation*}
                p(\mathcal{N}^c\Phi_{\nu_f,b}[f]) = b-c,
            \end{equation*} 
            which in particular implies that $\mathcal N^c \Phi_{\nu_f,b}[f] = 0$ for all $c \ge b$.
            Therefore, combining these observations to expand the summands in \eqref{eq: lem eq 1 crit}, we have, for \(1\leq i \leq j-1\),
            \begin{align*}
                X_{nt_i}[\psi_{n(t_j - t_i)}[\Phi_{\nu_f}[f]]]
                &= \e^{\lambda_f n(t_j-t_i)}X_{nt_i}[\e^{\mathcal{N}n(t_j-t_i)}\Phi_{\nu_f}[f]] \\
                &= \e^{\lambda_f n(t_j-t_i)}X_{nt_i}\left[\sum_{b = 1}^{p_f} \e^{\mathcal{N}n(t_j-t_i)}\Phi_{\nu_f, b}[f]\right] \\
                &= \e^{\lambda_f n(t_j-t_i)}X_{nt_i}\left[\sum_{b = 1}^{p_f} \sum_{c = 0}^\infty \frac{n^c(t_j - t_i)^c}{c!}\mathcal N^c\Phi_{\nu_f, b}[f]\right] \\
                &= \e^{\lambda_f n(t_j-t_i)}X_{nt_i}\left[\sum_{b = 1}^{p_f}  \sum_{c = 0}^{b-1} \frac{n^{c}(t_j - t_i)^c}{c!}\mathcal N^c\Phi_{\nu_f, b}[f]\right]\\
                &= \e^{\lambda_f nt_j} \sum_{b = 1}^{p_f}  \sum_{c = 0}^{b-1} \frac{n^{c}(t_j - t_i)^c}{c!}\e^{-\lambda_f nt_i} X_{nt_i}\left[\mathcal N^c\Phi_{\nu_f, b}[f]\right].
            \end{align*}
Then, putting this back into \eqref{eq: lem eq 1 crit} and renormalising by ${\rm e}^{-\lambda_f nt_j}n^{-(2p_f + p_1 - 2)/2}$, we obtain
    {\begin{align*}
                Y_{n,t_j}^1[\Phi_{\nu_f}[f]]
                &=Y^*_{n,t_{j-1},t_j}[\Phi_{\nu_f}[f]] \\ 
                &+ \sum_{b=1}^{p_f}n^{-(p_f-b)}\sum_{c=0}^{b-1}\frac{1}{c!}\left(\sum_{i=2}^{j-1} {(t_j-t_i)^c}Y^*_{n,t_{i-1},t_i}[\mathcal{N}^c\Phi_{\nu_f,b}[f]] +   {(t_j-t_1)^c}Y^1_{n,t_1}[\mathcal{N}^c\Phi_{\nu_f,b}[f]]\right).
            \end{align*}

            This and Lemma \ref{lemma: fdd small inter crit} imply that
            \begin{align}
                Y_{n,t_j}^1&[\Phi_{\nu_f}[f]] -Y^*_{n,t_{j-1},t_j}[\Phi_{\nu_f}[f]] \notag \\
                &-\sum_{c=0}^{p_f-1}\frac{1}{c!}\F(\sum_{i=2}^{j-1} {(t_j-t_i)^c}Y^*_{n,t_{i-1},t_i}[\mathcal{N}^c\Phi_{\nu_f,p_f}[f]] +  {(t_j-t_1)^c}Y^1_{n,t_1}[\mathcal{N}^c\Phi_{\nu_f,p_f}[f]]\R)
                \coninprob 0
                \label{eq: decomp for two terms}
            \end{align} 
           as \(n\rightarrow \infty\). 
 
           Now, for \(1\leq j \leq k\), all but the first term on the left-hand side of \eqref{eq: decomp for two terms} are random variables to which we may apply Lemma \ref{lemma: fdd small inter crit} to yield
            \begin{equation*}
                 Y_{n,t_j}^1[\Phi_{\nu_f}[f]]\conindis Z^{(j)}, \quad 1\leq j \leq k,
            \end{equation*}
            as \(n\rightarrow \infty\), where conditionally on \(\mathcal{W}(x)\), the \(Z^{(j)}\) are jointly Gaussian.} Thus, the lemma follows if we can show, conditionally on \(\mathcal{W}(x)\), the covariance matrix of \((Z^{(1)},\dots ,Z^{(k)})\) is equal to the covariance matrix of \((Z_C^f(t_1),\dots ,Z_C^f(t_k))\), almost surely. To show this, the same approach as in the proof of Theorem \ref{cor:small} is used, i.e.\ we consider two time points and use consistency of the limit. By \eqref{eq: decomp for two terms} for two time points, $t_i, t_j$, \(1\leq i\leq j \leq k\), we have
            \begin{align}
              &Y_{n,t_j}^1[\Phi_{\nu_f}[f]] -Y^*_{n,t_{i},t_j}[\Phi_{\nu_f}[f]]- \sum_{c=0}^{p_f-1}\frac{1}{c!}{(t_j-t_i)^c}Y^1_{n,t_i}[\mathcal{N}^c\Phi_{\nu_f,p_f}[f]] \coninprob 0,\label{eq: convex + generic function}
            \end{align}
        as \(n\rightarrow \infty\). Thus, again by Lemma \ref{lemma: fdd small inter crit}, we have that
            \begin{align*}
                &Y_{n,t_j}^1[\Phi_{\nu_f}[f]]\conindis Z_*^{(j)},\\
                &Y_{n,t_i}^1[\Phi_{\nu_f}[f]]\conindis Z_*^{(i)},
            \end{align*}
             as \(n\rightarrow \infty\), where conditionally on \(\mathcal{W}(x)\), \(Z_*^{(i)},Z_*^{(j)}\) are jointly Gaussian with conditional covariance matrix
             \begin{align*}
                 &\Cov(Z_*^{(j)},Z_*^{(i)}|\mathcal{W}(x)) = 
                 \bs{1}_{\lambda_f=\bar \lambda_f}C(f,i,j)\sum_{q=1}^{k_{1,1}}\tilde\varphi_{1,1}^{(q)}[\mathcal{N}^{p_1-1}V[\mathcal{N}^{p_f-1}\Phi_{\nu_f,p_f}[f],\mathcal{N}^{p_f-1}\Phi_{\nu_f,p_f}[f]]]W_{1,1}^{(q)}(x),
             \end{align*}
             where 
             \begin{align*}
                 C(f,i,j) &= \sum_{c=0}^{p_f-1}\int_0^{t_i}\frac{(t_j-t_i)^{c}(t_i-v)^{p_1-1}v^{2p_f-c-2}}{c!(p_1-1)!(p_f-c-1)!(p_f-1)!}\mathrm{d}v\\
                 &=\int_0^{t_i}\frac{(t_j-t_i+v)^{p_f-1}(t_i-v)^{p_1-1}v^{p_f-1}}{(p_1-1)!(p_f-1)!^2}\mathrm{d}v\\
                 &=\int_0^{t_i}\frac{(t_j-v)^{p_f-1}v^{p_1-1}(t_i-v)^{p_f-1}}{(p_1-1)!(p_f-1)!^2}\mathrm{d}v,
             \end{align*}
             where the second line follows from the binomial theorem. The formula for \(\Cov(Z_*^{(j)},\bar Z_*^{(i)}|\mathcal{W}(x))\) follows similarly. For \(f,g \in \mathrm{Ei}(\Lambda_C)\), joint convergence follows by replacing \(f\) with \(g\) in \eqref{eq: convex + generic function}. This concludes the proof.
             
        \end{proof}
        \begin{proof}[Proof of Lemma \ref{lemma: fdd small inter crit}]

                The proof of Lemma \ref{lemma: fdd small inter crit} is almost a line-by-line copy of Lemma \ref{lemma: fdd small inter}. Thus, we omit the majority of the details, but clarify three important steps for which the proofs differ: the analogues of \eqref{eq: small fdd bound 13}, the third moment control used to bound \eqref{eq: berry esseen phi}, and the proof of \eqref{eq: need to show small fdd 4}. Accounting for these differences, the proof of Lemma \ref{lemma: fdd small inter crit} can be constructed by following the proof of Lemma \ref{lemma: fdd small inter} line-by-line. 
        
        Firstly, the third moment control used to bound \eqref{eq: berry esseen phi} is now given by Theorem 2.4 of \cite{Moments_Chris} in place of Theorems 2.1 and 2.2 of \cite{Moments_Chris}. 
        
        We next show the analogue of \eqref{eq: small fdd bound 13}. Let us first introduce some further notation. Let 
        \begin{align*}
            &\bs Y_n^{(1)} := (Y_{n,t_1}^1[f_1^{(1)}],\dots,Y_{n,t_1}^1[f_{r_1}^{(1)}]),  \\
        &\bs Y_n^{(i)} := (Y^*_{n,t_{i-1},t_{i}}[f_1^{(i)}],\dots,Y_{n,t_{i-1},t_{i}}^*[f_{r_i}^{(i)}]), \quad 2 \leq i \leq k.
        \end{align*}
          Also, let 
    \begin{equation*}
       \kappa_{n,1} := \log(n), \quad \kappa_{n,2}=nt_1 - \log(n),
    \end{equation*}
    and
    \begin{align}
      \label{eq: renorm term crit}
          &{\bs S}_n^{(1,1)}:=\bs Y_n^{(1)}-\E_{\delta_x}[\bs Y_n^{(1)}|\mathcal{F}_{\kappa_{n,1}}],\quad {\bs S}_n^{(1,2)} = \E_{\delta_x}[\bs Y_n^{(1)}|\mathcal{F}_{\kappa_{n,1}}].
      \end{align} 
        
        With this notation, the analogue of \eqref{eq: small fdd bound 13} is to show that for \(1\leq j,\ell \leq r_1\), 
    \begin{multline}
    \mathrm{Cov}_x(\bs S_{n,j}^{(1, 1)}, \bs S_{n,\ell}^{(1, 1)}|\mathcal{F}_{\kappa_{n,1}})\\
    \coninprob \bs{1}_{\lambda_{1, j}=\bar\lambda_{1, \ell}}B(0,t_1,f_j^{(1)},f_{\ell}^{(1)})\sum_{q=1}^{k_{1,1}}\tilde\varphi_{1,1}^{(q)}[\mathcal{N}^{p_1-1}V[\mathcal{N}^{p_{1, j}-1}\Phi _{\nu_{1,j}}[f_j^{(1)}],\mathcal{N}^{p_{1, \ell}-1}\Phi _{\nu_{1,\ell}}[f_{\ell}^{(1)}]]]W_{1,1}^{(q)}(x),
    \label{eq: identical argument for crit 0.5}
    \end{multline} 
as \(n\rightarrow \infty\), and that for \(2\leq i \leq k\), \(1\leq j,\ell \leq r_i\), 
\begin{multline}
            \mathrm{Cov}_x( \bs Y_{n,j}^{(i)}, \bs Y_{n,\ell}^{(i)}|\mathcal{F}_{n(t_{i}-t_{i-1})}) \\
            \coninprob\bs{1}_{\lambda_{i, j}=\bar\lambda_{i, \ell}}B(t_{i-1},t_i,f_j^{(i)},f_{\ell}^{(i)})\sum_{q=1}^{k_{1,1}}\tilde\varphi_{1,1}^{(q)}[\mathcal{N}^{p_1-1}V[\mathcal{N}^{p_{i, j}-1}\Phi _{\nu_{i,j}}[f_j^{(i)}],\mathcal{N}^{p_{i, \ell}-1}\Phi _{\nu_{i,\ell}}[f_{\ell}^{(i)}]]]W_{1,1}^{(q)}(x)
            \label{eq: identical argument for crit}
        \end{multline}
       as \(n\rightarrow \infty\).

        We first show \eqref{eq: identical argument for crit}. For \(2\leq i \leq k\), \(1\leq j,\ell \leq r_i\), we have that
        \begin{align}
            &\mathrm{Cov}_x(\bs Y_{n,j}^{(i)}, \bs Y_{n,\ell}^{(i)}|\mathcal{F}_{n(t_{i}-t_{i-1})})\nonumber\\
            &= \e^{-(\lambda_{i,j}+\lambda_{i,\ell}) nt_i}n^{-\beta(f_j^{(i)},f_{\ell}^{(i)})}\sum_{m=1}^{N_{nt_{i-1}}}\Cov_{X_{nt_{i-1}}^m}(X_{n(t_i-t_{i-1})}[f_j^{(i)}],X_{n(t_i-t_{i-1})}[f_\ell^{(i)}])\notag \\
            &=: \e^{-(\lambda_{i,j}+\lambda_{i,\ell}) nt_i}n^{-\beta(f_j^{(i)},f_{\ell}^{(i)})}\sum_{m=1}^{N_{nt_{i-1}}}h_{n,j,\ell}^{(i)}(X_{nt_{i-1}}^m),
            \label{eq: that other covariance result}
        \end{align}
        where, for \(f,g \in \mathrm{Ei}(\Lambda_C)\), \(\beta(f,g):=p_1+p_f+p_g-2\), \(X_{nt_{i-1}}^m\) denotes the \(m\)-th particle alive at time \(N_{nt_{i-1}}\), and \(h_{n,j,\ell}^{(i)}(x)=\Cov_{x}(X_{n(t_i-t_{i-1})}[f_j^{(i)}],X_{n(t_i-t_{i-1})}[f_\ell^{(i)}])\). Now, by Theorem 2.4 of \cite{Moments_Chris}, we have that
        \begin{equation*}
            n^{-\beta(f_j^{(i)},f_{\ell}^{(i)})}\e^{-\lambda_1 n(t_i-t_{i-1})}\|h_{n,j,\ell}^{(i)}\|_{\infty}<\infty.
        \end{equation*}
        Using this bound, an identical argument to \eqref{eq: almost usre conv smal} can be used to show that
        \begin{equation}
             \e^{-\lambda_1 nt_i}n^{-\beta(f_j^{(i)},f_{\ell}^{(i)})}\F(X_{nt_{i-1}}[h_{n,j,\ell}^{(i)}]-\sum_{a=1}^{m_L}\sum_{b=1}^{p_a}\sum_{b=1}^{k_{a,b}}\tilde\varphi_{a,b}^{(c)}[\e^{(\lambda_i+\mathcal{N})nt_{i-1}}h_{n,j,\ell}^{(i)}]W_{a,b}^{(c)}(x)\R) \coninprob 0,
             \label{eq: critical covariance conv 1}
        \end{equation}
        as \(n\rightarrow \infty\). The aim is to show that the above triple sum converges to the right-hand side of \eqref{eq: identical argument for crit}. This and the above convergence will then allow us to conclude. 
        
        First note that by Lemma \ref{lem:evo-2}, we have that 
        \begin{align}
         h_{n,j,\ell}^{(i)}(x) &= \psi_{n(t_i-t_{i-1})}[f_j^{(i)}f_{\ell}^{(i)}](x) + \int_0^{n(t_i-t_{i-1})}\psi_{n(t_i-t_{i-1})-s}[\gamma\eta_s^{(2)}[f_j^{(i)},f_{\ell}^{(i)}]](x)\mathrm{d}s \nonumber\\
         &\quad - \psi_{n(t_i-t_{i-1})}[f_j^{(i)}](x)\psi_{n(t_i-t_{i-1})}[f_{\ell}^{(i)}](x). \label{eq: second moment critical fdd}
        \end{align}
        Then, by \ref{H1b} and \eqref{eq:gen-eigenfunctions}, we have as \(n\rightarrow \infty\), for \(1\leq a \leq m\),
        \begin{align}
            &\e^{-\lambda_1 n(t_i-t_{i-1})}n^{-\beta(f_j^{(i)},f_{\ell}^{(i)})}\|\e^{\mathcal{N}nt_{i-1}}\psi_{n(t_i-t_{i-1})}[f_j^{(i)}f_{\ell}^{(i)}]\|_{\infty} \rightarrow0, \nonumber\\
             &\e^{(\mathrm{Re}\lambda_a-\lambda_1)nt_{i-1}}\e^{-\lambda_1 n(t_i-t_{i-1})}n^{-\beta(f_j^{(i)},f_{\ell}^{(i)})}\|\Phi_a[\e^{\mathcal{N}nt_{i-1}}(\psi_{n(t_i-t_{i-1})}[f_j^{(i)}]\psi_{n(t_i-t_{i-1})}[f_{\ell}^{(i)}])]\|_{\infty} \rightarrow0. \label{eq: crit covariance 2}
        \end{align}
        For the integral term in \eqref{eq: second moment critical fdd}, using an identical argument to that presented in the proof of Lemma 3.4 of \cite{Moments_Chris} (see the equation after (3.16) in \cite{Moments_Chris} for the details), we have that as \(n\rightarrow \infty\), uniformly for \(x \in E\),
        \begin{align}
           &n^{-\beta(f_j^{(i)},f_{\ell}^{(i)})}\Bigg|\int_0^{n(t_i-t_{i-1})} \e^{-\lambda_1 n(t_i-t_{i-1})}\e^{\mathcal{N}nt_{i-1}}\psi_{n(t_i-t_{i-1})-s}[\gamma\eta_s^{(2)}[f_j^{(i)},f_{\ell}^{(i)}]](x)\nonumber\\
           &-\e^{(\lambda_{i, j}+\lambda_{i, \ell}-\lambda_1)s}\frac{s^{p_{i, j}+p_{i, \ell}-2}}{(p_{i, j}-1)!(p_{i, \ell}-1)!}\sum_{\alpha=0}^{p_1-1}\frac{(nt_{i-1})^{\alpha}\mathcal{N}^{\alpha}\e^{\mathcal{N}(n(t_{i}-t_{i-1})-s)}}{\alpha!}\times\dots\nonumber\\
           &\dots\times\Phi_{1,1}[V[\mathcal{N}^{p_{i, j}-1}\Phi_{\nu_{i,j}}[f_j^{(i)}],\mathcal{N}^{p_{i, \ell}-1}\Phi_{\nu_{i,\ell}}[f_{\ell}^{(i)}]]](x)\mathrm{d}s \Bigg|\rightarrow 0. \label{eq: crit covariance 3}
        \end{align}
         Next, taking an expansion of \(e^{\mathcal{N}(n(t_{i}-t_{i-1})-s)}\) in \eqref{eq: crit covariance 3} and applying the binomial theorem gives
         \begin{align*}
           &n^{-\beta(f_j^{(i)},f_{\ell}^{(i)})}\Bigg|\int_0^{n(t_i-t_{i-1})} \e^{-\lambda_1 n(t_i-t_{i-1})}\e^{\mathcal{N}nt_{i-1}}\psi_{n(t_i-t_{i-1})-s}[\gamma\eta_s^{(2)}[f_j^{(i)},f_{\ell}^{(i)}]](x)\nonumber\\
           &-\e^{(\lambda_{i, j}+\lambda_{i, \ell}-\lambda_1)s}\frac{(nt_i-s)^{p_1-1}s^{p_{i, j}+p_{i, \ell}-2}}{(p_1-1)!(p_{i, j}-1)!(p_{i, \ell}-1)!}\Phi_{1,1}[\mathcal{N}^{p_1-1}V[\mathcal{N}^{p_{i, j}-1}\Phi_{\nu_{i,j}}[f_j^{(i)}],\mathcal{N}^{p_{i, \ell}-1}\Phi_{\nu_{i,\ell}}[f_{\ell}^{(i)}]]](x)\mathrm{d}s \Bigg|\nonumber\\
           &\rightarrow 0, 
        \end{align*}
        as \(n\rightarrow \infty\), where for each summand we have used that this exponential expansion is dominated by the \((p_1-\alpha-1)\)th term. This, \eqref{eq: crit covariance 2}, \eqref{eq: crit covariance 3} and standard integration techniques then give
        \begin{multline*}
            \e^{-\lambda_1 n(t_i-t_{i-1})}n^{-\beta(f_j^{(i)},f_{\ell}^{(i)})}\left\|\e^{\mathcal{N}nt_{i-1}}h_{n,j,\ell}^{(i)}\right. \\
            \left.-\bs{1}_{\lambda_{i, j}=\bar\lambda_{i, \ell}} \Phi_{1,1}[\mathcal{N}^{p_1-1}V[\mathcal{N}^{p_{i, j}-1}\Phi_{\nu_{i,j}}[f_j^{(i)}],\mathcal{N}^{p_{i, \ell}-1}\Phi_{\nu_{i,\ell}}[f_{\ell}^{(i)}]]]B(t_{i-1},t_i,f_j^{(i)},f_{\ell}^{(i)})\right\|_{\infty} \rightarrow 0
        \end{multline*}
        as \(n\rightarrow \infty\). From this and \eqref{eq: critical covariance conv 1} we can conclude.
         
         For \(i=1\), an identical argument holds, we omit the proof.
         
Finally, for the analogue of \eqref{eq: need to show small fdd 4}, we apply the same argument as in \eqref{eq: small regime last condition for fdd} using \ref{H1b} and Theorem 2.4 of \cite{Moments_Chris} giving 
\begin{align*}
        \E_{\delta_x}\F[\F\|\boldsymbol{S}_{n}^{(1,2)}\R\|^2_2\R]  &=n^{-(2p_{1,j}+p_1-2)}\e^{-\lambda_1nt_1}\sum_{j=1}^{r_1}\E_{\delta_x}\F[|X_{\kappa_{n,1}}[\psi_{\kappa_{n,2}}[f_j^{(1)}]]|^2\R]\\
         &\leq C n^{-1}\kappa_{n,1}^{2p_{1,j}-1}\rightarrow 0 
     \end{align*}
as \(n\rightarrow \infty\). This and Markov's inequality imply that
\begin{equation*}
    {\bs S}_n^{(1,2)} \coninprob 0,
\end{equation*}
as \(n\rightarrow \infty\) as required.
\end{proof}

To show the third conditions of Theorem \ref{Theorem:Bill con} holds, we require the following intermediate lemma which is the critical analogue of Lemma \ref{lemma: interim large 1}.
\begin{lemma}
        \label{lemma: interim critical 1 2}
            Fix \(k\geq 1\). Assume that \ref{H2} holds and \(p_1=1\) in \ref{H1b}. Then, for \(1\leq \ell \leq k\), there exists a constant \(C_{\ell}\), such that
            \begin{equation}
            \label{eq: proof by induct result crit}
                \sup_{x \in E,f_1,\dots,f_{\ell} \in B_1(E)\setminus\bigoplus_{i=1}^{m_C}A_{\lambda_i},t\geq 0}{\rm e}^{-\ell\lambda_1t/2}\F|\E_{\delta_x}\F[\prod_{i=1}^{\ell}X_t[f_i]\R]\R| \leq C_{\ell}.
            \end{equation}
        \end{lemma}    
        The proof of this lemma follows the same argument as Lemma \ref{lemma: interim large 1}, thus we omit the details.
        \medskip

        \begin{lemma}
        \label{lemma: tightness small crit}
Under the assumptions of Theorem \ref{thm: main result crit}, for each \(T\geq 0\), \(f\in C_{\psi}(\Lambda_C)\), there exists a constant \(C_{f,T}>0\), such that, for any \(0\leq r \leq s \leq t \leq T\), \(n\geq 1\),
            \begin{equation}
    	       \sup_{x \in E}\E_{\delta_x}\big[|Y_{n,t}[f]-Y_{n,s}[f]|^{k_f}|Y_{n,s}[f]-Y_{n,r}[f]|^{k_f}\big] \leq C_{f,T}(t-r)^{3/2}. 
               \label{eq: bound tightness crit 1}
            \end{equation}
        \end{lemma}
   \begin{proof}
            To show \eqref{eq: bound tightness crit 1} we use one of two approaches depending on whether or not \(n(t-s)\) is less than 1. First, assume that \(0\leq n(t-s)<1\). In this case, the proof is almost a line-by-line copy of Lemma \ref{lem:tight}. Thus, we omit the majority of the proof, but clarify two details for which the proofs differ. Firstly, the sufficient bounds presented in \eqref{eq: tight control 1l}-\eqref{eq: tight control 3l} now take the form, for \(n\geq 1\), \(0\leq s \leq t \leq T\), 
            \begin{align*}
&\E_{\delta_x}\big[\E\big[|\alpha_{n,s,t}|^{k_f}|\mathcal{F}_{ns}]^2]\leq C_1(t-s)^2,\\
&\E_{\delta_x}\big[|\alpha_{n,s,t}|^{2{k_f}}]\leq C_1(t-s),\\
&\E_{\delta_x}\big[|\beta_{n,s,t}|^{2{k_f}}] \leq C_1(t-s)^2,     
\end{align*}
where
    \begin{align*}
      \alpha_{n,s,t} &= Y_{n, t}^1[f] - Y_{n, s}^1[f, f_2] - \mathbb E[Y_{n, t}^1[f, f_1] | \mathcal F_{ns}], \\
      \beta_{n,s,t} &= (\mathbb E[Y_{n, t}^1[f, f_1]|\mathcal F_{ns}] - Y_{n, s}^1[f, f_1]) - (Y_{n, t}^2[f] - Y_{n, s}^2[f]),
    \end{align*}
Secondly, the analogue of \eqref{eq: large assumption exp bound 12} is slightly more involved. The function which appears in this bound takes the form
\begin{equation*}
    h_{n,s,t} := \e^{-\lambda_f n(t-s)}\psi_{n(t-s)}[f_2]-f_2.
\end{equation*}
We split this function into two components: 
\begin{align}
   &h_{n,s,t}^{(1)} =  \e^{-\lambda_f n(t-s)}\psi_{n(t-s)}[\Phi_{\nu_f}[f_2]]-\Phi_{\nu_f}[f_2]=(\e^{\mathcal{N}n(t-s)}-1)\Phi_{\nu_f}[f_2],\\
   &h_{n,s,t}^{(2)} = \e^{-\lambda_f n(t-s)}\psi_{n(t-s)}[f_2-\Phi_{\nu_f}[f_2]]-(f_2-\Phi_{\nu_f}[f_2]).
   \label{eq: h functions split}
\end{align}
Using an identical argument to \eqref{eq: assum ext}, we have that 
\begin{equation*}
    \sup_{0\leq n(t-s)\leq 1}(n(t-s))^{-1/{k_f}}\|h_{n,s,t}^{(2)}\|_{\infty}<\infty.
\end{equation*}
Furthermore, by construction, \(h_{n,s,t}^{(2)}\in B(E)\setminus \bigoplus_{i=1}^{m_C}A_{\lambda_i}\). Therefore we can apply Lemma \ref{lemma: interim critical 1 2} to obtain
\begin{equation}
\label{eq: crit need 2221}
    n^{-2{k_f}(p_f-\frac{1}{2})}\e^{-{k_f}\lambda_1ns}\E_{\delta_x}\F[\F|X_{ns}\F[h_{n,s,t}^{(2)}\R]\R|^{2{k_f}}\R]\leq C(t-s)^2,
\end{equation}
where we have used that \(k\geq 2\) to obtain a bound independent of \(n\). Next, we can write
\begin{equation}
\label{eq: tightness crit need 1}
    h_{n,s,t}^{(1)} = \sum_{i=1}^{p_f-1}\frac{(n(t-s)\mathcal{N})^i\Phi_{\nu_f}[f_2]}{i!}.
\end{equation} 
For \(1\leq i \leq p_f-1\), we have that \(p(\mathcal{N}^i\Phi_{\nu_f}[f_2])=p_f-i\). Therefore, by Theorem 2.4 of \cite{Moments_Chris}, for \(1\leq i \leq p_f-1\), \(n\geq 1\)
\begin{equation}
  n^{-2{k_f}(p_f-i-\frac{1}{2})}\e^{-{k_f}\lambda_1ns}\E_{\delta_x}\F[\F|X_{ns}\F[\mathcal{N}^i\Phi_{\nu_f}[f_2]\R]\R|^{2{k_f}}\R] \leq C.
  \label{eq: tightness crit need 5}
\end{equation}
Since \(k\geq 2\), this and \eqref{eq: tightness crit need 1} imply 
\begin{equation}
    n^{-2{k_f}(p_f-\frac{1}{2})}\e^{-{k_f}\lambda_1ns}\E_{\delta_x}\F[\F|X_{ns}\F[h_{n,s,t}^{(1)}\R]\R|^{2{k_f}}\R] \leq C(t-s)^2. \label{eq: tightness crit need 10}
\end{equation}
Combining this and \eqref{eq: crit need 2221} gives the analogue of \eqref{eq: large assumption exp bound 12}. From here, the remainder of the proof can be constructed by using the same arguments as the proof of Lemma \ref{lem:tight} line-by-line. 

We now move onto the case where \(n(t-s)>1\). In this case, we set
      \begin{align*}
      \alpha_{n,s,t} &= Y_{n, t}^1[f, f_2] - Y_{n, s}^1[f, f_2], \\
      \beta_{n,s,t} &= Y_{n, t}[f, f_1] - Y_{n, s}[f, f_1].
    \end{align*}
   Again, for \(n\geq 1\), \(0\leq s \leq t \leq T\), the sufficient bounds take the form,  
\begin{align}
&\E_{\delta_x}\big[\E\big[|\alpha_{n,s,t}|^{k_f}|\mathcal{F}_{ns}]^2]\leq C_1(t-s)^2, \label{eq: crit tightness bound show 1}\\
&\E_{\delta_x}\big[|\alpha_{n,s,t}|^{2{k_f}}]\leq C_1(t-s),\label{eq: crit tightness bound show 2}\\
&\E_{\delta_x}\big[|\beta_{n,s,t}|^{2{k_f}}] \leq C_1(t-s)^2 \label{eq: crit tightness bound show 3}.     
\end{align}
We start by showing \eqref{eq: crit tightness bound show 1}. {Firstly, there exists constant \(C\), such that 
    \begin{align}
       \E\big[|\alpha_{n,s,t}|^{k_f}|\mathcal{F}_{ns}]^2
       &\leq C\E[|Y_{n,t}^1[f, f_2]-\E[Y_{n,t}^1[f, f_2]|\mathcal{F}_{ns}]|^{k_f}|\mathcal{F}_{ns}]^2
       +C\F|\E[Y_{n,t}^1[f, f_2]|\mathcal{F}_{ns}]-Y^1_{n,s}[f, f_2]\R|^{2{k_f}}.   \label{eq: cond exp crit}
    \end{align}
    Taking expectations in the final line we obtain 
    \begin{equation*}
    \E_{\delta_x}\F[\F|\E[Y_{n,t}^1[f, f_2]|\mathcal{F}_{ns}]-Y^1_{n,s}[f, f_2]\R|^{2{k_f}}\R]=n^{-2{k_f}(p_f-\frac{1}{2})}\e^{-{k_f}\lambda_1ns}\E_{\delta_x}\F[\F|X_{ns}\F[h_{n,s,t}\R]\R|^{2{k_f}}\R].  
    \end{equation*}
    As before, to bound the right-hand side, we take the decomposition \eqref{eq: h functions split}. By construction, \(h_{n,s,t}^{(2)} \in B(E)\setminus \bigoplus_{i=1}^{m_C}A_{\lambda_i}\), and \(\sup_{n\geq 1,0\leq s \leq t \leq 1}\|h_{n,s,t}\|_{\infty}<\infty\). Therefore, by Lemma \ref{lemma: interim critical 1 2},
\begin{equation}
\label{eq: tightness crit need 2}
     n^{-2{k_f}(p_f-\frac{1}{2})}\e^{-{k_f}\lambda_1ns}\E_{\delta_x}\F[\F|X_{ns}\F[h_{n,s,t}^{(2)}\R]\R|^{2{k_f}}\R]
        \leq Cn^{-{k_f}}\leq Cn^{-({k_f}-2)}(n(t-s))^{-2}(t-s)^2\leq C(t-s)^2.
\end{equation} Combining this with \eqref{eq: tightness crit need 1}-\eqref{eq: tightness crit need 10} gives us
\begin{equation}
    n^{-2{k_f}(p_f-\frac{1}{2})}\e^{-{k_f}\lambda_1(n+s)}\E_{\delta_x}\F[\F|X_{ns}\F[h_{n,s,t}\R]\R|^{2{k_f}}\R]
    \leq C(t-s)^2.
    \label{eq: critical tight need 4}
\end{equation}  For the first term on the right-hand side of \eqref{eq: cond exp crit}, we use the same argument as \eqref{eq: large part 2 exp bound} to obtain 
\begin{align}
      \E&[|Y_{n,t}^1[f, f_2]-\E[Y_{n,t}^1[f, f_2]|\mathcal{F}_{ns}]|^{k_f}|\mathcal{F}_{ns}]\nonumber \\
       &\leq C_1n^{-{k_f}(p_f-\frac{1}{2})}\e^{-{k_f}\lambda_1nt/2}N_{ns}^{ {k_f}/2 }\max_{i=1,\dots,N_{ns}}\E_{\delta_{X_{ns}^i}}\big[\big|X_{n(t-s)}[f, f_2]-\psi_{n(t-s)}[f, f_2](X_{ns}^i)\big|^{k_f}\big] \nonumber\\
       &\leq C_2\e^{-{k_f}\lambda_1ns/2}N_{ns}^{ {k_f}/2 }(t-s), \label{eq: large part 2 exp bound crit}
    \end{align}
    where the final inequality is due to Theorem 2.4 of \cite{Moments_Chris}. This and Theorem 2.1 of \cite{Moments_Chris} imply

\begin{equation}
    \E_{\delta_x}[\E[|Y_{n,t}^1[f, f_2]-\E[Y_{n,t}^1[f, f_2]|\mathcal{F}_{n+s}]|^{k_f}|\mathcal{F}_{ns}]^2]\leq C(t-s)^2.
    \label{eq: tech argu 1}
\end{equation}
Using this and \eqref{eq: critical tight need 4} in \eqref{eq: cond exp crit} gives \eqref{eq: crit tightness bound show 1}. Note that under the necessary renormalization for convergence of finite dimensional distributions  \eqref{eq: large part 2 exp bound crit} and \eqref{eq: tech argu 1} only hold if \(p_1=1\), the same bounds with \(p_1>1\) result in a factor of \(n^{k_f(p_1-1)}\) on the right-hand side of \eqref{eq: tech argu 1}, hence the need for this assumption. The proof of \eqref{eq: crit tightness bound show 2} is almost identical, so we omit the details. Finally, we show \eqref{eq: crit tightness bound show 3}. We have that
\begin{equation*}
    \E_{\delta_x}\big[|\beta_{n,s,t}|^{2{k_f}}]\leq C(\E_{\delta_x}[|Y_{n,t}[f, f_1]|^{2{k_f}}]+\E_{\delta_x}[|Y_{n,s}[f, f_1]|^{2{k_f}}]).
\end{equation*}
By the \(2{k_f}\)-th moment analogue of \eqref{eq: need for markov ineq}, we have that
\begin{align*}
   &\E_{\delta_x}[|Y_{n,t}[f,f_1]|^{2{k_f}}] \leq Cn^{-{k_f}} \leq C(t-s)^2(n(t-s))^{-2}\leq C(t-s)^2\\ 
   &\E_{\delta_x}[|Y_{n,s}[f,f_1]|^{2{k_f}}]\leq C(t-s)^2.
\end{align*}
Thus, \eqref{eq: crit tightness bound show 3} holds completing the proof.}
        \end{proof}         
        \medskip

\begin{proof}[Proof of Theorem \ref{thm: main result crit}]
Let \(f \in C_{\psi}(\Lambda_C)\) and for \(n\geq 1\), \(t\geq 0\), let \(Y_{n,t}[f]\) be as in \eqref{eq: compact ver of crit}. We are going to apply Theorem \ref{Theorem:Bill con} to \((Y_{n,\cdot}[f])_{n\geq 1}\) with limit \(Z_C^f\). Conditions 1 and 3 of the aforementioned theorem have already been shown in Theorem \ref{cor:crit} and Lemma \ref{lemma: tightness small crit} respectively. Therefore, it remains to check the second condition. We have, for \(x \in E\), and \(T\geq\delta \geq 0\), 
\begin{align}
   \E[\Var(Z_C^f(T)-Z_C^f(T-\delta)|\mathcal{W}(x))]\leq C_{\mathcal{W}(x)}(\mathcal{C}^{(4)}_{T,}(f,\bar f)+\mathcal{C}^{(4)}_{T-\delta,T-\delta}(f,\bar f)-2\mathcal{C}^{(4)}_{T-\delta,T}(f,\bar f)),
   \end{align}
 where \(C_{\mathcal{W}(x)}\) has finite first moment. The right-hand side is a sum of polynomial functions that are continuous at \((T,T)\), thus the left-hand side is continuous at \(T\). This means that we can apply Theorem \ref{Theorem:Bill con} to obtain functional convergence in \(\mathbb{D}[0,T]\). Since \(T\) was arbitrary, we have convergence in \(\mathbb{D}[0,\infty]\).
\end{proof}

\appendix
\section{Useful results}
In this appendix we house several results that are used throughout the article. The first result is an evolution equation for the joint moments of \(X\) that can be found in \cite[Lemma 4.1]{Moments_Chris}.

\begin{lemma}\label{lem:evo-2}
  Fix \(k\geq 2\). Assume that 
  \begin{equation}
  \label{eq: off bound apendix}
    \sup_{x\in E}\mathcal{E}_x[N^k]<\infty. 
  \end{equation}
  
  Then, for any \(f_1,\dots,f_k \in B(E)\), \(x\in E\), and \(t\geq 0\), we have that
\begin{equation}\label{eq:evo-23}
  \E_{\delta_x}\F[\prod_{i=1}^kX_t[f_i]\R] =\psi_t[f_1 \cdots f_k](x)+\int_0^t \psi_s\left[\gamma \eta_{t-s}^{(k)}[f_1,\dots,f_k]\right](x) \mathrm{d} s, \quad t \geq 0,
\end{equation}
where 
\[
   \eta_{t-s}^{(k)}[f_1,\dots,f_k](x) = \mathcal E_x \left[\sum_{\sigma \in \mathcal{P}^*(k)}\sum_{\bs i \in B_{|\sigma|,N}}\prod_{j = 1}^{|\sigma|} \E_{\delta_{x_{i_j}}}\F[\prod_{p \in \sigma_j}X_{t-s}[f_p]\R] \right],
\]
$\mathcal{P}^*(k)= \mathcal{P}(k)\setminus \{\{\{1, \dots, k\}\}\}$ with \(\mathcal{P}(k)\) denoting the set of partitions of \(\{1, \dots, k\}\) and 
\begin{align*}
    &B_{n,m}:= \{\bs i \in \{1, \dots, m\}^{n}: i_j \neq i_k, 1\leq j,k \leq n\}.
\end{align*}
\end{lemma}

\medskip

The next result shows that the regularity on the Markov process imposed by taking $f \in C_\psi(E)$ yields similar regularity of the branching process, which is used in several of the proofs.

\begin{lemma}
           \label{Lemma: conv of h5 over} 
               Assume that \ref{H1b} and \ref{H4} hold, and let $f \in C_{\psi}(E)$. Then, for any \(T\geq 0\), there exists \(C_{f,T}>0\), such that
\begin{align}
                 & \sup_{x \in E, 0 \leq t \leq {\color{black}T}}t^{-1/k_f}\ |\E_{\delta_x}[X_{t}[f]]-f(x)| \leq C_{f,T}, \label{eq: tight assum 111}\\
                 &\sup_{x \in E,0\leq t \leq {\color{black}T}} t^{-1}\E_{\delta_x}[(X_t[f]-\mathbf \E_{\delta_x}[X_t[{f}]])^{k_f}]\leq C_{f,T},  \label{eq: tight assum 211}\\
                &\sup_{x \in E,0\leq t \leq {\color{black}T}} t^{-1}\E_{\delta_x}[(X_t[f]-\mathbf \E_{\delta_x}[X_t[{f}]])^{2k_f}]\leq C_{f,T}.  \label{eq: tight assum 311}
                \end{align}
           \end{lemma}
           \begin{proof}
               
            Fix \(T\geq 0\). First note that since \(f \in B(E)\), the suprema in \eqref{eq: tight assum 1112} and \eqref{eq: tight assum 2112} can be extended to hold over \([0,T]\). Then, writing
            \begin{equation*}
                \E_{\delta_x}[X_t[f]]-f(x) =  \mathbf E_x[f(\xi_t)]-f(x)+\E_{\delta_x}[X_t[f]]-\mathbf E_x[f(\xi_t)],
            \end{equation*}
            it follows that in order to prove \eqref{eq: tight assum 111}, it is sufficient to prove that there exists \(\tilde C_{f,T}\geq 0\), such that
            \begin{align}
                &\sup_{x\in E,0\leq t \leq T}t^{-1/k_f}\big|\E_{\delta_x}[X_t[f]]-\mathbf E_x[f(\xi_t)]\big| \leq \tilde C_{f,T}. \label{eq: x xi switch 1}
            \end{align}
            To show this, for \(t\geq 0\), let
            \begin{equation*}
                A_t = \{\text{\(X\) branches at least once before time \(t\)}\}.
            \end{equation*}
            Until the first branch time, we can and do couple the paths of \(X\) and \(\xi\). Therefore, for any \(t\geq 0\),
            \begin{align*}
               | \E_{\delta_x}[X_t[f]]-\mathbf E_x[f(\xi_t)] |\leq \E_{\delta_x}[\bs{1}_{A_t} X_t[|f|]]+\|f\|_{\infty}\PP_{\delta_x}(A_t).
            \end{align*}
         Since \ref{H4} implies the branching rate is uniformly bounded, we have
        \begin{equation*}
             \PP_{\delta_x}(A_t) \leq 1-\e^{-\|\gamma \|_{\infty}t}\leq \|\gamma \|_{\infty}t \wedge 1.  
        \end{equation*}
        Furthermore, by splitting over the first branch time and using (M\(1\)) (which is implicitly assumed as \(C_{\psi}(E)\neq \emptyset\)), \ref{H1b} and the Markov property, there exists some \(\hat C_{f,T}\geq 0\), such that, for \(0\leq t \leq T\),
        \begin{align*}
            \E_{\delta_x}[\bs{1}_{A_t} X_t[|f|]] &= \bs E_{x}\F[\int_0^t \gamma(\xi_s)\e^{-\int_0^s \gamma(\xi_r)\mathrm{d}r} \mathcal{E}_{\xi_s}\F[\sum_{i=1}^{N} \E_{\delta_{x_i}}[X_{t-s}[|f|]]\R]\R]\mathrm{d}s\\
            &\leq \PP_{\delta_x}(A_t)\sup_{y \in E}\mathcal{E}_y[N]\sup_{z \in E,0\leq s \leq t}\E_{\delta_z}[X_s[|f|]] \leq \hat C_{f,T}\|f\|_{\infty}\PP_{\delta_x}(A_t).
        \end{align*}

        Thus,
        \begin{equation}
            \sup_{x\in E,0\leq t \leq T}t^{-1/k_f}\big|\E_{\delta_x}[X_t[f]]-\mathbf E_x[f(\xi_t)]\big| \leq (1+\hat C_{f,T})\|f\|_{\infty}\|\gamma\|_{\infty}, \label{eq: Functional assumption bp eq 1}
        \end{equation}
        as required.

        Hence, this and \eqref{eq: tight assum 1112} imply that there exists a $C_{f,T} > 0$ such that
            \begin{align}
                 \sup_{x \in E, 0 \leq t \leq {\color{black}T}}t^{-1/k_f}\ |\E_{\delta_x}[X_{t}[f]]-f(x)| \leq C_{f,T}. \label{eq:BPH5a}
                \end{align} 
        Both \eqref{eq: tight assum 211} and \eqref{eq: tight assum 311} follow similarly, thus we only present a sketch proofs for brevity. There exists a constant \(C\) such that
        \begin{align*}
            \E_{\delta_x}[(X_t[f]-\mathbf \E_{\delta_x}[X_t[{f}]])^{k_f}]\leq &C\E_{\delta_x}[(X_t[f]-f(\xi_t))^{k_f}]+C\mathbf{E}_{\delta_x}[(f(\xi_t)-\mathbf E_{x}[f(\xi_t)])^{k_f}]\\
            &+C(\mathbf E_{x}[f(\xi_t)]-\mathbf \E_{\delta_x}[X_t[{f}]])^{k_f},       
        \end{align*}
        where to house \(\xi\) on the same probability space as \(X\) we take the aforementioned coupling of \(X\) and \(\xi\) up to the first branching event of \(X\), then after we let them evolve independently. The second and third terms on the right-hand side are controlled by \eqref{eq: tight assum 2112} and \eqref{eq: Functional assumption bp eq 1} respectively. Thus, to prove \eqref{eq: tight assum 211}, it is left to show the existence of a constant \(C_{f,T}\), such that
        \begin{equation}
            \sup_{x \in E,0\leq t \leq T}t^{-1}\E_{\delta_x}[(X_t[f]-f(\xi_t))^{k_f}]\leq C_{f,T}, \label{eq: Functional assumption bp eq 2}
        \end{equation}
        and following identical steps that
        \begin{equation}
            \sup_{x \in E,0\leq t \leq T}t^{-1}\E_{\delta_x}[(X_t[f]-f(\xi_t))^{2k_f}]\leq C_{f,T} \label{eq: Functional assumption bp eq 3}
        \end{equation}
        to show \eqref{eq: tight assum 311}. These proofs now follow an identical approach to that of \eqref{eq: tight assum 111}. Indeed, both expectations are only non-zero on \(A_t\). On \(A_t\), the contribution from \(f(\xi_t)\) is bounded by \(\|f\|_{\infty}\), whereas the contribution from \(X_t[f]\) is controlled by
        \begin{equation*}
            \E_{\delta_x}[\bs{1}_{A_t} X_t[|f|]^{k_f}] = \bs E_{x}\F[\int_0^t \gamma(\xi_s)\e^{-\int_0^s \gamma(\xi_r)\mathrm{d}r} \mathcal{E}_{\xi_s}\F[\F(\sum_{i=1}^{N} \E_{\delta_{x_i}}[X_{t-s}[|f|]]\R)^{k_f}\R]\R]\mathrm{d}s,
        \end{equation*}
        where under (M\(k_f\)), for any fix \(0\leq t \leq T\), \cite{bmoments} shows
        \begin{equation*}
            \sup_{x \in E,0\leq t \leq T}\mathcal{E}_{x}\F[\F(\sum_{i=1}^{N} \E_{\delta_{x_i}}[X_{t}[|f|]]\R)^{k_f}\R]<\infty.
        \end{equation*}
        From these equations \eqref{eq: Functional assumption bp eq 2} follows using the same argument that showed \eqref{eq: tight assum 111}. Moreover, the previous two equations with \(k_f\) replaced by \(2k_f\) hold under (M\(2k_f\)) which gives us \eqref{eq: Functional assumption bp eq 3}.
\end{proof}

The next result is Billingsley \cite[Theorem 13.5]{bBill}, however since checking the conditions of this theorem is how we prove several functional convergence results, we state it here for the convenience of the reader. We write $\mathbb D[a, b]$ to denote the set of real-valued c\`adl\`ag paths on $[a, b]$ under the $J_1$-Skorokhod topology.

\medskip
    
\begin{theorem}
\label{Theorem:Bill con}
Let \(( Y_n)_{n\geq 1}\) be a sequence of processes in \(\mathbb D[a,b]\), and \(Y\) a process in~\(\mathbb D[a,b]\). Suppose the following conditions hold:
\begin{enumerate}
\item  For all increasing sequences of times, \(a\leq t_1 \leq \dots \leq t_k \leq b\), 
\[( Y_n(t_1),\dots, Y_n(t_k)) \conindis (Y(t_1),\dots, Y(t_k)),\] 
as \(n \rightarrow \infty\).
\item \( Y(b)- Y(b-\delta) \conindis 0\)
as \(\delta \rightarrow 0\).
\item There exists \(k\geq 2\), such that, for every \(n \geq 1\), and \(a \leq r \leq s \leq t \leq b\), \begin{equation*}
    \E \F [| Y_n(t)- Y_n(s)|^k |  Y_n(s)- Y_n(r)|^k\R] \leq C (t-r)^{3/2},
\end{equation*}
for some constant $C > 0$.
\end{enumerate}
Then, as $n \to \infty$,
\begin{equation*}
    Y_n(t) \conindis  Y(t)\text{ in } \mathbb D[a,b].
\end{equation*}
\end{theorem}

\medskip

\noindent 
The next result is an extension of the Berry-Esseen theorem. In the case where all the random variables have constant variance, then we do indeed recover the classical result. 

\medskip

\begin{theorem}
\label{theorem: Berry-Esseen}
    For \(n,k\geq 1\), let \(\bs{X}_1,\dots,\bs X_n\) be mean-zero independent random variables taking values in \(\mathbb{R}^k\). Assume that 
    \begin{equation*}
        \Sigma_n:=\sum_{i=1}^n \Var(\bs X_i) 
    \end{equation*}
    is non-degenerate. Then, letting \(\bs S_n = \sum_{i=1}^n \bs X_i\), there exists a constant $C_k > 0$ that depends only on \(k\), such that
    \begin{equation}
    \label{eq: general result}
        \sup_{A \in \mathcal{A}}\F|\mathbb{P}(\bs S_n \in A) - \mathbb{P}(\bs W_n \in A)\R| \leq C_k \sum_{i=1}^n \E\big[\|\Sigma_n^{-1/2}\bs X_i\|_2^3\big],
    \end{equation}
   where \(\mathcal{A}\) is the set of all measurable convex sets of \(\mathbb{R}^k\), and \(\bs W_n \sim \mathcal{N}(0,\Sigma_n)\).
    \end{theorem}
    \begin{proof}
    We will use the fact that the result was proved for $\Sigma_n = I$ in \cite{bGotze, bRaivc2018AMB}. 
    First note that, since $\Sigma_n$ is non-degenerate, \(\Sigma_n^{-1/2}\) is well-defined. Thus, we may apply the result for the case $\Sigma_n = I$ to the random variables \(\Sigma_n^{-1/2}\bs X_1,\dots,\Sigma_n^{-1/2} \bs X_n\) to yield
        \begin{equation}
        \label{eq: general result proof}
             \sup_{B \in \mathcal{A}}\F|\mathbb{P}(\Sigma_n^{-1/2}\bs S_n \in B) - \mathbb{P}(\Sigma_n^{-1/2} \bs W_n \in B)\R| \leq C_k \sum_{i=1}^n \E[\|\Sigma_n^{-1/2} \bs X_i\|_2^3],
        \end{equation}
        where indeed \(\Sigma_n^{-1/2}\bs W_n \sim \mathcal{N}(0,I)\). Finally, since linear transformations preserve both convexity and measurability, if \(A\in \mathcal{A}\), then \(\Sigma_n^{-1/2}A\in \mathcal{A}\). Therefore, since \(\Sigma_n^{-1/2}\) is not degenerate, for \(A\in \mathcal{A}\), taking \(B=\Sigma_n^{-1/2}A\) in \eqref{eq: general result proof} yields \eqref{eq: general result}.
    \end{proof}
    
\medskip

\noindent 
Next, we state and prove a result that gives an appropriate bound on moments of sums of independent mean-zero random variables.

\medskip

\begin{lemma}
        \label{lemma: combin}
           For {\color{black}\(k,n\geq 1\)}, let \(X_1,\dots,X_n\) be independent mean-zero random variables in \(\mathbb{C}\) with finite \(k\)-th moment. Then, there exists a constant \(C_k>0\) that depends only on \(k\), such that
           \begin{equation}
               \E\F[\F|\sum_{i=1}^nX_i\R|^k\R] \leq C_kn^{k/2}\max_{1\leq i \leq n}\E[|X_i|^k]. \label{eq: sum indep result}
           \end{equation}
           \end{lemma}
           \begin{proof}
           We first treat the case of real random variables, then the complex result immediately follows by splitting into real and imaginary parts. By the Marcinkiewicz–Zygmund inequality, we have
           \[
             \E\F[\F(\sum_{i=1}^nX_i\R)^k\R] 
             \le C_k \E\F[\left( \sum_{i = 1}^n |X_i|^2\right)^{k/2}\R] 
             = C_k n^{k/2}\E\F[\left( \frac1n\sum_{i = 1}^n |X_i|^2\right)^{k/2}\R],
           \]
           where $C_k$ is a positive constant that depends only on $k$. By Jensen's inequality, we have
           \[
           \E\F[\left( \frac1n\sum_{i = 1}^n |X_i|^2\right)^{k/2}\R]
           \le \frac1n \sum_{i = 1}^n \E[|X_i|^k] \le \max_{1 \le i \le n}\E[|X_i|^k].
           \]
           Putting these two inequalities together yields the result.
           \end{proof}

\section{Proof of Proposition \ref{prop: conditions for Linfinity spectrum}}\label{sec:H1+}
We split the proof into three steps. First, we record the \(L^2\)-regularising
properties of \(\psi\) which follow from the assumptions on the many-to-one
semigroup. Second, we invoke the spectral theory of compact operators to obtain
an \(L^2(E,\mu)\)-spectral decomposition of \(\psi_T\). Finally, we lift this
decomposition to the uniform \(B(E)\) asymptotic required in \ref{H1b}. 

By the assumptions of the proposition, (M$1$), \ref{H4} and the many-to-one formula, for each $t \ge T$, $f \in L^2(E) \cap B(E)$, we have 
\begin{equation}\label{eq:density}
 \psi_t[f](x) = \int_E q_t(x, y)f(y)\mu(dy),
\end{equation}
where
\begin{equation}\label{eq:B.1}
\|\psi_t[f]\|_\infty \le {\rm e}^{C t}\sup_{x \in E}\|p_t(x, \cdot)\|_{L^2(E, \mu)}\|f\|_{L^2(E, \mu)},
\end{equation}
with $\mu$ as in the first assumption of the proposition and $C = \sup_{x \in E}\gamma(x)(\mathcal E_x[N] - 1)$. Consequently, for each \(t\geq T\), \(\psi_t\) extends uniquely to a bounded
operator from \(L^2(E,\mu)\) to \(B(E)\). Since \(\mu(E)<\infty\), \(\psi_t\) is also a bounded operator on
\(L^2(E,\mu)\), for every \(t\geq T\).
 
We next extend the semigroup property. Let \(f\in L^2(E,\mu)\), \(t_1\geq T\),
and \(t_2\geq0\). Choose \(f_n\in B(E)\) such that
\(f_n\to f\) in \(L^2(E,\mu)\). For each \(n\), the semigroup property on
\(B(E)\) gives
\[
        \psi_{t_1+t_2}[f_n]
        =
        \psi_{t_2}[\psi_{t_1}[f_n]] .
\]
Using \eqref{eq:B.1}, we may pass to the limit \(n\to\infty\) and obtain
\begin{equation}
        \psi_{t_1+t_2}[f](x)
        =
        \psi_{t_2}[\psi_{t_1}[f]](x),
        \qquad x\in E.
        \label{eq:B.2}
\end{equation}
Thus \(\psi\) satisfies an {\it eventual semigroup property} on \(L^2(E,\mu)\).

Next we claim that, for every
\(f\in L^2(E,\mu)\) and every \(t\geq T\),
\begin{equation}
        \lim_{h\downarrow 0}
        \|\psi_{t+h}[f]-\psi_t[f]\|_{L^2(E,\mu)}
        =
        0 .
        \label{eq:B.3}
\end{equation}
Indeed, by assumption 2 of the proposition, (M$1$), \ref{H4}, and the many-to-one formula, for every
\(g\in C(E,\mu)\),
\begin{equation}
        \lim_{h\to0}
        \|\psi_h[g]-g\|_{L^2(E,\mu)}
        =
        0 .
        \label{eq:B.4}
\end{equation}
Fix \(f\in L^2(E,\mu)\), \(t\geq T\), and choose \(g_n\in C(E,\mu)\) such that
\(g_n\to f\) in \(L^2(E,\mu)\). Using \eqref{eq:B.2}, for \(h\) sufficiently small we
have
\[
\|\psi_{t+h}[f]-\psi_t[f]\|_{L^2}
\leq
  \|\psi_t[\psi_h[g_n]-g_n]\|_{L^2} 
  +\|\psi_h[\psi_t[f-g_n]]\|_{L^2}
  +\|\psi_t[f-g_n]\|_{L^2}.
\]
The first term tends to zero as \(h\to0\), for fixed \(n\), by \eqref{eq:B.1} and
\eqref{eq:B.4}. The second and third terms tend to zero as \(n\to\infty\), uniformly for
\(h\) in a bounded neighbourhood of zero, by \eqref{eq:B.1} and the boundedness of the
Feynman-Kac potential. This proves \eqref{eq:B.3}.

We now turn to the spectral decomposition. Since \(\psi_T\) has an
\(L^2(E,\mu)\) density, \(\psi_T\) is a Hilbert-Schmidt operator on
\(L^2(E,\mu)\), and hence compact. By the spectral theory of compact operators,
the non-zero spectrum of \(\psi_T\) consists of isolated eigenvalues of finite
algebraic multiplicity, with possible accumulation point only at \(0\). Fix
\(K^\ast>0\). Let \(\rho_1,\ldots,\rho_m\) be the eigenvalues of \(\psi_T\)
with \(|\rho_i|\geq K^\ast\), and let \(P_i\) be the corresponding Riesz
projections. Set
\[
        A_i:=P_iL^2(E,\mu),
        \qquad
        P:=\sum_{i=1}^m P_i .
\]
Then \(A_i\) is finite dimensional and invariant under \(\psi_T\), and the
restriction of \(\psi_T\) to \(A_i\) has the form
\[
        \psi_T|_{A_i}=\rho_i I+D_i,
\]
where \(D_i\) is nilpotent. Moreover, on the residual subspace
\((I-P)L^2(E,\mu)\), the spectral radius of \(\psi_T\) is strictly smaller than
\(K^\ast\). Therefore, for every \(K'>0\) such that the residual spectral
radius is smaller than \(K'\), there exists \(C<\infty\) such that
\begin{equation}
\left\|\psi_{nT}[f]-
\sum_{i=1}^m (\rho_i I+D_i)^n P_i f
\right\|_{L^2(E,\mu)}
\leq C (K')^n \|f\|_{L^2(E,\mu)} \label{eq:B.5}
\end{equation}
for all \(n\geq0\) and \(f\in L^2(E,\mu)\).

We next record how this \(L^2\) decomposition gives the finite-dimensional
objects appearing in \ref{H1b}. Since \(\psi_T\) maps \(L^2(E,\mu)\) into \(B(E)\)
and the eigenspaces associated with \(\rho_i\neq0\) are finite dimensional,
their generalised eigenvectors admit bounded representatives. We choose such
representatives and denote a Jordan basis of \(A_i\) by
\[
        \varphi^{(k)}_{i,j},
        \qquad
        1\leq j\leq p_i,\quad 1\leq k\leq k_{i,j},
\]
where
\[
        D_i\varphi^{(k)}_{i,1}=0,
        \qquad
        D_i\varphi^{(k)}_{i,j}=\varphi^{(k)}_{i,j-1},
        \quad j\geq2.
\]
The corresponding coordinate maps are bounded on \(B(E)\). Indeed, they are
given by integration against the dual generalised eigenvectors in
\(L^2(E,\mu)\), and \(\mu(E)<\infty\) implies \(B(E)\subset L^2(E,\mu)\)
continuously. We denote these coordinate maps by
\(\tilde\varphi^{(k)}_{i,j}\), and define
\[
        \Phi_i[f](x)
        =
        \sum_{j=1}^{p_i}\sum_{k=1}^{k_{i,j}}
        \tilde\varphi^{(k)}_{i,j}[f]\varphi^{(k)}_{i,j}(x),
        \qquad f\in B(E).
\]

The eventual semigroup property \eqref{eq:B.2} implies that each \(A_i\) is invariant
under \(\psi_t\), for every \(t\geq0\). Moreover, by \eqref{eq:B.3}, the restriction
\((\psi_t|_{A_i})_{t\geq0}\) is a continuous finite-dimensional semigroup.
Thus there exist \(\lambda_i\in\mathbb C\), with \(e^{\lambda_iT}=\rho_i\), and
a nilpotent operator \(\mathcal N_i:A_i\to A_i\), such that
\begin{equation}
        \psi_t|_{A_i}
        =
        e^{(\lambda_i+\mathcal N_i)t},
        \qquad t\geq0.
        \label{eq:B.6}
\end{equation}
Writing \(\mathcal N\) for the direct sum of the nilpotent operators \(\mathcal N_i\), this gives
the spectral terms \(e^{(\lambda_i+\mathcal N)t}\Phi_i[f]\) appearing in \ref{H1b}.

It remains to prove the uniform asymptotic. Let \(K\in\mathbb R\) be as in the
statement of the proposition and choose \(K^\ast=e^{KT}\). Let
\(\rho_1,\ldots,\rho_m\) be the eigenvalues of \(\psi_T\) with
\(|\rho_i|\geq K^\ast\), as above. Choose \(K_0<K\) such that the residual
spectral radius of \(\psi_T\) is at most \(e^{K_0T}\). Then \eqref{eq:B.5} gives
\begin{equation}
\left\|
        \psi_{nT}[f]
        -
        \sum_{i=1}^m
        e^{(\lambda_i+\mathcal N)nT}\Phi_i[f]
\right\|_{L^2(E,\mu)}
\leq
        C e^{K_0nT}\|f\|_{L^2(E,\mu)} .
        \label{eq:B.7}
\end{equation}
Now write \(t=nT+s\), where \(n\in\mathbb N\) and \(s\in[T,2T)\). By \eqref{eq:B.2} and
\eqref{eq:B.6},
\[
\psi_t[f]
-
\sum_{i=1}^m e^{(\lambda_i+\mathcal N)t}\Phi_i[f]  =
\psi_s
\left[
        \psi_{nT}[f]
        -
        \sum_{i=1}^m e^{(\lambda_i+\mathcal N)nT}\Phi_i[f]
\right].
\]
Using \eqref{eq:B.1}, uniformly over
\(s\in[T,2T)\), we obtain
\[
\sup_{x\in E}
\left|
        \psi_t[f](x)
        -
        \sum_{i=1}^m e^{(\lambda_i+\mathcal N)t}\Phi_i[f](x)
\right| 
\leq
C
\left\|
        \psi_{nT}[f]
        -
        \sum_{i=1}^m e^{(\lambda_i+\mathcal N)nT}\Phi_i[f]
\right\|_{L^2(E,\mu)} .
\]
Taking the supremum over \(f\in B_1(E)\), and using
\(\|f\|_{L^2(E,\mu)}\leq \mu(E)^{1/2}\), \eqref{eq:B.7} gives
\[
\sup_{x\in E,\ f\in B_1(E)}
\left|
        \psi_t[f](x)
        -
        \sum_{i=1}^m e^{(\lambda_i+\mathcal N)t}\Phi_i[f](x)
\right|
\leq
        C e^{K_0t}.
\]
Since \(K_0<K\), it follows that
\[
\sup_{x\in E,\ f\in B_1(E)}
t e^{-Kt}
\left|
        \psi_t[f](x)
        -
        \sum_{i=1}^m e^{(\lambda_i+\mathcal N)t}\Phi_i[f](x)
\right|
\longrightarrow0,
\qquad t\to\infty .
\]
This proves \eqref{eq: h1b assum example}.

Finally, ordering the resulting eigenvalues so that
\[
        \operatorname{Re}\lambda_1
        \geq
        \operatorname{Re}\lambda_2
        \geq
        \cdots
        \geq
        \operatorname{Re}\lambda_m
        \geq K,
\]
and discarding the empty list if no eigenvalue lies above the threshold, gives
the statement of the proposition.

\section*{Acknowledgements}
This work was partially supported by the EPSRC grant MaThRad EP/W026899/1.
 Part of this work was carried out while CD was a PhD student at the EPSRC Centre for Doctoral Training in Statistical Applied Mathematics at Bath (SAMBa), under the project EP/S022945/1.
 The authors would like to thank Dr C\'ecile Mailler for her involvement with the initial stages of this project.

\bibliographystyle{plain}
\bibliography{bibbp.bib}

\end{document}